\documentclass[12pt,reqno]{amsart}

\usepackage{amsmath}
\usepackage{amsfonts}
\usepackage{amscd}
\usepackage{amssymb}
\usepackage{amsthm}
\usepackage{mathdots}
\usepackage{mathtools}
\usepackage{multicol}
\usepackage{inputenc}
\usepackage{caption}
\usepackage{subcaption}
\usepackage{float}

\numberwithin{equation}{section}

\usepackage{MnSymbol}
\usepackage[linktocpage=true]{hyperref}
\usepackage{mathrsfs}

\usepackage{tikz,tikz-cd, color}
\usepackage{adjustbox}

\usepackage{appendix}
\usetikzlibrary{matrix}
\usetikzlibrary{decorations.pathmorphing}

\tikzset{
  symbol/.style={
    draw=none,
    every to/.append style={
      edge node={node [sloped, allow upside down, auto=false]{$#1$}}}
  }
}

\usepackage{stmaryrd}
\usepackage{enumerate}
\usepackage{xcolor}

\usepackage[all]{xy}

\newtheorem{thm}{Theorem}[section]
\newtheorem{lem}[thm]{Lemma}
\newtheorem{prop}[thm]{Proposition}
\newtheorem{cor}[thm]{Corollary}
\newtheorem{conj}[thm]{Conjecture}

\theoremstyle{definition}
\newtheorem{defn}[thm]{Definition}
\newtheorem{exmp}[thm]{Example}

\theoremstyle{remark}
\newtheorem{rem}[thm]{Remark}

\newcommand{\oo}{\mathcal{O}}
\newcommand{\bp}{\mathbb{P}}

\newcommand{\Hom}{{\rm Hom}}
\newcommand{\lHom}{\mathcal{H}om}
\newcommand{\Ext}{{\rm Ext}}
\newcommand{\lExt}{\mathcal{E}xt}
\newcommand{\ext}{{\rm ext}}
\newcommand{\Quot}{{\rm Quot}}

\newcommand{\rk}{{\rm rk\,}}
\newcommand{\ch}{{\rm ch}}
\newcommand{\im}{{\rm im\,}}
\newcommand{\id}{{\rm id}}

\newcommand{\ev}{{\rm ev}}
\newcommand{\E}{\mathcal{E}}
\newcommand{\dE}{{\vphantom{\E}}^{\vee}\!\E}
\newcommand{\basefield}{{k}}
\newcommand{\hilbl}{^{[\ell]}}
\newcommand{\syml}{^{(\ell)}}

\newcommand{\rcone}{{\hat{R}}}

\DeclareMathOperator{\gr}{gr}
\DeclareMathOperator{\SD}{SD}

\DeclareMathOperator{\coker}{coker}
\DeclareMathOperator{\db}{D^b}

\DeclareMathOperator{\TOR}{\mathcal{T}or}
\DeclareMathOperator{\K}{K}
\DeclareMathOperator{\Pic}{Pic}
\DeclareMathOperator{\GL}{GL}
\DeclareMathOperator{\Char}{Char}
\DeclareMathOperator{\sstable}{ss}
\DeclareMathOperator{\stable}{s}

\author[]{Thomas Goller}
\address{The College of New Jersey, NJ, USA}
\email{gollert@tcnj.edu}
\author[]{Yinbang Lin}
\address{Tongji University, Shanghai, China}
\email{yinbang\textunderscore lin@tongji.edu.cn}

\title[Gaeta resolutions and strange duality]{Gaeta resolutions and strange duality\\ over rational surfaces}
\begin{document}
\begin{abstract}
    Over the projective plane and at most two-step blowups of Hirzebruch surfaces, where there are strong full exceptional sequences of line bundles, we obtain foundational results about Gaeta resolutions of coherent sheaves by these line bundles. Under appropriate conditions, we show the locus of semistable sheaves not admitting Gaeta resolutions has codimension at least 2. We then study Le Potier's strange duality conjecture. Over these surfaces, for two orthogonal numerical classes where one has rank one and the other has sufficiently positive first Chern class,  we show that the strange morphism is injective. The main step in the proof is to use Gaeta resolutions to show that certain relevant Quot schemes are finite and reduced, allowing them to be enumerated using the authors' previous paper.
\end{abstract}

\maketitle

\section{Introduction}

In the moduli theory of sheaves over complex algebraic surfaces, there is a famous conjecture by Le Potier \cite{LP05}, called the {\em  strange duality conjecture}, which relates the global sections of two determinant line bundles on certain pairs of
moduli spaces of sheaves. Known results over rational surfaces are mostly in the cases where one of the moduli spaces parametrizes pure dimension 1 sheaves, see e.g. \cite{Dan02,Abe10,Yua21}. On other surfaces, the conjecture requires different formulations, see e.g. \cite{BolMarOpr17}.
In an attempt to provide a unified treatment of the conjecture over rational surfaces, the first author, Bertram, and Johnson \cite{BerGolJoh16} proposed to use Grothendieck's {\em Quot schemes} \cite{Gro60}, following Marian and Oprea's ideas \cite{MarOpr07,MarOpr07b} over curves.
A key tool in the study of Quot schemes over $\mathbb{P}^2$ in \cite{BerGolJoh16} is {\em Gaeta resolutions} of coherent sheaves in terms of the {\em strong full exceptional sequence} of line bundles  $(\oo(-2),\oo(-1),\oo)$.

This leads us to the study of exceptional sequences and Gaeta resolutions over rational surfaces.
A strong full exceptional sequence, if it exists, completely captures the derived category in an explicit way, by theorems of Baer \cite{Bae88} and Bondal \cite{Bon89}.
While coherent sheaves can always be resolved by locally free sheaves by the Hilbert Syzygy Theorem, we can bring the resolution under better control if the locally free sheaves are taken from an exceptional sequence.
Resolutions built from such exceptional sequences, called Gaeta resolutions, have been applied toward a variety of problems in the study of sheaves on rational surfaces \cite{Dre86,LeP94,CosHui18weakBN,CosHui20,coskun_existence_2019}. Though the existence of strong full exceptional sequences of line bundles in general is an open question, the answer is affirmative over a rational surface that can be obtained from a Hirzebruch surface $\mathbb{F}_e$ by blowing up at most two sets of points
\cite{hille-perling-exceptional11}, which we call a \emph{two-step blowup of $\mathbb{F}_e$}.

Over $\mathbb{P}^2$ or a two-step blowup of $\mathbb{F}_e$,
we choose a particular strong full exceptional sequence, determine when a sheaf admits Gaeta resolutions, and study general properties of such sheaves.
We then apply Gaeta resolutions to the study of strange duality, proving the injectivity of the strange morphism in some cases.
One of the key points in the proof is to show that relevant Quot schemes are finite and reduced, which we accomplish using Gaeta resolutions.
A parallel statement over $\mathbb{P}^2$ was proved in \cite{BerGolJoh16}.
Another crucial point is to enumerate the length of the finite Quot scheme, which was settled in our previous paper \cite{GolLin22} via the study of the moduli space of limit stable pairs \cite{Lin18}.

\medskip

We now set up the study. Let $S$ be a smooth projective algebraic surface over an algebraically closed field $k$ of characteristic $0$, with a strong full exceptional sequence of line bundles
\[(\E_1, \E_2,\dots, \E_n).\]
Given a coherent sheaf $F$ on $S$, we would like to find a resolution of $F$ of the form\footnote{To avoid clumsy notation, we drop the direct sum symbol from the exponents. If we want to denote tensor products, we will use $\otimes$.}
\begin{equation}
      0 \to \E_1^{a_1} \oplus \cdots \oplus \E_d^{a_d} \to \E_{d+1}^{a_{d+1}} \oplus \cdots \oplus \E_{n}^{a_{n}} \to F \to 0,
\end{equation}
which is called a {Gaeta resolution}.
If a Gaeta resolution exists, the exponents $a_1,\dots,a_n$ are uniquely determined by the numerical class of $F$.

One of our main technical results is a criterion to determine when a Gaeta resolution exists.
We state here the criterion for a two-step blowup of $\mathbb{F}_e$,
%a Hirzebruch surface, %$\mathbb{F}_e$
where we know there are strong full exceptional sequences of line bundles, see \S~\ref{ss:ex-sfes}.
%\todo{Y: reference to examples.}
Let $S_0$ denote the set of blown-up points in the first step, with corresponding exceptional divisors $E_i$ for $i \in S_0$. The second set of blown-up points, which lie on these exceptional divisors, is denoted $S_1$, and the corresponding exceptional divisors are $E_j$ for $j \in S_1$. Using the exceptional sequence \S~\ref{ss:ex-sfes}(c), we get the following result.

\begin{prop}[Proposition~\ref{prop:special-GTR-criterion}(b)]\label{prop:crit-blowup}
On a two-step blowup of $\mathbb{F}_e$, a torsion-free sheaf $F$ has a Gaeta resolution if and only if
\begin{enumerate}[(i)]
    \item $H^p(F)$ vanishes for $p \ne 0$;
    \item $H^p(F(D))$ vanishes for $p \ne 1$ and $D$ the divisors $-A+\sum_{i \in S_0} E_i$, $-B+\sum_{i \in S_0} E_i$, and $-A-B+\sum_{i \in S_0} E_i + \sum_{j \in S_0} E_j$; and
    % \item $H^p(F(-A-B+\sum_{i \in S_0} E_i + \sum_{j \in S_0} E_j))$ vanish for $p \ne 1$; and
    \item $H^1(F(E_i))=0$ and $H^1(F(E_j)) = 0$ for all $i \in S_0$ and $j \in S_1$.
\end{enumerate}
\end{prop}
Section~\ref{sect:gtr} contains a similar statement for $\bp^2$, which is known.

We will apply Gaeta resolutions towards moduli problems, using the connections to prioritary sheaves.
We impose a mild technical condition on $X$ a two-step blowup of $\mathbb{F}_e$, which we call {\em admissibility} (Definition~\ref{defn:adm-blowup}). We use $A$ to denote both the class of the fibers of the ruling $\mathbb{F}_e\to \mathbb{P}^1$ and the pullback of this class to $X$.
We denote a numerical class $f$ in the Grothendieck group $\K(X)$ by
\[f=(r,L,\chi),\] where $r$ is the rank, $L$ is the first Chern class (equivalent to the determinant line bundle), and $\chi$ is the Euler characteristic.
We have

\begin{prop}\label{prop:unirat-large-L}
 Let $X$ be an admissible blowup of $\mathbb{F}_e$ and $H$ be a polarization such that $H\cdot(K_X+A)<0$.
 For a numerical class $f\in \K(X)$ with fixed rank $r>0$ and fixed $\chi\geqslant 0$, suppose the first Chern class $L$ is sufficiently positive. Then a general $H$-semistable sheaf of class $f$ admits Gaeta resolutions.
\end{prop}

We mostly use Gieseker stability and use either ``$H$-semistable'' or ``semistable''. If we want to use slope stability, we will specify.
For the precise meaning of being sufficiently positive in Proposition \ref{prop:unirat-large-L}, see the conditions in Proposition~\ref{prop:chern-has-gr}(a) as well as Proposition \ref{prop:cokernel-properties}(b.ii).
There is another statement, Proposition~\ref{prop:unirat-large-c2}, in which the rank and first Chern class are fixed,
which asserts that if the discriminant is sufficiently large then general semistable sheaves admit Gaeta resolution up to a twist by a line bundle. In each case, we can immediately deduce that the moduli space $M(f)$ is unirational, which is known \cite{Bal87}.

By imposing stronger conditions on the numerical class $f$ and the polarization $H$, we prove a refinement of Proposition~\ref{prop:unirat-large-L}: %that the complement of the locus of semistable sheaves admitting Gaeta resolutions has codimension $\geqslant 2$.
% \todo{Y: We should assume $r\geqslant 2$.\\Y: Changed.}
\begin{thm}\label{thm:moduli-nongr-codim2}
Let $X$ be an admissible blowup. Assume the class $f$ is of rank $r\geqslant 2$, admits Gaeta resolutions in which the exponents $a_i$ are strictly positive and satisfy (\ref{eq:conditions-exponents}, \ref{eq:condition-alpha4}), that the polarization $H$ is general and satifies (\ref{eq:H-pos-lin-combo}, \ref{eq:conditions-on-H}), and that the discriminant of $f$ is sufficiently large in the sense of (\ref{eq:discriminant-bound}). % $a_i>0$, $i=1,\dots,n$.
    Then the closed subset $Z\subset M(f)$ of S-equivalence classes of semistable sheaves whose Jordan-H\"older gradings do not admit Gaeta resolutions has codimension $\geqslant 2$ in $M(f)$.
\end{thm}
If $X=\mathbb{F}_e$, the conditions (\ref{eq:conditions-exponents}, \ref{eq:H-pos-lin-combo}, \ref{eq:conditions-on-H}) are vacuous.
Using Proposition~\ref{prop:unirat-large-L} and Theorem~\ref{thm:moduli-nongr-codim2}, we deduce

\begin{cor}\label{cor:gen-stable-prop}%Suppose $r\geqslant 2$.
Assuming the conditions from Proposition \ref{prop:unirat-large-L} and appropriate conditions from Proposition \ref{prop:cokernel-properties}(b), a general sheaf in $M(f)$ is torsion-free, is locally-free if $r \geqslant 2$, satisfies the cohomological vanishing conditions in Proposition \ref{prop:special-GTR-criterion} (b.i-iii), and is globally generated if $\chi \geqslant r+2$. Assuming the stronger conditions from Theorem \ref{thm:moduli-nongr-codim2}, the same properties hold away from a locus of codimension $\geqslant 2$ in $M(f)$ (with locally free replaced by torsion-free).
\end{cor}

For the proof of Theorem~\ref{thm:moduli-nongr-codim2}, we will need the following statement which is of general interest.
\begin{thm}\label{thm:yoshioka}
Suppose $S$ is a rational surface other than $\mathbb{P}^2$ and $S\to \mathbb{P}^1$ is a morphism where a general fiber $D$ is isomorphic to $\mathbb{P}^1$. Let $H$ be a general ample divisor such that $H\cdot (K_S+2D)<0$. Assume there is a slope stable vector bundle with rank $r\geqslant 2$, first Chern class $c_1$, and second Chern class $c_2-1$ over $S$. Then
\[\Pic M(r,c_1,c_2)\cong \mathbb{Z}\oplus \Pic S.\]
\end{thm}

%\todo{T: The notation $(r,c_1,c_2)$ differs from our usual notation for numerical classes. There's probably no good way to avoid this. Do we need a comment to clarify?}

%\medskip

We will apply Gaeta resolutions towards the study of the strange duality conjecture.
Let $S$ be $\mathbb{P}^2$ or $X$ an admissible blowup of $\mathbb{F}_e$.
In the Grothendieck group $\K(S)$, let \[\sigma = (r,L,r\ell) \mbox{ for } r \geqslant 2 \qquad\mbox{and}\qquad \rho = (1,0,1-\ell) \mbox{ for } \ell \geqslant 1 \] be two numerical classes.
Notice that they are orthogonal: $\chi(\sigma\cdot \rho)=0$.
Let $H$ denote the hyperplane class on $\bp^2$ or a polarization satisfying (\ref{eq:H-pos-lin-combo}, \ref{eq:conditions-on-H}) on $X$. %, which in particular implies $ H \cdot (K_X+A) < 0$.
% \todo{T: Add references to assumptions on $H$ here.\\Y:Added. \\ T: Do we need to mention $H \cdot (K_X+A) < 0$?\\Y: I guess there is no need anymore.}
Let $M(\sigma)$ and $M(\rho)$ denote the moduli spaces of $H$-semistable sheaves, where $M(\rho)$ is isomorphic to the Hilbert scheme $S^{[\ell]}$ of points.
Let $\mathscr{Z}\subset S^{[\ell]}\times S$ be the universal subscheme and $I_{\mathscr{Z}}$ its ideal sheaf. For a coherent sheaf $W$ of class $\sigma$, consider the determinant line bundle \begin{equation*}
    \Theta_\sigma:=\det\left(p_{!}\left(I_{\mathscr{Z}}\stackrel{L}{\otimes} {q}^*W\right)\right)^{*}
\end{equation*}
on $S^{[\ell]}$, where ${p}$ and ${q}$ are the projections from $S^{[\ell]}\times S$ to the first and second factors, respectively.
There is also a similar line bundle $\Theta_\rho$ on $M(\sigma)$.
On $M(\sigma)\times M(\rho)$, the line bundle $\Theta_{\sigma,\rho}\cong \Theta_\rho\boxtimes \Theta_\sigma$ has a canonical section, which induces the {\em strange morphism}
\[\operatorname{SD}_{\sigma,\rho}\colon H^0(S^{[\ell]}, \Theta_\sigma)^*\to H^0(M(\sigma),\Theta_\rho).\]
The strange duality conjecture says that $\operatorname{SD}_{\sigma,\rho}$ is an isomorphism. For a more general setup, see \S~\ref{subsect:strange-mor}.

In this context, we prove the following result in support of the strange duality conjecture. In the statement, for $V$ a vector bundle on $S$, we write
\begin{equation}\label{eq:taut-bundle}
    V^{[\ell]}={p}_{*}(\oo_{\mathscr{Z}}\otimes {q}^* V).
\end{equation}

\begin{thm}\label{thm:sd-injective} Let $S$ be $\bp^2$ or $X$ an admissible blowup of $\mathbb{F}_e$, and $\sigma$, $\rho$, and $H$ as above.
If $L$ is sufficiently positive, then:
\begin{enumerate}[(a)]
\item The rank of the strange morphism is bounded below by $\int_{S^{[\ell]}} c_{2\ell}({V}^{[\ell]})$, for $V$ a vector bundle with numerical class $\sigma + \rho$;
\item %If Conjecture \ref{numbers-match} holds, for instance in the case $\ell \leqslant 11$, then
The strange morphism $\mathrm{SD}_{\sigma,\rho}$ is injective.
\end{enumerate}
\end{thm}

We sketch the proof.
Let $V$ be a vector bundle of class $\sigma + \rho$ that admits a general Gaeta resolution and consider quotients of $V^*$ of class $\rho$. Then the Quot scheme has expected dimension $0$. If it is finite and reduced, a simple argument shows that its length provides a lower bound for $\operatorname{SD}_{\sigma,\rho}$.
According to \cite{GolLin22}, in this case its length is $\int_{S^{[\ell]}}c_{2\ell}(V^{[\ell]})$. On the other hand, Theorem~\ref{numbers-match} relates this top Chern class to $\chi(S^{[\ell]}, \Theta_\sigma)$, and the determinant line bundle $\Theta_{\sigma}$ has no higher cohomology when $L$ is sufficiently positive, which finishes the proof. Thus, the crucial point is to establish that the Quot scheme is finite and reduced, which we prove by considering the relative Quot scheme over the space of Gaeta resolutions and calculating the dimension of the relative Quot scheme. The following theorem summarizes the key results related to the Quot scheme.

\begin{thm}\label{thm:finite-quot-scheme}
Let $S$, $\sigma$, $\rho$, and $H$ be as in Theorem~\ref{thm:sd-injective}, and $V$ be a vector bundle of class $\sigma+\rho$ that admits a general Gaeta resolution. If $L$ is sufficiently positive, then:
\begin{enumerate}[(a)]
    \item The Quot scheme $\Quot(V^*,\rho)$ parametrizing quotient sheaves of $V^*$ with numerical class $\rho$ is finite and reduced;
    \item For every point $[V^* \twoheadrightarrow F]$ of $\Quot(V^*,\rho)$, $F$ is an ideal sheaf $I_Z$ for general $Z \in S^{[\ell]}$ and the kernel is semistable;
    %\item If $M(\sigma)$ is nonempty, then the kernel of each quotient ;
    \item The length of $ \Quot(V^*,\rho)$ is $\int_{S^{[\ell]}} c_{2\ell}({V}^{[\ell]})$.
\end{enumerate}
\end{thm}
In \cite{BerGolJoh16}, parts (a) and (b) were proved over $\mathbb{P}^2$ and calculations were made that informed Johnson's expectation that the counting formula (c) should be true for del Pezzo surfaces \cite{Joh18}.
The formula (c) was proved in \cite{GolLin22} by the authors of the current paper, for a general smooth regular projective surface, assuming that the Quot scheme is finite and reduced.

The positivity conditions on $L$ in these theorems, which are stronger than for Proposition \ref{prop:unirat-large-L}, are summarized in the appendix. Theorem \ref{thm:finite-quot-scheme} requires (\ref{eq:first-three-conditions}, \ref{eq:sd-discriminant-bound}), and Theorem \ref{thm:sd-injective} requires (\ref{eq:fifth-condition}) as well.
% \todo{Y: li-qin \\
% T: Theorem 1.7 (c) assumes $M(\sigma)$ is non-empty, and we need this for Theorem 1.6. Shall we make the stronger assumptions on $H$ in these theorems to guarantee non-emptiness? \\
% T: I've replaced the reference to li-qin by (\ref{eq:sd-discriminant-bound}) (the discriminant bound).}

\medskip

We organize the paper as follows.
In \S~\ref{sect:hirzebruch}, we review basic facts about divisors and line bundles on Hirzebruch surfaces and their blowups.
In \S~\ref{sect:exc-sequence}, we review exceptional sequences in the derived category.
In \S~\ref{sect:gtr}, we obtain criteria for the existence of Gaeta resolutions, including %our main technical tool
Proposition~\ref{prop:crit-blowup}, and classify the numerical classes of sheaves admitting Gaeta resolutions.
In \S~\ref{sect:prop-gr}, we prove some general properties of such sheaves and relate them to prioritary sheaves.
In \S~\ref{sec:stability}, we discuss connections to semistable sheaves and prove Theorem~\ref{thm:moduli-nongr-codim2}.
In \S~\ref{sect:sd}, we set up the strange morphism and prove Theorem~\ref{thm:sd-injective}.
In \S~\ref{sect:finite-quot}, we prove Theorem~\ref{thm:finite-quot-scheme}.
Finally, the appendix contains a summary of the positivity conditions on $L$ required in the proofs of Theorems \ref{thm:sd-injective} and \ref{thm:finite-quot-scheme}.
% \ref{thm:sd-injective} and \ref{thm:finite-quot-scheme}.

{\em Acknowledgment. }YL would like to thank Alina Marian and Dragos Oprea for helpful correspondences. TG would like to thank Lothar G\"{o}ttsche for providing updates on his work with Anton Mellit.
% question and encouragement, which initiated the study. YL also would like to thank  for pointing out the recent work by G\"ottsche-Mellit before it was posted on arXiv.
YL is supported by grants from the Fundamental Research Funds for the Central Universities and Applied Basic Research Programs of Science and Technology Commission Foundation
of Shanghai Municipality.

\section{Hirzebruch surfaces and blowups}\label{sect:hirzebruch}

We review some basic results on divisors and cohomology of line bundles on Hirzebruch surfaces and their blow-ups. Along the way, we introduce two technical assumptions (\ref{cond:avoid-b}) and (\ref{cond:avoid-fiber-dir}). The first is not a restriction, while the second is a mild condition.

\subsection{Divisors on blowups of Hirzebruch surfaces}
Let $\mathbb{F}_e$ denote the Hirzebruch surface $\bp(\oo_{\bp^1} \oplus \oo_{\bp^1}(e))$ with $e\geqslant 0$, which is a rational surface that is ruled over $\bp^1$. Note that $\mathbb{F}_0 \cong \bp^1 \times \bp^1$. Letting $A$ denote the class of a fiber and $B$ the 0-section with self-intersection $-e$, $A$ and $B$ generate the effective cone of $\mathbb{F}_e$, and $A^2=0$, $A \cdot B = 1$,  $B^2 = -e$.
The canonical divisor is
$K_{\mathbb{F}_e}=-(e+2)A-2B$.
The divisor classes $eA+B$ and $A$ generate the nef cone. To simplify notation, let \[C=eA+B.\] The linear system $|A|$ induces the morphism to $\bp^1$ giving the ruling, while $|C|$ induces a morphism $\mathbb{F}_e \to \bp^{e+1}$; if $e=0$, this is the other ruling of $\bp^1 \times \bp^1$, while if $e > 0$, this contracts $B$ to a point and maps the fibers of the ruling to distinct lines through that point in $\bp^{e+1}$.
In particular, suppose $x,y$ are distinct points on $\mathbb{F}_e$, where we allow $y$ to be infinitely near to $x$. Then $|A|$ separates $x,y$ unless $x,y$ are distinct points on the same fiber or $y$ corresponds to the tangent direction along the fiber at $x$. Similarly, $|C|$ separates $x,y$ unless $x,y$ are distinct points on $B$ or $x \in B$ and $y$ is the tangent direction along $B$ at $x$.

\begin{rem}\label{rem:blowups}
The blowup of $\bp^2$ at any point is isomorphic to $\mathbb{F}_1$. The blowup of $\bp^1 \times \bp^1$ at any point is isomorphic to the blowup of $\mathbb{F}_1$ at a point not on $B$.
For $e>0$, the blowup of $\mathbb{F}_e$ at a point on $B$ is isomorphic to the blowup of $\mathbb{F}_{e+1}$ at a point not on $B$. Thus, when considering the surfaces that arise from blowing up $\bp^2$ or Hirzebruch surfaces, it suffices to consider blowups of $\mathbb{F}_e$ for $e>0$ where the blown-up points are not on $B$ \cite[p.519]{GriHar78}. So, quite often, in the case $e>0$, we impose the condition that
\begin{equation}\label{cond:avoid-b}
    \text{the blowup avoids $B$}
\end{equation}
in the sense that none of the blown-up points $p_1,\dots,p_s$ is on $B$.
\end{rem}

Let $X$ be obtained from a sequence of blowups
\[
    X=X_t \to X_{t-1} \to \cdots \to X_1 \to X_0 = \mathbb{F}_e,
\]
where $b_i \colon X_i \to X_{i-1}$ is the $i$th blowup at a point $p_i\in X_{i-1}$. Assume that the indices $\{1,\dots,t\}$ can be partitioned into two sets \[S_0=\{1,\dots,s\} \quad \mbox{and} \quad S_1 = \{s+1,\dots,t\}\] such that the $b_i$ within each of the sets $\{b_1,\dots,b_s\}$ and $\{b_{s+1},\dots,b_t\}$ commute. In other words, $X$ can be obtained from $\mathbb{F}_e$ by up to two blowups, each possibly at multiple points.
By Remark \ref{rem:blowups}, it suffices to consider the case $e>0$ and that the blowup avoids $B$.

We define a partial ordering on the set $\{ p_1,\dots,p_t \}$ by $p_j \succ p_i$ if $p_j$ is on the exceptional divisor of $b_i$, and say that the height of $p_i$ is 0 if $p_i$ is minimal with respect to $\succ$, while otherwise $p_i$ has height 1. We can consider points of height 0 as lying on $\mathbb{F}_e$, while if $p_j \succ p_i$ then $p_j$ is infinitely near to $p_i$ and can be viewed as a tangent direction at $p_i$ on $\mathbb{F}_e$. For simplicity, we choose the partition $S_0$ and $S_1$ such that
\begin{align*}
    i \in S_0 &\quad \mbox{ iff } p_i \mbox{ has height 0, while }\\
    j \in S_1 &\quad \mbox{ iff } p_j \mbox{ has height }1.
\end{align*} Thus, we think of $X$ as being obtained from $\mathbb{F}_e$ by first blowing up a collection of points $\{p_1,\dots,p_s \}$ on $\mathbb{F}_e$ and then blowing up a collection of points $\{ p_{s+1},\dots,p_t \}$ on the exceptional divisors of the first blowup.

Let $E_i$ denote the total transform in $X$ of the exceptional divisor of $b_i$. We abuse notation by writing $A,B$ for the pullbacks of the divisors $A,B$ on $\mathbb{F}_e$. The Picard group of $X$ is generated by $A,B,E_1,\dots,E_t$, with the following intersections:
\[
    A^2=0, \quad A \cdot B=1, \quad B^2=-e, \quad A \cdot E_i=0, \quad B \cdot E_i=0, \quad E_i \cdot E_j=-\delta_{i,j}.
\]
Here, $\delta$ is the Kronecker delta function. The canonical divisor is
\[
    K_X = -(e+2)A-2B + \sum_{i=1}^t E_i.
\]
Let $\tilde{E}_i$ denote the strict transform of $E_i$. If $p_j$ has height 1, then $\tilde{E}_j = E_j$, while if $p_i$ has height 0, then $\tilde{E}_i = E_i - \sum_{j \colon p_j \succ p_i} E_j$. Note that $\tilde{E}_i^2 = -1 - \#\{ j \colon p_j \succ p_i \}$.

The following classification of base loci of certain linear systems on $X$ will be useful. When describing the linear systems below, we write a divisor in parentheses, as in $(D)$, to indicate that it is a fixed part of the linear system.

\begin{lem}\label{lem:LS_base_locus}
Suppose $D$ on $X$ is the pullback of an effective divisor on $\mathbb{F}_e$. Then if $p_j \succ p_i$,
\[
    |D-E_j| = |D-E_i| + (E_i-E_j).
\]
\end{lem}

\begin{proof} Tensoring the short exact sequence
\[
    0 \to \oo(-E_i+E_j) \to \oo \to \oo_{(E_i-E_j)} \to 0
\]
by $\oo(D-E_j)$ and taking cohomology, we get an exact sequence
\[
    0 \to H^0(\oo(D-E_i)) \to H^0(\oo(D-E_j)) \to H^0(\oo(D-E_j)|_{(E_i-E_j)}).
\]
As $\tilde{E}_i \cdot (D-E_j) = -\tilde{E}_i \cdot E_j = -1$, we see that $H^0(\oo(D-E_j)|_{(E_i-E_j)}) = 0$ as $(E_i-E_j)$ is a connected (possibly reducible) curve and every section of this line bundle must be 0 on the component $\tilde{E}_i$.
\end{proof}

\begin{rem}\label{rem:basepoints} If $D$ is a divisor on $\mathbb{F}_e$ and $p_i$ is a point of height 0, then the curves in the linear series $|D-E_i|$ on $X$ are in bijection with curves in $|D|$ on $\mathbb{F}_e$ that contain $p_i$. Similarly, if $p_j \succ p_i$, then curves in $|D-E_i-E_j|$ on $X$ are in bijection with curves in $|D|$ on $\mathbb{F}_e$ that contain $p_i$ and have tangent direction $p_j$ at $p_i$. The curves in the linear system on $X$ are obtained as pullbacks of the corresponding curves on $\mathbb{F}_e$, with one copy of the appropriate exceptional divisors removed.
\end{rem}

The linear systems on $X$ in the following examples will play an important role. We assume that $e>0$ and that the blowup avoids $B$.
\begin{exmp}\label{ex:LS_fiber}
If $j \in S_1$, let $i \in S_0$ denote the index such that $p_j \succ p_i$. Then
\[
    |A-E_j|=(A-E_i) + (E_i-E_j).
\]
Moreover, if $p_j$ is the tangent direction along the fiber $A$ at $p_i$, then
\[
    |A-E_j|=(A-E_i-E_j) + (E_i)
\]
The curve $(A-E_i)$ can be obtained by considering $b_i$ to be the first blown-up point, taking the strict transform of the unique fiber $A$ containing $p_i$ under $b_i$, and then taking the pullback of that strict transform under the remaining blowups, which may be reducible if other $p_i$ lie on that fiber or are infinitely near to points on that fiber. In particular, if $p_j$ is the tangent direction along the fiber $A$ at $p_j$, then $p_j$ lies on $(A-E_i)$, and taking the strict transform with respect to $b_{j}$ yields the curve $(A-E_i-E_j)$.
\end{exmp}

\begin{exmp}\label{ex:LS_section}
For $j \in S_1$, let $i \in S_0$ denote the index such that $p_j \succ p_i$. Then
\[
    |C-E_j| = |C-E_i| + (E_i-E_j),
\]
where $|C-E_i|$ is basepoint-free. This follows from the fact that $|C|$ separates points on $\mathbb{F}_e$ (including infinitely near points) as long as the points are not contained on $B$, and by assumption the blown-up points are not on $B$.
\end{exmp}

\begin{exmp}\label{ex:LS_ample} For $i \in S_0$ and $j \in S_1$, the base locus of $|C+A-E_i-E_j|$ can be described as follows:
\begin{enumerate}[(a)]
    \item If $p_j \succ p_i$, then $|C+A-E_i-E_j|$ is basepoint-free unless $p_j$ is the tangent direction along the fiber $A$ containing $p_i$, in which case
    \[
        |C+A-E_i-E_j| = |C| + (A-E_i-E_j).
    \]
    \item If $p_j \not \succ p_i$, let $i' \in S_0$ denote the index such that $p_j \succ p_{i'}$. Then
    \[
        |C+A-E_i-E_j|=|C+A-E_i-E_{i'}| + (E_{i'}-E_j),
    \]
    and $|C+A-E_i-E_{i'}|$ is basepoint-free unless $p_i$ and $p_{i'}$ lie on the same fiber $A$, in which case
    \[
        |C+A-E_i-E_j|=|C| + (A-E_i-E_{i'}) + (E_{i'}-E_j).
    \]
\end{enumerate}

We explain the two parts. For (b), the linear system contains the union of $(A-E_i)$ and a curve in $|C-E_{i'}|$. As the latter is basepoint-free, the only possible basepoints are on $(A-E_i)$. By the same argument with the roles of $i$ and $i'$ reversed, we see that the linear system is basepoint-free unless $p_i$ and $p_{i'}$ lie on the same fiber.

For (a), we note that if the linear system has a basepoint $p$, which we may assume is a point of $\mathbb{F}_e$, then the linear system of curves in $|C+A|$ that contain $p_i$ and $p$ must have $p_j$ as an (infinitely near) basepoint, which, by a similar argument as the one for (b), implies that $p_i$ and $p$ lie on the same fiber and $p_j$ is the tangent direction along that fiber.
\end{exmp}

Some of the calculations in later sections are simplified if the linear systems $|C+A-E_i-E_j|$ for $p_j \succ p_i$ are basepoint-free. For this purpose, we will often assume that
\begin{equation}\label{cond:avoid-fiber-dir}\mbox{the blowup avoids fiber directions}\end{equation}
in the sense that $p_{s+1},\dots,p_t$ are distinct from the point where the strict transform of the fiber $A$ containing $p_i$ meets the exceptional divisor of $b_i$, for all $i \in S_0$.

We summarize the assumptions on $X$ in the following definition:
\begin{defn}\label{defn:adm-blowup} The rational surface $X$ is an {\em admissible blowup} of $\mathbb{F}_e$ if it is an at most two-step blowup and the following conditions hold:
\begin{itemize}
    \item if $e=0$, then $S_0=S_1 = \emptyset$;
    \item if $e>0$, then the blowup avoids $B$ and avoids fiber directions.
\end{itemize}
% \todo{Y: I think it's better to combine the last two items. \\ T: Done.}
\end{defn}
We emphasize that $S_0$ or $S_1$ can be empty. In particular, the definition includes $\mathbb{F}_e$.
Then, by the above discussion and particularly Example \ref{ex:LS_ample}, we have shown that if $X$ is admissible, then every divisor in the set
\begin{equation} \label{eq:bpf-divisors}
    \mathcal{D} = \{A,C\} \cup \{ C-E_i \mid i \in S_0 \} \cup \{C+A-E_i-E_j \mid p_j \succ p_i \}
\end{equation}
is basepoint-free. This leads to the following result.

\begin{prop}\label{prop:very-ample} Suppose $X$ is an admissible blowup of $\mathbb{F}_e$. Suppose $L$ is the line bundle associated to a positive integral linear combination of all divisors in the set
\[
    \mathcal{D} = \{A,C\} \cup \{ C-E_i \mid i \in S_0 \} \cup \{C+A-E_i-E_j \mid p_j \succ p_i \}.
\]
Then $L$ is very ample.
\end{prop}

%\todo{T: Is the following proof too wordy? Should we reduce it to just a basic sketch?}

\begin{proof} The linear system associated to $L$ contains unions of divisors in the linear systems associated to the divisors in $\mathcal{D}$, so since these divisors are all basepoint-free, it suffices to show that the divisors in $D$ collectively separate points and tangents on $X$ in the following sense:
\begin{enumerate}[1)]
    \item For any distinct points $q_1,q_2 \in X$, there is a divisor $D$ linearly equivalent to a divisor in $\mathcal{D}$ such that $q_1 \in D$ and $q_2 \notin D$;
    \item For any $q \in X$, there are divisors $D_1,D_2$ that contain $q$, are each linearly equivalent to a divisor in $\mathcal{D}$, and whose images in $\mathfrak{m}_q/\mathfrak{m}_q^2$ are linearly independent.
\end{enumerate}

For 1), $|C|$ can be used to separate points on the complement of $B \cup \bigcup_{i \in S_0} E_i$ since it is very ample there, $|A|$ can separate two points on $B$, $|C-E_i|$ can separate two points on $\tilde{E}_i$, and $|C+A-E_i-E_j|$ separates any two points on $E_j$. Separating points on different exceptional curves or a point on the exceptional locus from a point on the complement is similarly easy.

For 2), as $C$ is very ample on the complement of $B \cup \bigcup_{i \in S_0} E_i$, it suffices to consider the cases $q \in D$, where $D$ is $B$, $\tilde{E}_i$ for $i \in S_0$, or $E_j$ for $j \in S_1$. In each case, it suffices to choose $D_1$ transversal to $D$ at $q$ and $D_2$ which is a union of $D$ with another divisor that does not contain $q$. These divisors can be chosen to be general in the following linear subsystems:
\begin{table}[H]
    \centering
    \begin{tabular}{c|c|c}
        & $D_1$ & $D_2$ \\ \hline
        $q \in B$ & $|A|$ & $(B) + |eA| \subset |C|$  \\
        $q \in E_i \setminus \bigcup_{j \colon p_j \succ p_i} E_j$ & $|C-E_i|$ & $(E_i) + |C-E_i| \subset |C|$ \\
        $q \in E_j$ & $|C+A-E_i-E_j|$ & $(E_j) + |C-E_i-E_j| \subset |C-E_i|$
    \end{tabular}
    %\caption{Caption}
    %\label{tab:my_label}
\end{table}
\noindent This completes the proof.
\end{proof}

As very ample is equivalent to 1-very ample and $m$-very ampleness is additive under tensor products \cite{HTT05}, we immediately see that if $L$ is the line bundle associated to a positive integral linear combination of all divisors in $\mathcal{D}$ in which the weight of each divisor is $\geqslant m$, then $L$ is $m$-very ample.

\subsection{Cohomology of line bundles}\label{subsect:coh-line-bdl}

First, we summarize how to calculate the
cohomology groups of line bundles on the Hirzebruch surface $\mathbb{F}_e$, following \cite{CosHui18weakBN}.
By Hirzebruch-Riemann-Roch,
\[\chi(\oo(aA+bB))=(a+1)(b+1)-\frac{1}{2}eb(b+1).\]
Since the effective cone of $\mathbb{F}_e$ is generated by $A$ and $B$,
\[H^0(\oo(aA+bB))\not=0 \quad \mbox{if and only if}\quad a,b\geqslant 0.\]
Then Serre duality implies that
\[H^2(\oo(aA+bB))\not=0\quad \mbox{if and only if}\quad a\leqslant -(e+2)\mbox{ and }b\leqslant -2.\]
It suffices to assume that $b\geqslant -1$, as other cases can then be obtained via Serre duality. In this case, as $h^2$ vanishes and the Euler characteristic is known, it suffices to calculate $h^0$, which can be done as follows:
\begin{enumerate}[(a)]
    \item $h^i(\oo(aA-B))=0$, for all $i$ and $a$.
    \item $h^0(\oo(aA))=a+1$ if $a\geqslant -1$ and $0$ otherwise.
    \item Let $b\geqslant 1$. If $a \geqslant be-1$, then
    \[h^0(\oo(aA+bB))=\chi(\oo(aA+bB)), \]
    %and $a\geqslant 0$,
    while if $a \leqslant be-2$, then the equality
    \[h^0(\oo(aA+bB)=h^0(\oo(aA+(b-1)B))\]
    allows $h^0$ to be determined by induction on $b$.
\end{enumerate}

In particular, we deduce the following:
\begin{lem}\label{lem:lb-van-higher-cohom} Let $L = \oo_{\mathbb{F}_e}(aA + bB)$. Then
\begin{enumerate}[(a)]
    \item $H^2(L)=0$ if and only if $b \geqslant -1$ or $a \geqslant -1-e$;
    \item $H^1(L) = 0$ if $b=-1$, if $b=0$ and $a \geqslant -1$, or if $b \geqslant 1$ and $a \geqslant be-1$.
\end{enumerate}
\end{lem}

In order to use these calculations on $\mathbb{F}_e$ to obtain information about the cohomology of line bundles on $X$ a two-step blowup of $\mathbb{F}_e$, we use the following general result.

\begin{lem}\label{lem:line-bundle-pullback}Let $\pi\colon \tilde{Y}\to Y$ be a blowup of a smooth projective surface at distinct (possibly infinitely near) points.
For a line bundle $L$ on $Y$, $H^i(\tilde{Y},\pi^*L)\cong H^i(Y, L)$.
\end{lem}
\begin{proof}
According to \cite[V. Proposition 3.4]{Har77}, $\pi_*\oo_{\tilde{Y}}\cong \oo_Y$ and $R^i\pi_*\oo_{\tilde{Y}}=0$ for $i>0$. Then $\pi_*\pi^*L\cong L$ and $R^i\pi_*(\pi^*L)=0$ for $i>0$ by the projection formula.
The spectral sequence $H^i(Y,R^j\pi_*(\pi^*L))\Rightarrow H^{i+j}(\tilde{Y},\pi^*L)$ gives the result.
\end{proof}

Then, letting $X$ denote a two-step blowup of $\mathbb{F}_e$, we have:
\begin{lem}\label{lem:reduce-higher-cohom-blowup} Let $L$ be a line bundle on $X$ such that $L|_{E_i} \cong \oo_{E_i}$ for some $1 \leqslant i \leqslant t$. Then
\begin{enumerate}[(a)]
\item $H^p(L(E_i)) \cong H^p(L)$ for all $p$;
\item $H^2(L(-E_i)) \cong H^2(L)$;
\item If the base locus of $|L|$ does not contain $E_i$, then $H^1(L(-E_i)) \cong H^1(L)$.
\end{enumerate}
\end{lem}

\begin{proof} For (a), consider the short exact sequence
\[
    0 \to L \to L(E_i) \to L(E_i)|_{E_i} \cong \oo_{E_i}(-1) \to 0.
\]
Since $\oo_{E_i}(-1) \cong \oo_{\bp^1}(-1)$ has no cohomology, we get the result. For (b) and (c), consider the short exact sequence
\[
    0 \to L(-E_i) \to L \to L|_{E_i} \cong \oo_{E_i} \to 0.
\]
Since $H^2(\oo_{E_i})=0$, we immediately obtain (b). If the base locus of $|L|$ does not contain $E_i$, then $H^0(L) \to H^0(\oo_{E_i}) \cong k$ is surjective, which gives the result since $H^1(\oo_{E_i})=0$.
\end{proof}

\section{Exceptional sequences}\label{sect:exc-sequence} We review basic facts about Hom functors and exceptional sequences in the bounded derived category of a smooth projective variety $Y$ over $k$. In particular, we discuss how to replace a general complex with a complex built from a strong full exceptional sequence $\mathfrak{E}$, which we call an $\mathfrak{E}${\em -complex}.

%We introduce the notion of {\em cohomological discernibility} for exceptional sheaves and list the exceptional sequences most important for this study.

\subsection{Hom functors}

If $A^{\bullet}$ and $B^{\bullet}$ are bounded complexes of coherent sheaves, then $\Hom^{\bullet}(A^{\bullet},B^{\bullet})$ is the complex of vector spaces defined by
\[
    \Hom^i(A^{\bullet},B^{\bullet}) = \bigoplus_q \Hom(A^q,B^{q+i}) \quad \text{and} \quad d(f) = d_B \circ f - (-1)^i f \circ d_A.
\]
The degree-0 cohomology of this complex is the vector space of chain maps $A^{\bullet} \to B^{\bullet}$ modulo chain homotopy. This complex is especially useful when it computes the derived functor
\[
    R\Hom(A^{\bullet},B^{\bullet}) = \bigoplus_{j \in \mathbb{Z}} \Hom(A^{\bullet},B^{\bullet}[j]),
\]
as in the lemma below. The graded summands $R^j\Hom(A^{\bullet},B^{\bullet}) = \Hom(A^{\bullet},B^{\bullet}[j])$ are denoted $\Ext^j(A^{\bullet},B^{\bullet})$. In the case when $A^{\bullet}$ and $B^{\bullet}$ are sheaves $A$ and $B$ in degree 0,
\[
    R\Hom(A,B) = \bigoplus_{j\geqslant 0} \Ext^j(A,B)[-j]
\]
consists of the usual $\mathrm{Ext}^j$ groups for sheaves in each degree $j$.

Similarly, $\lHom^{\bullet}(A^{\bullet},B^{\bullet})$ is defined by
\[
    \lHom^i(A^{\bullet},B^{\bullet}) = \bigoplus_q \lHom(A^q,B^{q+i}) \quad \text{and} \quad d(f) = d_B \circ f - (-1)^i f \circ d_A.
\]
If either ${A^\bullet}$ or $B^{\bullet}$ is a complex of locally free sheaves, then $\lHom^{\bullet}(A^{\bullet},B^{\bullet})$ represents $R\lHom(A^{\bullet},B^{\bullet})$.

We also recall a few facts about derived functors. The cohomology groups of $R \Gamma$, often denoted $\mathbb{H}^i$, are called hypercohomology, and hypercohomology of a sheaf is just sheaf cohomology. The derived functor $R \Hom$ has the property that $R^i\Hom(A^{\bullet},B^{\bullet}) = \Hom(A^{\bullet},B^{\bullet}[i])$.
Moreover, $R \Gamma R\lHom = R(\Gamma \circ \lHom) = R\Hom$ %(Huybrechts p.85),
and there is a spectral sequence
\[
    E_1^{p,q}= H^p(A^q) \implies \mathbb{H}^{p+q}(A^{\bullet}).
\]
See \cite{Huy06} for more details.

\begin{lem} \label{lem:chain-maps} If $A^{\bullet}$ and $B^{\bullet}$ are any complexes composed of locally free sheaves such that all higher Exts between $A^i$ and $B^j$ vanish for all $i,j$, then $\Hom^{\bullet}(A^{\bullet},B^{\bullet})$ computes $R\Hom(A^{\bullet},B^{\bullet})$. In particular, $\Hom(A^{\bullet},B^{\bullet})$ is the space of chain maps $A^{\bullet} \to B^{\bullet}$ modulo chain homotopy.
\end{lem}

\begin{proof} We calculate $R\Hom(A^{\bullet},B^{\bullet})$ as follows. First, we represent $R\lHom(A^{\bullet},B^{\bullet})$ by $\lHom(A^{\bullet},B^{\bullet})$. Then, since $A^i$ and $B^j$ have no higher Exts between them, $\lHom(A^i,B^j)$ has no higher cohomology, so by the above spectral sequence we can calculate $R\Gamma R\lHom(A^{\bullet},B^{\bullet})$ simply as $\Gamma \lHom(A^{\bullet},B^{\bullet}) = \Hom^{\bullet}(A^{\bullet},B^{\bullet})$. Thus, the complex $\Hom^{\bullet}(A^{\bullet},B^{\bullet})$ represents $R\Hom(A^{\bullet},B^{\bullet})$, which gives the result.
\end{proof}

\subsection{Exceptional sequences} The material reviewed in this section can be found in \cite{GK04}.

\begin{defn}
An object $\mathcal{E}\in \db(Y)$ is {\em exceptional} if \[\Hom(\E,\E[\ell])=\begin{cases} \basefield, &\ell=0 \\
0, & \text{otherwise}.\end{cases}\]
An {\em exceptional sequence} is a sequence $(\E_1,\dots,\E_n)$ of exceptional objects such that \[\Hom(\E_i,\E_j[\ell])= 0,  \mbox{ for } i>j \mbox{ and all } \ell.\]
It is {\em strong} if in addition \[\Hom(\E_i,\E_j[\ell])=0, \mbox{ for all }i,j \mbox{ and }\ell\not=0.\]
It is {\em full} if $\{\E_i\}_{i=1}^n$ generates $\db(Y)$ as a triangulated category.
\end{defn}

%\todo{T: Shall we introduce notation for the exceptional sequence to improve the notation "$\E$-complex"? Maybe a different capital script E, such as $\mathfrak{E}$ or $\mathbb{E}$? Maiorana uses the first of these; Bridgeland-Stern use the second.
%}

Let $\mathfrak{E}=(\E_1,\dots,\E_{n})$ be a strong full exceptional sequence of locally free sheaves on $Y$.
The full triangulated subcategories $\langle \E_i \rangle$ generated by individual $\E_i$
% \todo{Y: Is it better to replace "each" by "individual"? \\ T: Changed.}
yield a semi-orthogonal decomposition of $D^b(Y)$. Thus, for each object $T$ in $D^b(Y)$, there is a diagram of morphisms
\begin{equation}\label{eq:filtration}\xymatrix{
    & C_{n} \ar[dl] && C_{n-1} \ar[dl] & \cdots & C_2 \ar[dl] && C_1 \ar[dl] \\
    T=T_n \ar[rr] && T_{n-1} \ar[rr] \ar[ul]_{[1]} && \cdots \ar[rr] \ar[ul]_{[1]} && T_1 \ar[rr] \ar[ul]_{[1]} && T_0 \cong 0 \ar[ul]_{[1]}
}\end{equation}
in which each triangle $C_i \to T_i \to T_{i-1}$ is distinguished and $C_i$ is in $\langle \E_i \rangle$. Each $T_{i-1}$ can be constructed as the left mutation of $T_i$ through $\E_i$, namely as the cone
\[
    %\bigoplus_j \Hom(\E_i[-j],T^i) \otimes \E_i[-j]
    R\Hom(\E_i,T_i) \otimes \E_i \to T_i \to T_{i-1}.
\]
We call $C_i = R\Hom(\E_i,T_i) \otimes \E_i$ %= \bigoplus_j \Hom(\E_i[-j],T^i) \otimes \E_i[-j]$
the \emph{factors} of $T$ with respect to the exceptional sequence. The diagram is functorial and in particular, the factors of $T$ are unique up to isomorphism.

Using the diagram, the factors of $T$ can be assembled to produce a complex isomorphic to $T$. We define an \emph{$\mathfrak{E}$-complex} to be a bounded complex $A^{\bullet}$ such that each $A^i$ is a direct sum of sheaves in the exceptional sequence.

\begin{lem}\label{lem:E-cpx}
%Each $T^i$ in the $\E$-cofiltration of $T$, and in particular $T$, is isomorphic to an $\E$-complex, with the $\E_i$'s appearing in the same positions as in the factors of $T$. The same is true for each $T_i$ in the $\E$-filtration of $T$.
$T$ is isomorphic to an $\mathfrak{E}$-complex whose sheaves are the same as in the complex $\bigoplus_{i=1}^n C_i$ (but with different maps). %The same is true for each $T_i$ in the $\E$-filtration of $T$.
\end{lem}

\begin{proof}
We prove that each $T_i$ is isomorphic to an $\mathfrak{E}$-complex by induction on $i$. Assume $T_{i-1}$ is isomorphic to an $\mathfrak{E}$-complex $A^{\bullet}$. By the previous lemma, the morphism $T_{i-1}[-1] \to C_i$ can be represented by a chain map $A^{\bullet}[-1] \to C_i$, whose mapping cone is an $\mathfrak{E}$-complex whose sheaves are the same as $A^{\bullet} \oplus C_i$ and which represents $T_i$.
%is an $\E$-complex as well. Because of the shift by 1, the mapping cone contains $C_i$ in the same degrees, which gives the result. %The proof for the $T_i$ is similar by downward induction on $i$.
\end{proof}

\begin{exmp}
Suppose $n=3$ and $T$ is a sheaf in degree 0 that has a resolution $0 \to \E_1^{a_1} \to \E_2^{a_2} \oplus \E_3^{a_3} \to T \to 0$. Then the diagram (\ref{eq:filtration}) can be realized as
\[\xymatrix{
    & \E_3^{a_3} \ar[dl] && \E_2^{a_2} \ar[dl] & & \E_1^{a_1}[1] \ar[dl] \\
    T \cong [\E_1^{a_1} \to \E_2^{a_2} \oplus \E_3^{a_3}] \ar[rr] && [\E_1^{a_1} \to \E_2^{a_2}] \ar[rr] \ar[ul]_{[1]} && \E_1^{a_1}[1] \ar[rr] \ar[ul]_{[1]} && 0 \ar[ul]_{[1]}
}\]
\end{exmp}

The proof gives an inductive algorithm for assembling the factors of $T$ into an $\mathfrak{E}$-complex isomorphic to $T$. We say that an $\mathfrak{E}$-complex $A^{\bullet}$ is \emph{minimal} if, among all $\mathfrak{E}$-complexes isomorphic to $A^{\bullet}$, the total number of sheaves in the complex is as small as possible. If each $C_i$ is represented by the minimal $\mathfrak{E}$-complex described above, then the complex obtained from this algorithm is minimal as well, and we call it the \emph{minimal $\mathfrak{E}$-complex} of $T$. By Lemma \ref{lem:chain-maps}, it is unique up to quasi-isomorphism.

Using $\mathfrak{E}$-complexes to represent objects in the derived category is useful because the sheaves in an exceptional sequence have no higher Exts between them, so Lemma $\ref{lem:chain-maps}$ implies that any morphism between two objects represented by $\mathfrak{E}$-complexes can be realized as a chain map between the $\mathfrak{E}$-complexes.

There is a direct way to identify the factors of $T$ by making use of the dual of the exceptional sequence. The (left) \emph{dual} of a full exceptional sequence $(\E_1,\dots,\E_{n})$ is a full exceptional sequence $(\dE_n,\dots,\dE_1)$
%\todo{T: $({}^\vee \E_n,\dots,\E_1)$}
with the property that
\begin{equation}\label{eq:dual-excl-seq}
    \Hom(\dE_i,\E_j[\ell]) = \begin{cases}
        k, & \text{if $\ell=0$ and $i=j$;} \\
        0, & \text{otherwise.}
        \end{cases}
\end{equation}
The dual sequence always exists, can be constructed from $(\E_1,\dots,\E_n)$ by mutations, and is characterized up to isomorphism by (\ref{eq:dual-excl-seq}).

\begin{lem}\label{lem:cohom-discernible} Suppose $(\E_1,\dots,\E_n)$ is a strong full exceptional sequence and $(\dE_n,\dots,\dE_1)$ is its dual. Then the factors $C_i$ of an object $T$ satisfy
\[
    C_i \cong R\Hom(\dE_i,T) \otimes \E_i.
    %\bigoplus_{j} \Hom(\E_i'[-j],T) \otimes \E_i[-j].
\]
\end{lem}

\begin{proof}
Let $A^{\bullet}$ be the minimal $\mathfrak{E}$-complex of $T$, whose sheaves are the same as $\bigoplus_i C_i$. Then, as the higher Exts between $\dE_i$ and the sheaves in $(\E_1,\dots,\E_n)$ vanish, $\Hom^{\bullet}(\dE_i,A^{\bullet})$ computes $R\Hom(\dE_i,T)$. Thus, by property (\ref{eq:dual-excl-seq}),
\[
    R\Hom(\dE_i,T) \cong \Hom^{\bullet}(\dE_i,C_i) \cong R\Hom(\E_i,T_i),
\]
%and comparing cohomologies in degree $j$ yields $\Hom(\E_i[-j],T^i) \cong \Hom(\E_i'[-j],T)$.
and tensoring by $\E_i$ gives the result.
\end{proof}

\subsection{Main examples of strong full exceptional sequences}\label{ss:ex-sfes}

In later sections, we will focus on the following choices of strong full exceptional sequences of line bundles on $\bp^2$, Hirzebruch surfaces, and, more generally, two-step blowups of Hirzebruch surfaces.
\begin{enumerate}[(a)]
    \item On $\bp^2$, the exceptional sequence $(\oo(-2),\oo(-1),\oo)$ is strong and full (see \cite[Corollary 8.29, Exercise 8.32]{Huy06} for a more general result on $\bp^n$).
    \item On $\mathbb{F}_e$, the exceptional sequence $(\oo(-C-A),\oo(-C),\oo(-A),\oo)$ is strong and full (\cite{Orl92projective-bundles}, \cite{hille-perling-exceptional11} Proposition 5.2).
Note that the exceptional and strong properties can easily be checked using Lemma \ref{lem:lb-van-higher-cohom}.

    \item On $X$ a two-step blowup of $\mathbb{F}_e$, consider the exceptional sequence
    \begin{multline}\label{sequence-blowup}
        \oo(-C-A),\oo(-C-A+E_{s+1}),\dots, \oo(-C-A+E_{t})  ,\\ \oo(-C),\oo(-A), \oo(-E_1),\dots , \oo(-E_s),\oo.
    \end{multline}
    This sequence is obtained from the sequence on $\mathbb{F}_e$ in (b) by standard augmentations, so it is strong and full (\cite{hille-perling-exceptional11} Theorem 5.8). The exceptional and strong properties can easily be checked by using Lemma \ref{lem:reduce-higher-cohom-blowup} to reduce to calculations on $\mathbb{F}_e$.
    Note that this example specializes to (b) by allowing the set of blown-up points to be empty.
\end{enumerate}

The dual sequences are as follows and can be verified by checking (\ref{eq:dual-excl-seq}):
\begin{itemize}
    \item On $\bp^2$, the dual exceptional sequence is
    \[
        \oo, T(-1)[-1], \oo(1)[-2],
    \]
    where $T$ is the tangent sheaf, which can be checked by Bott's formula \cite{Bot57}.
    \item On $\mathbb{F}_e$, the dual exceptional sequence is
    \[
        \oo, \oo(A)[-1], \oo(B)[-1], \oo(A+B)[-2],
    \]
    which follows from the calculations in \S~\ref{subsect:coh-line-bdl}.
    \item On the two-step blowup of $\mathbb{F}_e$, the dual exceptional sequence is
    \begin{align*}
    &\oo, \oo_{E_s}[-1],\dots,\oo_{E_1}[-1], \oo(A - \sum_{i \in S_0} E_i)[-1], \oo(B - \sum_{i \in S_0} E_i)[-1], \\
    & \qquad \qquad \qquad \oo_{E_t}[-2],\dots,\oo_{E_{s+1}}[-2],\oo(A+B - \sum_{i \in S_0} E_i - \sum_{j \in S_1} E_j)[-2].
\end{align*}
This can be seen by using the fact that $\oo(-C-2A-B + \sum_{i \in S_0} E_i + \sum_{j \in S_1} E_j) \cong K_X$, the short exact sequences $0 \to \oo(-E_i) \to \oo \to \oo_{E_i} \to 0$, and Lemma \ref{lem:reduce-higher-cohom-blowup}(a) to reduce to the calculations on $\mathbb{F}_e$.
\end{itemize}
% \todo{Y: Shall we keep the table before? By doing this, we won't have to list the divisors in various places. \\ T: How is the table better than the dual exceptional sequence?}

\section{Gaeta resolutions}\label{sect:gtr}
We study two-step resolutions of coherent sheaves by the exceptional sheaves in the previous section. We call such resolutions {\em Gaeta resolutions}. We provide a general criterion and a criterion specialized for rational surfaces that detect when a sheaf admits a Gaeta resolution. %, via the study of various cohomology groups.
We classify numerical classes over two-step blowups of $\mathbb{F}_e$ of sheaves admitting Gaeta resolutions, including the case of allowing a twist by a line bundle on $\mathbb{F}_e$.
These results lay the foundation for our applications of Gaeta resolutions in later sections.

\subsection{Definition of Gaeta resolutions}\label{subsect:def-gaeta}

Let $\mathfrak{E}=(\E_1,\dots,\E_n)$ be a strong full exceptional sequence on a smooth projective variety $Y$ over $k$. We are particularly interested in minimal $\mathfrak{E}$-complexes of the following form.
\begin{defn}\label{defn:gaeta}
For a coherent sheaf $F$,
a resolution of $F$ of the form
\[
    0 \to \E_1^{a_1} \oplus \cdots \oplus \E_d^{a_d} \to \E_{d+1}^{a_{d+1}} \oplus \cdots \oplus \E_{n}^{a_{n}} \to F \to 0
\]
is a \emph{Gaeta resolution}. If the minimal $\mathfrak{E}$-complex of $F$ is of this form, then we say that $F$ \emph{admits a Gaeta resolution}. The non-negative integers $a_1,\dots,a_n$ are called the \emph{exponents} of the Gaeta resolution.
\end{defn}
Clearly, the exponents $a_i$ determine the numerical class of $F$. Conversely, the class of $F$ determines the exponents inductively using semi-orthogonality as
\[
    a_i = \begin{cases}
    \chi(\E_i,F) - \sum_{j=i+1}^n a_j \hom(\E_i,\E_j) & \text{for $i > d$;} \\
    -\chi(F,\E_i) - \sum_{j = 1}^{i-1} a_j \hom(\E_j,\E_i) & \text{for $i \leqslant d$}.
    \end{cases}
\]
The dual sequence $(\dE_n,\dots,\dE_1)$ can be used to obtain a criterion for when a sheaf admits a Gaeta resolution.

\begin{prop}[General criterion] \label{prop:general-GTR-criterion} Let $F$ be a coherent sheaf. Then $F$ admits a Gaeta resolution if and only if
\[
    \Hom(\dE_i[j],F)
\]
vanish for all $i$ and $j$ except possibly for
\begin{align*}
    &\Hom(\dE_i,F), \quad d+1 \leqslant i \leqslant n; \\
    &\Hom(\dE_i[1],F), \quad 1 \leqslant i \leqslant d.
\end{align*}
Moreover, if $F$ admits a Gaeta resolution, then the exponents are
\[
    a_i = \begin{cases}
        \hom(\dE_i,F), & d+1 \leqslant i \leqslant n; \\
        \hom(\dE_i[1],F), & 1 \leqslant i \leqslant d.
    \end{cases}
\]
\end{prop}

\begin{proof} By Lemma \ref{lem:cohom-discernible} and (\ref{eq:dual-excl-seq}), the condition on $\Hom(\dE_i[j],F)$ is equivalent to $C_i$ being a direct sum of copies of $\E_i$ in degree 0 if $d+1 \leqslant i \leqslant n$ and in degree $-1$ if $1 \leqslant i \leqslant d$. By Lemma \ref{lem:E-cpx}, this proves the first statement. The second statement also follows from the calculation of the $C_i$.
\end{proof}

\subsection{Gaeta resolutions on rational surfaces}

In the context of the strong full exceptional sequences in \S~\ref{ss:ex-sfes}, we will focus on the following Gaeta resolutions.
\begin{exmp}\label{exmp:gaeta}
\begin{enumerate}[(a)]
    \item On $\bp^2$, we consider Gaeta resolutions of the form
\begin{align}\label{gtr-P2}
    0 \to \oo(-2)^{\alpha_1} \to \oo(-1)^{\alpha_2} \oplus \oo^{\alpha_3} \to & F \to 0.
\end{align}

    \item On $X$ a two-step blowup of $\mathbb{F}_e$, we consider Gaeta resolutions of the form
    \begin{align}\label{gtr-blowup}
        0 \to \oo(-C-A)^{\alpha_1} &\oplus \bigoplus_{j \in S_1} \oo(-C-A+E_j)^{\gamma_j} %\to
        \nonumber \\
        &\to \oo(-C)^{\alpha_2} \oplus \oo(-A)^{\alpha_3} \oplus \bigoplus_{i \in S_0} \oo(-E_i)^{\gamma_i} \oplus \oo^{\alpha_4} \to F \to 0.
    \end{align}

\end{enumerate}
Note that in the case when the set of blown-up points is empty, (\ref{gtr-blowup}) specializes to the Gaeta resolutions on $\mathbb{F}_e$ that were considered in \cite{CosHui18weakBN}.
\end{exmp}

For these examples,
%by making use of cohomological discernibility,
by applying the general criterion for having a Gaeta resolution (Proposition \ref{prop:general-GTR-criterion})
with the explicit dual sequences in \S~\ref{ss:ex-sfes},
we deduce a more explicit criterion.

\begin{prop}\label{prop:special-GTR-criterion}
\begin{enumerate}[(a)]
\item On $\bp^2$, a sheaf $F$ admits a Gaeta resolution of the form (\ref{gtr-P2})
    if and only if
\[
    \text{$H^p(F)=0$ for $p \ne 0$, \quad $H^p(F(-1)) = 0$ for $p \ne 1$, \quad and \quad $\Hom(T(-1),F)=0$.}
\]
%and the natural map $\Hom(\oo(-1),\oo) \otimes H^0(F) \to H^0(F(1))$ is injective.

\item On $X$ a two-step blowup of $\mathbb{F}_e$, a torsion-free sheaf $F$ admits a Gaeta resolution of the form (\ref{gtr-blowup})
    if and only if
\begin{enumerate}[(i)]
    \item $H^p(F)$ vanishes for $p \ne 0$;
    %\item $H^p(F(-A+\sum_{i \in S_0} E_i))$, $H^p(F(-B+\sum_{i \in S_0} E_i))$, and $H^p(F(-A-B+\sum_{i \in S_0} E_i + \sum_{j \in S_0} E_j))$ vanish for $p \ne 1$; and
    \item $H^p(F(D))$ vanishes for $p \ne 1$ and $D$ the divisors $-A+\sum_{i \in S_0} E_i$, $-B+\sum_{i \in S_0} E_i$, and $-A-B+\sum_{i \in S_0} E_i + \sum_{j \in S_0} E_j$; and
    %all $D \in \{\, -A+\sum_{i \in S_0} E_i,-B+\sum_{i \in S_0} E_i,-A-B+\sum_{i \in S_0} E_i + \sum_{j \in S_0} E_j \, \}$; and
    \item $H^1(F(E_i))=0$ and $H^1(F(E_j)) = 0$ for all $i \in S_0$ and $j \in S_1$.
\end{enumerate}
\end{enumerate}
\end{prop}

% \todo{T: I believe there is an error in (b.i). $H^p(F(-A+\sum_{i \in S_0} E_i))$ and $H^p(F(-B+\sum_{i \in S_0} E_i))$ should vanish for $p \ne 1$ as in (b.ii). They should be moved into (ii).\\Y: Corrected. \\T: (iii) is a proposed restatement of (ii).}

\begin{proof} For (a), Proposition \ref{prop:general-GTR-criterion} includes the first two conditions as well as the condition that $\Ext^p(T(-1),F)$ vanishes for $p=0,2$. Applying $\Hom(-,F)$ to the short exact sequence $0 \to T(-1) \to \oo(1)^3 \to \oo(2) \to 0$ and using the vanishing of $H^2(F(-1))$ shows that the first two conditions already guarantee $\Ext^2(T(-1),F)=0$. %Applying $\Hom(-,F)$ to the short exact sequence $0 \to \oo(-1) \to \oo^3 \to T(-1) \to 0$ shows that the vanishing of $\Hom(T(-1),F)$ is equivalent to the injectivity of $\Hom(\oo^3,F) \to \Hom(\oo(-1),F)$, which can be identified with the map in (iii).

For (b), Proposition \ref{prop:general-GTR-criterion} includes (i) and (ii) as well as the condition that $\Ext^p(\oo_{E_i},F)$ vanishes for $p=0$ and $p=2$ for all $i \in S_0 \cup S_1$. The vanishing $\Hom(\oo_{E_i},F)=0$ is guaranteed since $F$ is torsion-free, while applying $\Hom(-,F)$ to the short exact sequence $0 \to \oo(-E_i) \to \oo \to \oo_{E_i} \to 0$ and using $H^1(F)=H^2(F)=0$ shows that $\Ext^1(\oo(-E_i),F) \cong \Ext^2(\oo_{E_i},F)$, hence $H^1(F(E_i))=0$ is an equivalent condition.
\end{proof}

\subsection{Chern characters and Gaeta resolutions}\label{ss:chern-characters}

Recall that the exponents in the Gaeta resolutions are determined by the numerical class of $F$. The results in this subsection classify numerical classes that arise as cokernels of Gaeta resolutions. First, we review some useful numerical invariants.

On a smooth projective surface $S$ over $k$, if $F$ is a coherent sheaf of positive rank $r$, first Chern class $c_1$, and second Chern character $\ch_2$, set
\[
    \nu = \frac{c_1}{r} \qquad \text{and} \qquad \Delta = \frac{1}{2} \nu^2 - \frac{\ch_2}{r},
\]
which are called the $\emph{total slope}$ and $\emph{discriminant}$, respectively, of $F$. A simple calculation shows that the discriminant of a line bundle is 0 and the discriminant of $F$ is unchanged when $F$ is tensored by a line bundle.
Using these invariants, the second Chern class of $F$ can be written as
\[
    c_2 = \binom{r}{2} \nu^2 + r \Delta.
\]
Set $P(\nu) = \chi(\oo_S)+\frac{1}{2} \nu(\nu-K_S)$. Then by Riemann-Roch the Euler characteristic can be written as
\[
    \chi(F) = r(P(\nu)-\Delta),
\]
and similarly, if $F_1,F_2$ are sheaves and $r_i$, $\nu_i$, $\Delta_i$ are the rank, total slope, and discriminant of $F_i$, then the Euler pairing is
\begin{equation}\label{eq:euler-pair}
    \chi(F_1,F_2) = r_1 r_2(P(\nu_2-\nu_1) - \Delta_1 - \Delta_2).
\end{equation}

On $X$ the two-step blowup of a Hirzebruch surface, writing $c_1 = \alpha A + \beta B - \sum_{i \in S_0} \gamma_i E_i - \sum_{j \in S_1} \gamma_j E_j$, we compute
\[
    P(\nu) = \left( \frac{\alpha}{r}+1- \frac{e \beta}{2r} \right) \left( \frac{\beta}{r}+1 \right) - \frac{1}{2} \sum_{i \in S_0 \cup S_1} \frac{\gamma_i}{r}\left( \frac{\gamma_i}{r}+1 \right).
    %- \sum_j \frac{\gamma_j}{2r}\left( \frac{\gamma_j}{r}+1 \right).
\]

We write the numerical class of a sheaf as a triple $(r,c_1,\chi)$ in which $r$ is a non-negative integer, $c_1$ is an integral divisor class, and $\chi$ is an integer.
We say that a numerical class $f = (r,c_1,\chi)$ \emph{admits Gaeta resolutions} if there is a sheaf $F$ of rank $r$, first Chern class $c_1$, and Euler characteristic $\chi$, such that $F$ admits a Gaeta resolution. For a sheaf $F$ of class $f$ and a line bundle $L$, we denote the class of $F\otimes L$ as $f(L)$, and we write $c_2(f)$, $\nu(f)$ and $\Delta(f)$ for the second Chern class, total slope, and discriminant of $F$, which depend only on $f$.

\begin{prop}\label{prop:chern-has-gr}
On $X$ a two-step blowup of $\mathbb{F}_e$, consider the numerical class
\[
    f=\Big(r,\alpha A + \beta B - \sum_{i \in S_0 \cup S_1} \gamma_i E_i,
    %- \sum_{j \in S_1} \gamma_j E_j,
    \chi \Big)
\]
of positive rank. Then
\begin{enumerate}[(a)]
\item $f$ admits Gaeta resolutions (\ref{gtr-blowup}) if and only if $\gamma_i$, $\gamma_j$, $\alpha_4:=\chi$, and the following three integers are all $\geqslant 0$:
\begin{align*}
    \alpha_1 &:= -\chi\Big(f\Big(-A-B+ \sum_{i \in S_0} E_i + \sum_{j \in S_1} E_j\Big)\Big) \\
    \alpha_2 &:= -\chi\Big(f\Big(-B+\sum_{i \in S_0} E_i\Big)\Big) \\
    \alpha_3 &:= -\chi\Big(f\Big(-A+\sum_{i \in S_0} E_i\Big)\Big) %\\
    %\alpha_4 &= \chi
\label{eq:chi-ineq-blowup}
\end{align*}

\item Assume $\gamma_i \geqslant 0$ and $\gamma_j \geqslant 0$. If the discriminant $\Delta(f)$ is sufficiently large,
%satisfies \[\Delta(f) \gg 1 + \max\Big(\sum_{i \in S_0} \gamma_i, \sum_{j \in S_1} \gamma_j\Big)/r,\]
% \todo{Y: I think the expression is unnecessarily explicit. \\
% T: Changed.}
then there is a line bundle $L$ pulled back from $\mathbb{F}_e$ such that $f(L)$ admits Gaeta resolutions.
\end{enumerate}
\end{prop}

\begin{proof}
For (a), assuming $f$ admits Gaeta resolutions, by comparing first Chern classes, the exponent of $\oo(-C-A+E_j)$ must be $\gamma_j$ and the exponent of $\oo(-E_i)$ must be $\gamma_i$. The remaining exponents can be easily calculated. % using cohomological discernability and are exactly $\alpha_1,\alpha_2,\alpha_3$ as stated and $\alpha_4 = \chi$.
Conversely, the inequalities show that we can define the exponents in the same way, and a simple calculation shows that the numerical class of the cokernel must be $f$.

For (b), consider the numerical class $f' = \Big(r, \alpha A + \beta B, \chi + \sum_{i \in S_0} \gamma_i \Big)$
on $\mathbb{F}_e$. An elementary calculation shows that
\[
    \Delta(f') = \Delta(f) + \sum_{i \in S_0} \binom{\gamma_i/r}{2} + \sum_{j \in S_1} \binom{(\gamma_j/r)+1}{2}.
\]
% \todo{Y: Is it common to use the binomial symbol for rational numbers? \\ T: I think it's pretty standard.}
Let $M = 1 + \max(\sum_{i \in S_0} \gamma_i, \sum_{j \in S_1} \gamma_j)/r$. Since $\Delta(f') \geqslant \Delta(f) \gg M$, the following lemma ensures that we can choose a line bundle $L$ on $\mathbb{F}_e$ such that $f'(L)$ admits Gaeta resolutions
\[
    \oo(-C-A)^{\alpha_1'} \to \oo(-C)^{\alpha_2'} \oplus \oo(-A)^{\alpha_3'} \oplus \oo^{\alpha_4'}
\]
in which $\alpha_4' = \chi(f'(L)) \geqslant rM = r + \max\Big(\sum_{i \in S_0} \gamma_i, \sum_{j \in S_1} \gamma_j\Big)$. Here, the inequality follows from the following lemma.  Hence, $\alpha_1' \geqslant \max(\sum_{i \in S_0} \gamma_i,\sum_{j \in S_1} \gamma_j)$ as well by comparing ranks. Then a simple calculation shows that the cokernels of Gaeta resolutions
\[
    \oo(-C-A)^{\alpha_1' - \sum_j \gamma_j} \oplus \bigoplus_{j \in S_1} \oo(-C-A+E_j)^{\gamma_j} \to \oo(-C)^{\alpha_2'} \oplus \oo(-A)^{\alpha_3'} \oplus \bigoplus_{i \in S_0} \oo(-E_i)^{\gamma_i} \oplus \oo^{\alpha_4'-\sum_i \gamma_i}
\]
have numerical class $f(L)$.
\end{proof}

\begin{lem}
On $\mathbb{F}_e$, fix a rank $r > 0$ and a first Chern class $c_1$ and consider the numerical class $f = (r,c_1,\chi)$. Let $M$ be a positive real number. Then there are constants $C_1,C_2,C_3$ depending only on $e$ such that for all $\chi$ such that
\[
    \Delta(f) \geqslant \frac{e+2}{2}M^2 + C_1 M^{3/2} + C_2 M + C_3,
\]
there is a line bundle $L$ such that $\chi(f(L)) \geqslant rM$ and the class $f(L)$ admits Gaeta resolutions.
\end{lem}

\begin{proof} We use a setup similar to \cite[Lemma 4.5]{CosHui20}.
Consider the curve $Q \colon \chi(f(L_{x,y}))=0$ in the $xy$-plane, where $L_{x,y}$ is the (in general non-integral) line bundle
\[
    L_{x,y} = xB + yA -\nu(f) + \frac{1}{2} K_{\mathbb{F}_e} .
\]
Set $\Delta = \Delta(f)$. By Riemann-Roch,
\[
    \frac{\chi(f(L_{x,y}))}{r} = %((x-1)+1)((y-(e+2)/2) + 1 - e(x-1)/2) - \Delta =
    x\left(y-\frac{e}{2}x \right)-\Delta,
\]
so $Q$ is the hyperbola $\Delta = x \left(y-\frac{e}{2}x \right)$, or, as a function of $x$, $y = Q(x) = \frac{\Delta}{x} + \frac{e}{2}x$.

Let $\Lambda$ denote the lattice in the plane of points such that $L_{x,y}$ is integral, which is a shift of the standard integral lattice. We say that a point $(x,y) \in \Lambda$ is \emph{minimal} if
\begin{itemize}
    \item $(x,y)$ is on or above the upper branch $Q_1$ of $Q$, and
    \item $(x-1,y)$ and $(x,y-1)$ are both on or below $Q_1$.
\end{itemize}
The minimal points exactly correspond to the line bundles $L_{x,y}$ such that $f(L_{x,y})$ admits Gaeta resolutions, according to Proposition~\ref{prop:chern-has-gr}(a).

We need to find a minimal point $(x,y)$ such that $\chi(f(L_{x,y}))/r \geqslant M$. For this, consider the line $y=(e+1)x$, which intersects $Q_1$ at the point
\[
    (x',y') = \left(\sqrt{2\Delta/(e+2)}, (e+1)\sqrt{2 \Delta / (e+2)}\right).
\]
The tangent line to $Q_1$ at this point has equation $y = -x + \sqrt{2(e+2)\Delta}$
and $Q_1$ lies above the tangent line. See Figure 1.

\begin{figure}[h]
\begin{minipage}{.5\textwidth}
\centering
\includegraphics[height=3in]{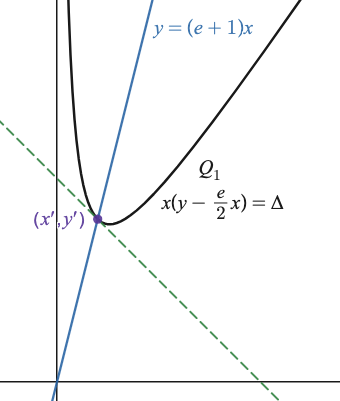}
\captionof{figure}{The hyperbola $Q_1$}
\label{fig:1}
\end{minipage}%
\begin{minipage}{.5\textwidth}
\centering
\includegraphics[height=3in]{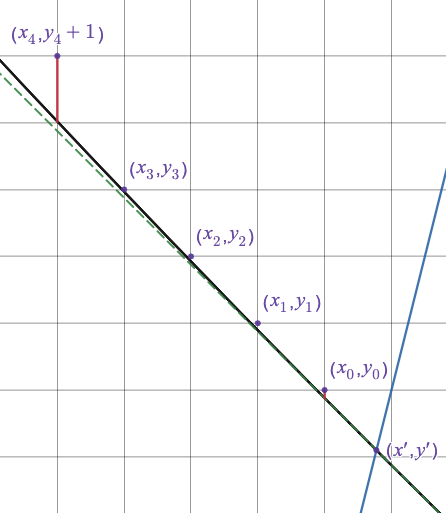}
\captionof{figure}{Example of minimal points near $(x',y')$}
\label{fig:2}
\end{minipage}
\end{figure}

Let $(x_0,y_0)$ be the unique minimal point such that $\epsilon_x := x'-x_0$ satisfies $0 \leqslant \epsilon_x < 1$. which lies between $Q_1$ and the shift $Q_1 + (0,1)$. Then $\epsilon_y := y_0 - y'$ satisfies $0 \leqslant \epsilon_y \leqslant 2$. Let $\epsilon := y_0 - Q(x_0)$ denote the vertical distance from $(x_0,y_0)$ to $Q_1$, which satisfies $0 \leqslant \epsilon < 1$. Then $\chi(f(L_{x_0,y_0}))/r = \epsilon x_0$. % >\epsilon(x'-1)$$
If $\epsilon x_0 \geqslant M$, then $L_{x_0,y_0}$ gives the result, so assume on the contrary that
$\epsilon < M/{x_0}$.
Then let $m$ be a positive integer $<x_0$ and consider the lattice point $(x_m,y_m) = (x_0 - m, y_0 + m)$. We wish to find $m$ as small as possible such that $Q_1$ lies above $(x_m,y_m)$, as then the point $(x_m,y_m+1)$ will be minimal and $y_m+1 - Q(x_m)$ will be close to 1. See Figure 2 for an example in which this is achieved with $m=4$.

To find $m$ such that $Q_1$ lies above $(x_m,y_m)$, the rise in $Q$ between $x_0$ and $x_0-m$ should exceed $m+\epsilon$, namely
\begin{align*}
    Q(x_0-m) - Q(x_0) - m - \epsilon
    = \frac{m \Delta}{x_0^2}\left( 1 + \frac{m}{x_0} + \left(\frac{m}{x_0}\right)^2 + \cdots \right) - \frac{e+2}{2}m - \epsilon
\end{align*}
should be positive. As this quantity exceeds the approximation obtained by truncating the geometric series at the first two terms, and as $x_0 \leqslant x'$, it suffices to take $m$ such that
\[
    \frac{m\Delta}{(x')^2}\left(1+\frac{m}{x_0} \right)-\frac{e+2}{2}m - \epsilon \geqslant 0,
\]
which yields $m \geqslant \sqrt{2\epsilon x_0/(e+2)}$.
Setting $m_0 = \lceil \sqrt{2\epsilon x_0/(e+2)} \rceil$, we then have
\begin{align*}
    \frac{\chi(L_{x_{m_0},y_{m_0+1}})}{r}
    &= (1+\epsilon) x_0 - \epsilon_x(e+1)m_0 - \epsilon_y m_0 - \frac{e+2}{2}m_0^2 - m_0,
\end{align*}
and for this to be $\geqslant M$ we need
\[
    (1+\epsilon) x_0 \geqslant \frac{e+2}{2}m_0^2 + (\epsilon_x(e+1)+\epsilon_y+1)m_0 + M.
\]
Let $\delta = m_0 - \sqrt{2 \epsilon x_0/(e+2)}$, which satisfies $0 \leqslant \delta < 1$. Then the inequality simplifies to
\[
    x_0 \geqslant (\delta(e+2)+\epsilon_x(e+1)+\epsilon_y+1)\sqrt{2\epsilon x_0/(e+2)} + (\epsilon_x(e+1)+\epsilon_y+1)\delta + \frac{e+2}{2}\delta^2 + M.
\]
Replacing $x_0$ by $x'-\epsilon_x$, $\epsilon x_0$ by its upper bound $M$, $\epsilon_x$, $\epsilon_y$, and $\delta$ by their upper bounds, and solving for $\Delta$, we get a sufficient bound for $\Delta$:
\[
    \Delta \geqslant \frac{e+2}{2}M^2 + C_1 M^{3/2} + C_2 M + C_3
\]
for constants $C_1,C_2,C_3$ that depend only on $e$.
\end{proof}

\section{Properties of sheaves with general Gaeta resolutions}
\label{sect:prop-gr}

Given an exceptional sequence $(\E_1 ,\E_2,\dots,\E_n)$ from \S~\ref{ss:ex-sfes} with $d$ chosen as in Example \ref{exmp:gaeta} and a sequence of non-negative integers $\vec{a} = (a_1,\dots,a_n)$,
consider the vector space
\[
    H_{\vec{a},d}=\Hom(\E_1^{a_1} \oplus \cdots \oplus \E_d^{a_d},\E_{d+1}^{a_{d+1}} \oplus \cdots \oplus \E_n^{a_n}).
\]
We let $\rcone \subset H_{\vec{a},d}$ denote the open subset of injective maps
\[
    \phi \colon \E_1^{a_1} \oplus \cdots \oplus \E_d^{a_d} \to \E_{d+1}^{a_{d+1}} \oplus \cdots \oplus \E_n^{a_n},
\]
which is non-empty if and only if
$r := a_{d+1} + \cdots + a_n - a_1 - \cdots - a_d \geqslant 0$. In this case, we set \begin{equation*}
    F_\phi:=\coker \phi,
\end{equation*}
let $f$ denote the numerical class of these cokernels, which have rank $r$, and call the following projectivization the \emph{space of Gaeta resolutions}:
\begin{equation}\label{eq:res-space}
R_f := \bp \rcone \subset \bp(H_{\vec{a},d}).
\end{equation}
Then $F_{\phi}$ satisfies various cohomology vanishing conditions (Proposition \ref{prop:special-GTR-criterion}). The purpose of this section is to prove additional properties of sheaves admitting {\em general} Gaeta resolutions, including the prioritary condition and a weak Brill-Noether result.

\subsection{Basic properties}

We begin by proving some basic properties of sheaves admitting a general Gaeta resolution.

\begin{prop}\label{prop:cokernel-properties}
{\ }\begin{enumerate}[(a)]
    \item On $\bp^2$ and $\mathbb{F}_e$, we have the following for general $\phi$:
    \begin{enumerate}[(i)]
    \item If $r = 0$, then $F_{\phi}$ is torsion supported on the determinant of $\phi$ (true even if $\phi$ is not general).
    \item If $r = 1$, then $F_\phi$ is torsion-free.
    \item If $r \geqslant 2$, then $F_\phi$ is locally free.
    \item If $a_n \geqslant r+2$, then $F_\phi$ is globally generated.
    \end{enumerate}
\item On $X$ an admissible blowup of $\mathbb{F}_e$, using the notation in (\ref{gtr-blowup}), the same cases are true if we assume, for each $i \in S_0$ such that $\{j \colon p_j \succ p_i \}$ is nonempty:
\begin{align*}
    \gamma_i &\geqslant \sum_{j \colon p_j \succ p_i} \gamma_j &\text{for (ii)}, \\
    \gamma_i &\geqslant 1 + \sum_{j \colon p_j \succ p_i} \gamma_j &\text{for (iii)}, \\
    \gamma_i &\geqslant r - \alpha_4 + \sum_{j \colon p_j \succ p_i} \gamma_j &\text{for (iv)}.
\end{align*}
\end{enumerate}
\end{prop}

The key to the proof of this proposition is a Bertini-type statement concerning the codimension on which a general map between vector bundles $\mathcal{A}$ and $\mathcal{B}$ drops rank. The case when $\mathcal{B} \otimes \mathcal{A}^*$ is globally generated is well known.

\begin{prop}\label{prop:Bertini-type} On a smooth projective variety $Y$, consider maps
$%\[
    \phi \colon \mathcal{A} \to \mathcal{B}
$, %\]
where $\mathcal{A}$, $\mathcal{B}$ are fixed vector bundles of ranks $a,b$ such that $a \leqslant b$. Let
\[
    Z \subset Y \times \bp \Hom(\mathcal{A},\mathcal{B})
\]
denote the locus of pairs $(y,\phi)$ such that the rank of $\phi$ at $y$ is $<a$.
\begin{enumerate}[(a)]
    \item Suppose $\mathcal{A}^* \otimes \mathcal{B}$ is globally generated.
    Then the codimension of $Z$ in $Y \times \bp \Hom(\mathcal{A},\mathcal{B})$ is $b-a+1$.
    \item Suppose there are non-trivial decompositions
    %\[
        $\mathcal{A} = \mathcal{A}_1 \oplus \mathcal{A}_2$
        and $\mathcal{B}=\mathcal{B}_1 \oplus \mathcal{B}_2$
    %\]
    such that each $\mathcal{A}_i^* \otimes \mathcal{B}_j$ is globally generated except that $H^0(\mathcal{A}_2^* \otimes \mathcal{B}_1)=0$. Set $a_i = \mathrm{rk}(\mathcal{A}_i)$ and $b_j = \mathrm{rk}(\mathcal{B}_j)$. If $b_2 < a_2$, then $Z = Y \times \bp \Hom(\mathcal{A},\mathcal{B})$.
    If $b_2 \geqslant a_2$,
    then the codimension of $Z$ in $Y \times \bp \Hom(\mathcal{A},\mathcal{B})$ is $\min(b_2-a_2,b-a)+1$.
\end{enumerate}
\end{prop}

\begin{proof} Part (a) is a special case of
\cite[Teorema 2.8]{ottaviani1995varieta}  (see also \cite[Proposition 2.6]{Hui16-interpolation})
and can be proved as follows. Consider the map of global sections
\[
    \pi \colon \Hom(\mathcal{A},\mathcal{B}) \otimes \oo_Y \to \lHom(\mathcal{A},\mathcal{B}),
\]
which is surjective since $\mathcal{A}^* \otimes \mathcal{B}$ is globally generated. This map induces a map of projective bundles
\[
    \ev \colon Y \times \bp \Hom(\mathcal{A},\mathcal{B}) \to \bp \lHom(\mathcal{A},\mathcal{B}), \qquad (y,\phi \colon \mathcal{A} \to \mathcal{B}) \mapsto (y,\phi|_y \colon \mathcal{A}|_y \to \mathcal{B}|_y).
\]
Let $\Sigma$ denote the locus in the target of points $(y,\phi_y \colon \mathcal{A}|_y \to \mathcal{B}|_y)$ such that $\phi_y \colon \mathcal{A}|_y \to \mathcal{B}|_y$ drops rank, and let $Z = \ev^{-1}(\Sigma)$. As $\Sigma$ has codimension $b-a+1$ in each $\bp \Hom(\mathcal{A}|_y,\mathcal{B}|_y)$ and $\pi|_y$ is surjective, $Z$ has codimension $b-a+1$ in $Y \times \bp \Hom(\mathcal{A},\mathcal{B})$.

We prove (b) by adapting this argument. In this case, $\pi$ is not surjective, so the codimension of $Z$ may drop if the image of $\ev$ is not transversal to $\Sigma$. At each point $y$, fixing bases, $\phi|_y$ is a $b \times a$ matrix of the form
\[
    \begin{bmatrix}
    M & 0 \\
    N & P
    \end{bmatrix}
\]
in which $M,N,P$ are general if $\phi$ is general. Since $b \geqslant a$, this matrix drops rank if and only if the columns are linearly dependent. If $b_2 < a_2$, the columns of $P$ are always linearly dependent, while if $b_2 \geqslant a_2$, there are two cases to consider in which the columns are linearly dependent:
\begin{enumerate}[(i)]
    \item The columns of $P$ are linearly dependent. This occurs in codimension $1+b_2-a_2$ in the space of such block matrices.
    \item The columns of $P$ are linearly independent. Then the linear dependence involves a column $\vec{c}$ of $\begin{bmatrix} M \\ N \end{bmatrix}$, hence $\vec{c}$ is in the span of the remaining columns, which occurs in codimension $1+b-a$. (For such $\phi$, $\Sigma$ intersects the image of $\ev$ transversally.)
\end{enumerate}
Thus, $Z$ has codimension $1+\min(b_2-a_2,b-a)$ in $Y \times \bp\Hom(\mathcal{A},\mathcal{B})$.

\end{proof}

\begin{proof}[Proof of Proposition \ref{prop:cokernel-properties}]

For (a), as each $\E_i \otimes \E_j^*$ for $j \leqslant d$ and $i > d$ is globally generated,  Proposition \ref{prop:Bertini-type} (a) implies that for general $\phi$ the locus where $\phi$ drops rank is either empty or has codimension $r+1$, which proves $(i)$ and $(iii)$. Then $(ii)$ follows from the fact that $F_\phi$ cannot have zero-dimensional torsion as it has a two-step resolution by locally free sheaves. For (iv), there is a commutative diagram
\[\xymatrix{
    && \oo^{a_n} \ar[d] \ar@{=}[r] & \oo^{a_n} \ar[d] \\
    0 \ar[r] & \E_1^{a_1} \ar[r]^-{\phi} & \E_2^{a_2} \oplus \cdots \oplus \E_{n-1}^{a_{n-1}} \oplus \oo^{a_n} \ar[r] & F_{\phi} \ar[r] & 0
}\]
with exact rows, hence $\oo^{a_n} \to F_\phi$ is surjective if and only if the induced map
\[
    \E_1^{a_1} \to \E_2^{a_2} \oplus \cdots \oplus \E_{n-1}^{a_{n-1}}
\]
is surjective. This is a map between a bundle of rank $a_1$ and a bundle of rank
\[
    a_2 + \cdots + a_{n-1} = a_1 + r - a_n \leqslant a_1 - 2,
\]
hence by dualizing the same Bertini-type statement, a general such map is surjective on all fibers.

For (b), set
\begin{align*}
    \mathcal{F}&= \oo(-C-A)^{\alpha_1} \oplus \bigoplus_{j \in S_1} \oo(-C-A+E_j)^{\gamma_j}\quad \mbox{and}\\
    \mathcal{G}&= \oo(-C)^{\alpha_2} \oplus \oo(-A)^{\alpha_3} \oplus \bigoplus_{i \in S_0} \oo(-E_i)^{\gamma_i} \oplus \oo^{\alpha_4}.
\end{align*}
The complication is that for each $j \in S_1$, the global sections of the line bundle
\[
    L \otimes \oo(C+A-E_j)
\]
vanish on $E_i - E_j$, where $i$ is the index such that $p_i \prec p_j$ and $L$ is any line bundle in $\mathcal{G}$ except for $\oo(-E_i)$. Still, we can adapt the proof of the Bertini-type statement as follows.

Consider
\[
    \pi \colon \Hom(\mathcal{F},\mathcal{G}) \otimes \oo_X \to \lHom(\mathcal{F},\mathcal{G}),
\]
which is not surjective on each $E_i-E_j$, and the induced map
\[
    \ev \colon X \times \bp \Hom(\mathcal{F},\mathcal{G}) \to \bp \lHom(\mathcal{F},\mathcal{G}), \qquad (y,\phi \colon \mathcal{F} \to \mathcal{G}) \mapsto (y,\phi|_y \colon \mathcal{F}|_y \to \mathcal{G}|_y).
\]
Let $\Sigma$ denote the locus in the target where the linear maps drop rank and $Z=\ev^{-1}(\Sigma)$. At points $p$ in the open complement $U$ of the exceptional locus $\bigsqcup_{i \in S_0} E_i$, $\pi|_p$ is surjective, hence the codimension of $Z|_U$ in $U \times \bp \Hom(\mathcal{F},\mathcal{G})$ is $r+1$. On each $\tilde{E}_i$, we apply Proposition \ref{prop:Bertini-type} with
\[
    Y=\tilde{E}_i, \quad \mathcal{A}_2 = \bigoplus_{j \colon p_j \succ p_i} \oo(-C-A+E_j)^{\gamma_j}|_{\tilde{E}_i}, \quad \mathcal{B}_2 = \oo(-E_i)^{\gamma_i}|_{\tilde{E}_i}
\]
and $\mathcal{A}_1,\mathcal{B}_1$ the restrictions of the remaining line bundles. Assuming $\gamma_i \geqslant \sum_{j \colon p_j \succ p_i} \gamma_j$ and using the fact that the global sections of each line bundle summand of $\mathcal{A}^* \otimes \mathcal{B}$ lift to $X$, we deduce that the codimension of $Z|_{\tilde{E}_i}$ in $\tilde{E}_i \times \bp \Hom(\mathcal{F},\mathcal{G})$ is
\[
    \min(\gamma_i - \sum_{j \colon p_j \succ p_i} \gamma_j, r)+1.
\]
On $E_j$, a similar argument with
\[
    Y = {E}_j, \quad \mathcal{A}_2 = \bigoplus_{j' \ne j \colon p_{j'} \succ p_i} \oo(-C-A+E_{j'})^{\gamma_{j'}}|_{E_j}, \quad \mathcal{B}_2 = \oo(-E_i)^{\gamma_i}|_{E_j}
\]
shows that $Z|_{E_j}$ has codimension
\[
    \min(\gamma_i + \gamma_j - \sum_{j \colon p_j \succ p_i} \gamma_j, r) + 1
\]
in $E_j \times \bp\Hom(\mathcal{F},\mathcal{G})$. Altogether, since $\tilde{E}_i$ and $E_j$ have codimension 1 in $X$, we deduce that the general fiber of $Z$ over $\bp \Hom(\mathcal{F},\mathcal{G})$ is empty or has codimension at least
\[
    \min(r,\gamma_i-\sum_{j \colon p_j \succ p_i}\gamma_j + 1) + 1
\]
in $X$. As this codimension is $\geqslant 2$ in the case $r=1$ and $\geqslant 3$ in the case $r \geqslant 2$ assuming $\gamma_i \geqslant 1 + \sum_{j \colon p_j \succ p_i} \gamma_j$ for all $i \in S_0$, we get the result for (b.ii) and (b.iii).

For (b.iv), we dualize and use a similar argument.
\end{proof}

\subsection{Prioritary sheaves}\label{sect:prioritary}

We relate sheaves which admit general Gaeta resolutions to prioritary sheaves, which will facilitate the study of stability in \S~\ref{sec:stability}. % and stable sheaves.
We begin by reviewing the prioritary condition.
\begin{defn}
Let $S$ be a smooth surface and $D$ be a divisor on it. A coherent sheaf $F$ on $S$ is $D$-{\em prioritary} if it is torsion-free and $\Ext^2(F,F(-D))=0$.
If $S$ is a Hirzebruch surface or its blowup, $A$-prioritary sheaves are simply called {\em prioritary sheaves}.
\end{defn}
We will need the following lemma (\cite[Lemma 3.1]{coskun_existence_2019}) comparing prioritary conditions with respect to different divisors.
\begin{lem}\label{lem:prioritary-div}
    Let $S$ be a smooth surface, $D_1$ and $D_2$ be two divisors such that $D_1\geqslant D_2$. Then $D_1$-prioritary sheaves are $D_2$-prioritary.
\end{lem}

Now let $X$ denote an admissible blowup of $\mathbb{F}_e$ and $f$ be a fixed class of positive rank admitting Gaeta resolutions.

\begin{prop}\label{prop:prioritary} For a divisor $D$, consider the locus
\[
    \{ \phi \in R_f \mid \text{$F_{\phi}$ is torsion-free and not $D$-prioritary} \}.
\]
If $D = A$, then the locus is empty. If $D = A+C$, then the locus is empty if $e=0$ and otherwise has
 codimension $\geqslant \sum_{i \in S_0} \gamma_i + \alpha_4 - r + 1$ in $R_f$.
\end{prop}

\begin{proof} We begin by proving the second statement. Let $F=F_\phi$ and take %$D=C+A - \sum_{i \in I} E_i$
    $D = C+A$. We need to study $\Ext^2(F,F(-D)) \cong \Hom(F(-K-D),F)^\vee$. We twist the Gaeta resolution of $F$ by $\mathcal{O}(-K-D)$ and then apply $\Hom(-,F)$, obtaining the exact sequence
    \begin{align*}
        0\to &\Hom(F(-K-D),F) \\
        \to &H^0(F(C+K+D))^{\alpha_2} \oplus H^0(F(A+K+D))^{\alpha_3}\\
     &\oplus \bigoplus_{i \in S_0} H^0(F(E_{i}+K+D))^{\gamma_{i}} \oplus H^0(F(K+D))^{\alpha_4} \\
      \xrightarrow{h} & H^0(F(C+A+K+D))^{\alpha_1} \oplus \bigoplus_{j \in S_1} H^0(F(C+A+K+D-E_{j}))^{\gamma_{j}}.
    \end{align*}
    A calculation using the Gaeta resolution for $F$ shows that $H^0(F(K+D))$ vanishes as, after tensoring by $K+D$, the line bundles in degree 0 have no global sections and the line bundles in degree $-1$ have no $H^1$:
    \[
        H^1(\oo(-A-C+K+D)) \cong H^1(\oo)^{\vee}
        \cong 0
    \]
and, for $j \in S_1$,
    \[
        H^1(\oo(-A-C+E_j + K + D)) \cong H^1(\oo(-E_j))^{\vee} \cong 0.
    \]
Similarly, for $i \in S_0$, we get $H^0(F(E_i+K+D)) \cong 0$ since $H^1(-\oo(E_i)) \cong 0$ and $H^1(\oo(-E_i - E_j)) \cong 0$.
Thus, the map $h$ in the exact sequence reduces to the map obtained by applying $\Hom(-,F(K+D))$ to the map
    \[
        \oo(-C-A)^{\alpha_1} \oplus \bigoplus_{j \in S_1} \oo(-C-A+E_{j})^{\gamma_{j}} \xrightarrow{\psi} \oo(-C)^{\alpha_2} \oplus \oo(-A)^{\alpha_3}
    \]
from the Gaeta resolution.
If this map is surjective, then $h$ must to be injective. The locus in $\bp(H_{\vec{a},d})$ of $\phi$ such that $\psi$ is not surjective has codimension $\geqslant \alpha_1 + \sum_{j \in S_1} \gamma_j - \alpha_2 - \alpha_3 + 1 = \sum_{i \in S_0} \gamma_i + \alpha_4 -r + 1$, giving the desired estimate. If $e=0$, then $B=C$ and $S_0$ and $S_1$ are empty since $X$ is admissible, from which it follows that $H^0(F(C+K+D)) \cong H^0(F(A+K+D)) \cong 0$.

The first statement follows from a similar argument with $D=A$ by checking that $H^0(F(C+K+D))$, $H^0(F(A+K+D))$, $H^0(F(E_i+K+D))$ for $i \in S_0$, and $H^0(F(K+D))$ all vanish.
\end{proof}

The case $e=0$ (where $X=\mathbb{F}_0$) is in \cite[Proposition 2.20]{Ped21}.

\begin{prop}\label{prop:prioritary-complete}
 Fix exponents such that $r \geqslant 1$, the condition in Proposition~\ref{prop:cokernel-properties}(b.ii) holds, and $\sum_{i \in S_0} \gamma_i + \alpha_4 \geqslant r$. Then the open family of Gaeta resolutions whose cokernels are torsion-free and $(C+A)$-prioritary is a complete family of $(C+A)$-prioritary sheaves.
\end{prop}

By Lemma \ref{lem:prioritary-div}, the same statement holds with $C+A$ replaced by any divisor $D \in \mathcal{D}$, where $\mathcal{D}$ is defined in (\ref{eq:bpf-divisors}). The inequality $\sum_{i \in S_0} \gamma_i + \alpha_4 \geqslant r$ is not needed for the case $D=A$.

\begin{proof} If $F$ has a general Gaeta resolution, then Proposition~\ref{prop:cokernel-properties}(b) implies that $F$ is torsion-free, and Proposition~\ref{prop:prioritary} implies that $F$ is $(C+A)$-prioritary.
Thus, the argument proving \cite[Proposition 3.6]{CosHui18weakBN} applies here: the sequence being strong full exceptional implies the Kodaira-Spencer map is surjective.
\end{proof}

\subsection{Weak Brill-Noether result}
\label{weak-brill-noether}

Let $S$ be $\bp^2$ or $X$ an admissible blowup of $\mathbb{F}_e$. Consider Gaeta resolutions $\E_1^{a_1} \oplus \cdots \oplus \E_d^{a_d} \to \E_{d+1}^{a_{d+1}} \oplus \cdots \oplus \E_n^{a_n}$ on $S$ of the types in Example~\ref{exmp:gaeta}.

For integers $\ell \geqslant 1$ and $r \geqslant 1$, set $a_n = \ell r$ and assume the remaining exponents $a_1,\dots,a_{n-1} \geqslant 0$ satisfy
\begin{equation}\label{eq:WBN-exponents}
    a_1 + \cdots + a_d = a_{d+1} + \cdots + a_{n-1} + (\ell-1)r.
\end{equation}
On $X$, we further assume
\begin{equation}\label{eq:WBN-gammas}
    \gamma_i \geqslant \sum_{j \colon p_j \succ p_i} \gamma_j, \qquad \text{for all $i \in S_0$}.
\end{equation}
Consider a general map
\[
    \phi\colon \E_1^{a_1} \oplus \cdots \oplus \E_d^{a_d} \to \E_{d+1}^{a_{d+1}} \oplus \cdots \oplus \E_{n-1}^{a_{n-1}} \oplus \oo_X^{\ell r}.
\]

By Proposition~\ref{prop:cokernel-properties}, $\phi$ is injective and the cokernel $F_{\phi}$ is torsion-free of rank $r$. Furthermore, it has vanishing higher cohomology, and \[\oo_S^{\ell r} \cong H^0(F_{\phi}) \otimes \oo_S.\]
Then, if $Z$ is a zero-dimensional subscheme of length $\ell$, $\chi(F_{\phi}\otimes I_Z)=0$ and we have the following weak Brill-Noether result.

\begin{prop}\label{orth-ker-quot} Let $S$ denote $\bp^2$ or $X$ an admissible blowup of $\mathbb{F}_e$. Suppose the sequence of exponents $(a_1,\dots,a_n)$ satisfies $a_n = \ell r$ and (\ref{eq:WBN-exponents}) for some $r,\ell \geqslant 1$, as well as (\ref{eq:WBN-gammas}) in the case $S=X$. If $\phi$ is general and $Z \in S^{[\ell]}$ is general, then $F_{\phi} \otimes I_Z$ has vanishing cohomology in all degrees.
\end{prop}

\begin{proof}
As the vanishing of the cohomology of $F_{\phi} \otimes I_Z$ is an open condition on families of $\phi$ and of $Z$, it suffices to prove the claim for a single choice of $\phi$. We construct $\phi$ as the direct sum of $r$ maps $\{ \phi_m \}_{1 \leqslant m \leqslant r}$, each of which is a general map of the form
\[
    \E_1^{a_1'} \oplus \cdots \oplus \E_d^{a_d'} \to \E_{d+1}^{a_{d+1}'} \oplus \cdots \oplus \E_{n-1}^{a_{n-1}'} \oplus \oo_X^\ell,
\]
where the exponents, which depend on $m$, are non-negative and satisfy
\[
    a_1' + \cdots + a_d' = a_{d+1}' + \cdots + a_{n-1}' + \ell - 1.
\]
On $X$, we need the additional condition that the exponent of $\oo(-E_i)$ is at least as large as the sum of the exponents of $\oo(-C-A+E_j)$ for $j \in S_1$ such that $p_j \succ p_i$, which can be ensured due to (\ref{eq:WBN-gammas}).
Then, by Proposition \ref{prop:cokernel-properties}, the cokernel of $\phi$ is of the form $F_{\phi} \cong \bigoplus_{m=1}^r L_m \otimes I_{Z_m'}$,
where, for each $m$, $L_m$ is a line bundle with vanishing higher cohomology, $Z_m'$ is a 0-dimensional subscheme, vanishing on $Z_m'$ imposes independent conditions on $H^0(L_m)$, and $H^0(L_m \otimes I_{Z_m'})$ is $\ell$-dimensional. Choose distinct points $Z = \{ q_1,\dots,q_{\ell} \}$ inductively so that each $q_u$ avoids the base loci of the linear systems of curves in $|L_m|$ that vanish on $Z_m' \sqcup \{ q_1,\dots,q_{u-1} \}$. Then
$F_{\phi} \otimes I_Z \cong \bigoplus_{m=1}^r L_m \otimes I_{Z_m' \sqcup Z}$ has no cohomology.
\end{proof}

\subsection{No sections vanishing on curves}
The following result will be used to apply \cite{GolLin22} in the study of finite Quot schemes in \S~\ref{sect:finite-quot}.

\begin{prop}\label{prop:no-sec-van-curves} Suppose $\phi$ is a general Gaeta resolution
and $F_\phi$ is the cokernel.
\begin{enumerate}[(a)]
    \item On $\mathbb{P}^2$ or $\mathbb{F}_e$, $F_\phi$ has no sections vanishing on curves.
    \item On $X$ an admissible blowup of $\mathbb{F}_e$, assume \begin{enumerate}[(i)]
    \item $\gamma_j \geqslant \alpha_4-r$ for all $j \in S_1$,
    \item $\gamma_i \geqslant \sum_{j \in S_1 \colon p_j \succ p_i} \gamma_j + \max(0,\alpha_4-r)$ for all $i \in S_0$.
\end{enumerate}
Then $F_\phi$ has no sections vanishing on curves.
\end{enumerate}
\end{prop}

\begin{proof}
Let $F=F_\phi$. It suffices to prove that $H^0(F(-D))=0$ for all minimal nonzero effective divisors $D$. Twisting the Gaeta resolution by $-D$ and taking cohomology, we get a long exact sequence
\[
    0 \to H^0(F(-D)) \to H^1(\E_1(-D)^{a_1} \oplus \cdots \oplus \E_d(-D)^{a_d}) \xrightarrow{\tilde{\phi}} H^1(\E_{d+1}(-D)^{a_{d+1}} \oplus \cdots \oplus \E_{n}(-D)^{a_{n}}).
\]
We need to show that the induced map $\tilde{\phi}$ is injective when $\phi$ is general.

(a) On $\bp^2$, the only minimal effective divisor is the hyperplane class $H$, and $H^1(\oo(-3))$
vanishes, so injectivity of $\tilde{\phi}$ is trivial. Similarly, on $\mathbb{F}_e$, the minimal effective divisors are $A$ and $B$, and an easy calculation using \S~\ref{subsect:coh-line-bdl} shows that $\oo(-C-2A)$
and $\oo(-C-A-B)$ both have vanishing $H^1$.

(b) Now consider the sequences on the two-step blowup of $\mathbb{F}_e$. We argue by cases.

\underline{\emph{Case 1:}} $D$ is minimal nonzero effective not equal to any $\tilde{E}_i$ for $i \in S_0$ or $E_j$ for $j \in S_1$. In this case, as above, we show that the domain of $\tilde{\phi}$ vanishes. Fix a curve $Y$ in the linear equivalence class of $D$. As $D$ is minimal effective, $Y$ must be connected. As $C+A$ is ample on $\mathbb{F}_e$ and $D$ is not contained in $\bigcup E_i$, there is a curve $Y'$ in the linear series $|C+A|$ that is connected,
does not contain $Y$, and intersects $Y$. Then $Y \cup Y'$ is a connected curve in the linear series $|C+A+D|$, so $H^1(\oo(-C-A-D)) = 0$ by Lemma~\ref{lem:conn-curve}.

We use a similar argument to show the vanishing of $H^1(\oo(-C-A-E_j-D))$ for $j \in S_1$.
Let $i \in S_0$ be the index such that $p_j \succ p_i$. By Lemma~\ref{lem:LS_base_locus},
\[
    |C+A-E_j| = |C+A-E_i| + (E_i-E_j),
\]
where $|C+A-E_i|$ is basepoint-free (Example~\ref{ex:LS_section}) and contains connected curves that intersect $\tilde{E}_i$, so the union of such a curve and $(E_i-E_j)$ is connected. Thus, for $D$ minimal effective not equal to any $\tilde{E}_i$ or $E_j$, $H^1(\oo(-C-A+E_j-D))=0$ (and this vanishing holds for $D=E_j$ as well).

\underline{\emph{Case 2:}} $D = E_j$ for $j \in S_1$. Note that if $L$ is a line bundle with no cohomology whose restriction to $E_j$ is trivial, then taking cohomology of the short exact sequence
\[
    0 \to L(-E_j) \to L \to L|_{E_j} \to 0
\]
yields an isomorphism $H^1(L(-E_j)) \cong H^0(\oo_{E_j})$.
Moreover, because $H^1(\oo(-E_j))=0$ and $H^1(\oo(-C-A+E_j-E_j))=0$, we can view $\tilde{\phi}$ as the map obtained from
\[
    \oo(-C-A)^{\alpha_1} \oplus \bigoplus_{j' \in S_1 \setminus \{j\}} \oo(-C-A+E_{j'})^{\gamma_{j'}} \to \oo(-C)^{\alpha_2} \oplus \oo(-A)^{\alpha_3} \oplus \bigoplus_{i \in S_0} \oo(-E_i)^{\gamma_i}
\]
by restricting to $E_j$ and then taking the induced map on global sections. As the restriction to $E_j$ of each of these exceptional sheaves is trivial, $\tilde{\phi}$ is of the form
\[
    \tilde{\phi} \colon \basefield^{\alpha_1} \oplus \bigoplus_{j' \in S_1 \setminus \{j\}} \basefield^{\gamma_{j'}} \to \basefield^{\alpha_2} \oplus \basefield^{\alpha_3} \oplus \bigoplus_{i \in S_0} \basefield^{\gamma_i}.
\]
Thus, a necessary condition for $\tilde{\phi}$ to be injective is $\alpha_1 - \gamma_j + \sum_{j' \in S_1} \gamma_{j'} \leqslant \alpha_2 + \alpha_3 + \sum_{i \in S_0} \gamma_i$,
or equivalently, as $\alpha_1 + \sum_{j' \in S_1} \gamma_{j'} + r = \alpha_2 + \alpha_3 + \sum_{i \in S_0} \gamma_i + \alpha_4$,
\[
    \gamma_j \geqslant \alpha_4 - r.
\]

To find additional sufficient conditions for $\tilde{\phi}$ to be injective, we observe that certain blocks of $\tilde{\phi}$ are injective when $\phi$ is general.
Let $i$ denote the index such that $p_j \succ p_i$. Then:
\begin{itemize}%[(a)]
    \item The block $\bigoplus_{j' \ne j \colon p_{j'} \succ p_i} \basefield^{\gamma_{j'}} \to \basefield^{\gamma_i}$ is injective for general $\phi$ if
\[
    \sum_{j' \ne j \colon p_{j'} \succ p_i} \gamma_{j'} \leqslant \gamma_i
\]
because the linear series $|C+A-E_i-E_{j'}|$ is basepoint-free.
In fact, this inequality is necessary because $E_j$ is in the base locus of $|A-E_{j'}|$, $|eA+B-E_{j'}|$, and $|C+A-E_{i'}-E_{j'}|$ for all $j' \ne j$ such that $p_{j'} \succ p_i$ and all $i' \in S_0 \setminus \{i\}$, hence all other blocks involving these $k^{\gamma_{j'}}$ are 0.

\item For a similar reason, for each $i' \in S_0 \setminus \{i\}$, $\bigoplus_{j' \colon p_{j'} \succ p_{i'} %\ne p_i
    } \basefield^{\gamma_{j'}} \to \basefield^{\gamma_{i'}}$ is injective for general $\phi$ if
\[
    \sum_{j' \colon p_{j'} \succ p_{i'}} \gamma_{j'} \leqslant \gamma_{i'}
\]
(though this inequality may not be necessary for $\tilde{\phi}$ to be injective).
\end{itemize}
As all blocks corresponding to $\basefield^{\alpha_1}$ are general if $\phi$ is general because the linear systems $|A|$, $|C|$, and $|C+A-E_{i'}|$ for $i' \in S_0$ are all basepoint-free, these inequalities suffice to ensure that $\tilde{\phi}$ is injective when $\phi$ is general.

\underline{\emph{Case 3:}} $D = \tilde{E}_i$ for $i \in S_0$. Similar to the previous case, if $L$ is a line bundle on $X$ with no cohomology whose restriction to $\tilde{E}_i$ is trivial, then
\[
    H^1(L(-\tilde{E}_i)) \cong H^0(L|_{\tilde{E}_i}),
\]
while $H^1(\oo(-\tilde{E}_i))=0$. Thus, $\tilde{\phi}$ can be viewed as the map obtained by restricting
\begin{align*}
    \oo(-C-A)^{\alpha_1} \oplus \bigoplus_{j \in S_1} &\oo(-C-A+E_{j})^{\gamma_{j}} \tag{$*$} \\
    & \to \oo(-C)^{\alpha_2} \oplus \oo(-A)^{\alpha_3} \oplus \bigoplus_{i' \in S_0 \setminus \{i\}} \oo(-E_{i'})^{\gamma_{i'}} \oplus \oo(-E_i)^{\gamma_i}
\end{align*}
to $\tilde{E}_i$ and then taking the induced map on global sections. The restriction of each of the exceptional bundles in $(*)$ to $\tilde{E}_i$ is $\oo_{\tilde{E}_i}$, except for $\oo(-E_i)|_{\tilde{E}_i} \cong \oo_{\tilde{E}_i}(1)$
and $\oo(-C-A+E_j)|_{\tilde{E}_i} \cong \oo_{\tilde{E}_i}(1)$ for $p_j \succ p_i$,
each of which has a two-dimensional space of global sections. Thus, the restriction of ($*$) to $\tilde{E}_i$ is
\[
    \oo_{\tilde{E}_i}^{\alpha_1} \oplus \bigoplus_{j \colon p_j \not \succ p_i} \oo_{\tilde{E}_i}^{\gamma_j} \oplus \bigoplus_{j \colon p_j \succ p_i} \oo_{\tilde{E}_i}(1)^{\gamma_j} \to \oo_{\tilde{E}_i}^{\alpha_2} \oplus \oo_{\tilde{E}_i}^{\alpha_3} \oplus \bigoplus_{i' \in S_0 \setminus \{i\}} \oo_{\tilde{E}_i}^{\gamma_{i'}} \oplus \oo_{\tilde{E}_i}(1)^{\gamma_i}
\]
and hence $\tilde{\phi}$ is of the form
\[
    \tilde{\phi} \colon \basefield^{\alpha_1} \oplus \bigoplus_{j \colon p_j \not \succ p_i} \basefield^{\gamma_j} \oplus \bigoplus_{j \colon p_j \succ p_i} H^0(\oo_{\tilde{E}_i}(1))^{\gamma_j} \to \basefield^{\alpha_2} \oplus \basefield^{\alpha_3} \oplus \bigoplus_{i' \in S_0 \setminus \{i\}} \basefield^{\gamma_{i'}} \oplus H^0(\oo_{\tilde{E}_i}(1))^{\gamma_i}.
\]
A necessary inequality for $\tilde{\phi}$ to be injective is $\alpha_1 + \sum_{j \in S_1} \gamma_j + \sum_{j \in S_1 \colon p_j \succ p_i} \gamma_j \leqslant \alpha_2 + \alpha_3 + \sum_{i' \in S_0} \gamma_{i'} + \gamma_i$, or equivalently,
\[
    \gamma_i \geqslant \alpha_4 - r + \sum_{j \colon p_j \succ p_i} \gamma_j.
\]
As in the previous case, we obtain sufficient conditions for $\tilde{\phi}$ to be injective when $\phi$ is general by looking at various blocks.

\begin{itemize}%[(a)]
    \item The block $\bigoplus_{j \colon p_j \succ p_i} H^0(\oo_{\tilde{E}_i}(1))^{\gamma_j} \to H^0(\oo_{\tilde{E}_i}(1))^{\gamma_i}$
is injective for general $\phi$ if
\[
    \gamma_i \geqslant \sum_{j \colon p_j \succ p_i} \gamma_j
\]
because the linear system $|C+A-E_i-E_j|$ is basepoint-free.
As there are no nonzero maps $\oo_{\tilde{E}_i}(1) \to \oo_{\tilde{E}_i}$, this condition is necessary for $\tilde{\phi}$ to be injective as all other blocks involving $\bigoplus_{j \colon p_j \succ p_i} H^0(\oo_{\tilde{E}_i}(1))^{\gamma_j}$ are 0.

\item For $i' \in S_0 \setminus \{i\}$, the block $\bigoplus_{j \colon p_j \succ p_{i'}} \basefield^{\gamma_j} \to \basefield^{\gamma_{i'}}$ is injective for general $\phi$ if
\[
    \gamma_{i'} \geqslant \sum_{j \colon p_j \succ p_{i'}} \gamma_j
\]
(but this condition may not be necessary).
\end{itemize}
The blocks involving $\basefield^{\alpha_1}$ are all general for general $\phi$ as the linear systems $|A|$, $|C|$, $|C+A-E_{i'}|$ for $i' \ne i$ are all basepoint-free and the curves in $|C+A|$ containing $p_i$ have no fixed tangent direction at $p_i$, so the above inequalities are sufficient.
\end{proof}

\begin{lem}\label{lem:conn-curve} Let $D$ be a nonzero effective divisor on a rational surface and $|D|$ denote its linear series. The following are equivalent:
\begin{enumerate}[(a)]
    \item $|D|$ contains a connected curve;
    \item Every curve in $|D|$ is connected;
    \item $H^1(\oo(-D))=0$.
\end{enumerate}
\end{lem}

\begin{proof}
Let $Y$ be a curve in the linear equivalence class of $D$ and consider the corresponding short exact sequence
\[
    0 \to \oo(-D) \to \oo \to \oo_Y \to 0.
\]
Taking cohomology and using $H^1(\oo)=0$, we see that $H^1(\oo(-D))=0$ if and only if $H^0(\oo) \to H^0(\oo_Y)$ is an isomorphism, which is true if and only if $H^0(\oo_Y)$ is one-dimensional, which holds exactly when $Y$ is connected. (If $D$ is ample, $(c)$ also follows from Kodaira vanishing.)
\end{proof}

\section{Gaeta resolutions and stability}\label{sec:stability}

For $X$ an admissible blowup, we study the connection between the existence of Gaeta resolutions and stability of a sheaf, which allows Gaeta resolutions to be applied in the study of moduli problems. First, we describe conditions ensuring that general stable sheaves in $M(f)$ admit Gaeta resolutions (Proposition~\ref{prop:unirat-large-L}). Then, by imposing stronger conditions on $f$ and on the polarization $H$, we show the locus of maps in the resolution space whose cokernels are unstable has codimension $\geqslant 2$ (Proposition~\ref{prop:bad-locus-res}), as well as the parallel statement in the moduli space that the locus of sheaves not admitting Gaeta resolutions has codimension $\geqslant 2$ (Theorem~\ref{thm:moduli-nongr-codim2}). In particular, the latter results imply that $M(f)$ is non-empty and that general stable sheaves away from a locus of codimension $\geqslant 2$ satisfy various nice properties (Corollary~\ref{cor:gen-stable-prop}).

We have discussed about prioritary conditions in the previous section.
One of the motivations to consider the prioritary condition is the following statement, which is essentially in the proof of \cite[Theorem 1]{Wal98}.
\begin{lem}\label{lem:stable-is-priotary}
Over a smooth projective surface $S$, if $D$ is a divisor and $H$ is a polarization such that $H\cdot (K_S+D)<0$, then any $H$-semistable torsion-free sheaf is $D$-prioritary.
\end{lem}

\begin{proof}
If $F$ is torsion-free and $H$-semistable, then by Serre duality, $\Ext^2(F,F(-D)) \cong \Hom(F,F(K_S+D)) \cong 0$ by semistability and the fact that $H \cdot (K_S+D)<0$.
\end{proof}

As a warm-up, we prove Proposition~\ref{prop:unirat-large-L}.
Recall that $X$ is an admissible blowup of $\mathbb{F}_e$, $H$ satisfies $H \cdot (K_X+A) < 0$, and $f \in K(X)$ is a class of rank $r > 0$, Euler characteristic $\chi \geqslant 0$, and first Chern class satisfying the inequalities in Propositions \ref{prop:chern-has-gr}(a) and \ref{prop:cokernel-properties}(b.ii).
\begin{proof}[Proof of Proposition \ref{prop:unirat-large-L}]
When non-empty, the moduli space $M(f)$ is irreducible (\cite[Theorem 1]{Wal98}), using the fact that the stack of prioritary sheaves is irreducible. Since semistability is an open property in families, Propositions~\ref{prop:chern-has-gr}(a), \ref{prop:prioritary-complete}, and Lemma~\ref{lem:stable-is-priotary} show that a general semistable sheaf in $M(f)$ has a Gaeta resolution.
\end{proof}

Using Proposition~\ref{prop:chern-has-gr}(b), we obtain a formulation similar to Proposition \ref{prop:unirat-large-L}.

\begin{prop}\label{prop:unirat-large-c2}
 Let $X$ and $H$ be as in Proposition~\ref{prop:unirat-large-L}.
 Assume the class $f$ has fixed rank $r>0$ and first Chern class, and its discriminant is sufficiently large.
 Then there is a line bundle $L$ pulled back from $\mathbb{F}_e$ such that for %all
 general semistable sheaves $F$ of class $f$, $F \otimes L$ admits a Gaeta resolution.
\end{prop}

Since the space $R_f$ of Gaeta resolutions is rational, we immediately deduce the following special cases of \cite[Theorem 2.2]{Bal87}.

\begin{cor}\label{cor:unirational} For $X$, $H$, and $f \in K(X)$ as in Propositions \ref{prop:unirat-large-L} or \ref{prop:unirat-large-c2}, $M(f)$ is unirational if it is nonempty.
\end{cor}

Over Hirzebruch surfaces, similar results were proved in \cite[Theorem 2.10, Proposition 4.4]{CosHui20}.

The rest of this section is dedicated to the proof of Theorem~\ref{thm:moduli-nongr-codim2},
which requires additional technical conditions on the exponents of the Gaeta resolutions and on the polarization. Before stating these conditions, we set some notation. For simplicity, we write the Gaeta resolutions of the form (\ref{gtr-blowup}) on $X$ as
\[
    0 \to \E_1^{a_1} \oplus \cdots \oplus \E_d^{a_d} \xrightarrow{\phi} \E_{d+1}^{a_{d+1}} \oplus \cdots \oplus \E_{d+1}^{a_{d+1}} \to F_{\phi} \to 0.
\]
Let $\rcone \subset H_{\vec{a},d}$ denote the locus of injective maps $\phi$ as in \S~\ref{sect:prop-gr}.
Let $p$ and $q$ denote the projections from $\rcone \times X$ to the first and second factors, respectively.
Let $\mathbb{F}$ denote the cokernel of the tautological map over $\rcone \times X$:
\begin{align}\label{eq:univ-res}
        0 \to \bigoplus_{i=1}^d q^*\E_i^{a_i}  \to  \bigoplus_{j=d+1}^n q^*\E_j^{a_j} \to \mathbb{F}\to 0.
    \end{align}
Because it is the family of cokernels of injective maps, $\mathbb{F}$ is flat over $\rcone$.

The technical conditions on the exponents are as follows. We assume
\begin{equation}\label{eq:conditions-exponents}
r \geqslant 2 \quad \text{and} \quad \gamma_i \geqslant \sum_{j \colon p_j \succ p_i} \gamma_j + 1 \; \quad \forall i \in S_0
\end{equation}
to ensure that the complement of $\rcone$ in $H_{\vec{a},d}$ has codimension $\geqslant 2$ and that the cokernels $F_\phi$ are torsion-free away from a locus of codimension $\geqslant 2$ in $\rcone$ (the proof is similar to the proof of Proposition \ref{prop:cokernel-properties}, but we consider the image of the locus $Z$ in Proposition \ref{prop:Bertini-type} in $\bp \Hom(\mathcal{A},\mathcal{B})$). We make the additional assumption
\begin{equation}\label{eq:condition-alpha4}
    \sum_{i \in S_0} \gamma_i + \alpha_4 \geqslant r+1,
\end{equation}
which ensures that the $F_\phi$ are $(C+A)$-prioritary away from a locus of codimension $\geqslant 2$ in $\rcone$ (Proposition \ref{prop:prioritary}), require all of the exponents $a_i$ to be strictly positive, and require the discriminant of the sheaves $F_\phi$ to be sufficiently large in the sense of (\ref{eq:discriminant-bound}).

To state the conditions on the ample divisor $H$, which we assume is general, we write %$H$ as
\[H=uA+vC-\sum_{i \in S_0 \cup S_1} d_i E_i \]
for rational numbers $u,v,d_i > 0$ \footnote{As scaling $H$ does not affect stability, these weights could be taken as integers as well.}.
We assume that
\begin{equation}\label{eq:H-pos-lin-combo}
    u > \sum_{j \in S_1} d_j \mbox{ and } v > \sum_{i \in S_0} d_i > \sum_{j \in S_1} d_j,
\end{equation}
namely that $H$ is a positive linear combination of all divisors in $\mathcal{D}$, where $\mathcal{D}$ is defined in (\ref{eq:bpf-divisors}). Note that by Proposition \ref{prop:very-ample}, (\ref{eq:H-pos-lin-combo}) implies that $H$ is ample.
This condition implies that
\[H\cdot (K_X+D)<0,\ \forall D\in \mathcal{D}.\]
as well as the condition $H \cdot (K_X+2A) < 0$ that appears in Theorem \ref{thm:yoshioka}.
Moreover, we assume that no $d_i$ can be too close to $v$ in the sense that
\begin{equation}\label{eq:conditions-on-H}
\frac{d_i}{v} \leqslant \frac{\lambda}{\lambda+1} \quad \text{and}
 \quad \frac{d_i}{v} \leqslant \sqrt{\frac{2\lambda}{t}} \quad \text{for all $i$, where $\lambda = \frac{u}{v} + \frac{e}{2}$ \\
}.
\end{equation}

\begin{rem} Our arguments can be extended to allow $\mathcal{D}$ to include other divisors whose linear systems are basepoint-free with general member isomorphic to $\bp^1$, for instance other divisors of the form $C+A- \sum_{i \in I} E_i$ for $I \subset \{1,\dots,t\}$ such that $\# I \leqslant e+2$, which could expand the range of $H$ depending on the configuration of the blown-up points.
\end{rem}

\subsection{Locus of unstable sheaves in the space of Gaeta resolutions.}\label{subsect:bad-locus-res}

In this subsection we prove the following proposition.

\begin{prop}\label{prop:bad-locus-res}
Assume that $\vec{a}$ and $H$ satisfy the conditions (\ref{eq:conditions-exponents}-\ref{eq:conditions-on-H}). %(\ref{eq:H-pos-lin-combo}) and (\ref{eq:conditions-on-H}).
Then the complement of $\{\phi\in \rcone \mid F_\phi% =\coker \phi
\mbox{ is $H$-semistable}\}$ in $H_{\vec{a},d}$ has codimension $\geqslant 2$. In particular, the moduli space $M(f)$ is non-empty.
\end{prop}

As observed above, the complement of $\rcone$
has codimension at least $2$ in $H_{\vec{a},d}$,
as does the locus of $\phi$ such that $F_\phi$ is not torsion-free. Moreover, by Proposition \ref{prop:prioritary}, the locus of $\phi$ such that $F_{\phi}$ is not $(C+A)$-prioritary also has codimension $\geqslant 2$ in $H_{\vec{a},d}$,
so it suffices to show that the locus of unstable torsion-free sheaves that are $(C+A)$-prioritary
has codimension at least 2. We will do this by showing various Harder-Narasimhan strata have codimension $\geqslant 2$.

For $\phi\in \rcone$ whose cokernel $F_\phi$ is $(C+A)$-prioritary
and unstable, let
\[0=G_0\subsetneqq G_1\subsetneqq \cdots \subsetneqq G_\ell =F_\phi\]
be the Harder-Narasimhan filtration of $F_\phi$ with respect to Gieseker stability, let $\gr_i=G_i/G_{i-1}$ denote the grading and let $r_i$, $\nu_i$, and $\Delta_i$ denote the rank, total slope, and discriminant of $\gr_i$. Then $\nu_1 \cdot H \geqslant \cdots \geqslant \nu_\ell \cdot H$ and $r_i \geqslant 1$ and $\Delta_i \geqslant 0$ for all $i$.

\begin{lem}\label{large-hn}
The following subset of $\rcone$ has codimension at least $2$:
\begin{equation*}
\left\{
\phi\in \rcone \, \middle\vert \begin{array}{l}
     F_\phi \mbox{ is %locally free,
     unstable and $(C+A)$-prioritary};\\
    \mbox{$(\nu_i-\nu_j)\cdot D>2$ for some $i<j$ and some $D \in \mathcal{D} \cup \{C+A\}$}
\end{array}
\right\}.\end{equation*}
\end{lem}
\begin{proof}
All $D \in \mathcal{D} \cup \{C+A\}$ are basepoint-free, so Bertini's Theorem implies that general divisors in the corresponding linear systems are nonsingular. Moreover, these general divisors are isomorphic to $\mathbb{P}^1$. Let $D$ be general in the linear system $|D|$ so that it avoids the singularities of $F_{\phi}$.
The restrictions of the cokernels in a neighborhood of $\phi$, which are locally free on $D$.
As $(C+A)$-prioritary implies $D$-prioritary, $F_\phi$ is $D$-prioritary, so the Kodaira-Spencer map
\[\operatorname{T}_\phi \rcone\to \Ext_X^1(F_\phi,F_\phi)\to \Ext^1_D(F_\phi|_D,F_\phi|_D)\]
is surjective, according to \cite[Corollary 15.4.4]{LeP97} or \cite[Proposition 2.6]{CosHui20}.  The restrictions to $D$ provide a complete family of vector bundles. % parametrized by $\rcone$.
On the other hand, the inequality implies that over $D$, $\mu_{\max{}} (F_\phi|_D)-\mu_{\min{}} (F_\phi|_D)>2$.
Thus, the subset has codimension at least 2, according to \cite[Corollary 15.4.3]{LeP97}.
\end{proof}

The last step is to show that the locus of $\phi$ in $\rcone$ satisfying the following two conditions has codimension $\geqslant 2$:
\begin{enumerate}[(a)]
    \item $F_\phi$ is unstable and $(C+A)$-prioritary;
    \item For all $i<j$, the inequality $(\nu_i-\nu_j)\cdot D\leqslant 2$ holds for all $D \in \mathcal{D} \cup \{C+A\}$.
\end{enumerate}
We will do this for each locally-closed stratum of this locus of fixed Harder-Narasimhan type, namely for an integer $\ell \geqslant 2$ and polynomials $P_1,\dots,P_\ell$, we let $Y_{P_1,\dots,P_\ell}$ denote the locus of $\phi$ in $\rcone$ such that (a) and (b) hold and the Harder-Narasimhan filtration of $F_{\phi}$ has length $\ell$ and the Hilbert polynomial of $\gr_i$ is $P_i$. Our strategy for showing the codimension of $Y_{P_1,\dots,P_\ell}$ in $\rcone$ is $\geqslant 2$ is based on similar ideas for $\mathbb{P}^2$ in \cite[Chapter 15]{LP05}.

We begin with the following observation:
\begin{lem}\label{lem:gr-ext}
\begin{enumerate}[(i)]
    \item For $i<j$, $\Hom(\gr_i,\gr_j)=0$.
    \item Under the conditions (a) and (b) above,
 $\Ext^2(\gr_i,\gr_j)=0$ for all $i$ and $j$.
\end{enumerate}
\end{lem}

\begin{proof}Part (i) follows from semistability.

For (ii), (b) implies that when $i<j$, $(\nu_i - \nu_j + K) \cdot D \leqslant 0$ for all $D \in \mathcal{D}$ and $(\nu_i-\nu_j+K) \cdot (C+A) \leqslant -e-2 < 0$.
As $H$ is a positive linear combination of all divisors in $\mathcal{D}$, letting $m$ denote the minimum of the weights of $A$ and $C$, $H$ can be written as $m(C+A)$ plus a non-negative linear combination of divisors in $\mathcal{D}$. Thus, $(\nu_i - \nu_j + K)\cdot H<0$ and this inequality is also true when $i \geqslant j$ since (\ref{eq:H-pos-lin-combo}) implies $H \cdot K < 0$. Thus, $\Ext^2(\gr_i,\gr_j)\cong\Hom(\gr_j,\gr_i\otimes K)^\vee=0 $ for all $i$ and $j$.
\end{proof}

Let $\operatorname{Flag} = \operatorname{Flag}(\mathbb{F}/\rcone;P_1,\dots,P_\ell) \to \rcone$ be the relative flag scheme of filtrations whose grading $\gr_i$ has Hilbert polynomial $P_i$. Given $\phi \in Y_{P_1,\dots,P_\ell}$ and a point $p=(\phi,(G_1,\dots,G_\ell))$ of the fiber over $\phi$, there is an exact sequence
\begin{equation}\label{sequence-tangent}
    0\to \Ext^0_+(F_\phi,F_\phi)\to \operatorname{T}_p\operatorname{Flag}\to \operatorname{T}_\phi \rcone\xrightarrow[]{\omega_+}\Ext^1_+(F_\phi,F_\phi),
\end{equation}
\cite[Proposition 15.4.1]{LeP97} realizing the vertical tangent space as the group $\Ext^0_+(F_\phi,F_\phi)$ and the normal space of $Y_{P_1,\dots,P_\ell}$ in $\rcone$ at $\phi$ as the image of $\omega_{+}$.
Here the groups $\Ext^i_+$ are defined with respect to the filtration of $F_{\phi}$, and $\omega_+$ is the composite map
\begin{equation}\label{eq:omega-plus}\operatorname{T}_\phi \rcone \xrightarrow{\operatorname{KS}} \Ext^1(F_\phi,F_\phi)\xrightarrow{h_+} \Ext^1_+(F_\phi,F_\phi).\end{equation}
The Kodaira-Spencer map $\operatorname{KS}$ is surjective since $\rcone$ parametrizes a complete family of prioritary sheaves. The idea is to show that $h_+$ is also surjective and that $\ext^1_+(F_\phi,F_\phi) \geqslant 2$, which imply that $Y_{P_1,\dots,P_\ell}$ has codimension $\geqslant 2$ in $\rcone$.
For foundational material on the groups $\Ext^i_\pm$, see \cite{DreLeP-85}.

There is a canonical exact sequence
\[\dots \to \Ext^{1}(F_\phi,F_\phi)\xrightarrow{h_+} \Ext^{1}_+(F_\phi,F_\phi)\to \Ext^{2}_-(F_\phi,F_\phi)\to \dots \]
showing that surjectivity of $h_+$ is guaranteed by the vanishing of $\Ext^2_-(F_\phi,F_\phi)$. This group can be calculated using the spectral sequence
\begin{equation*}
    E_1^{p,q}=\begin{cases}
    0 & \text{if $p<0$,} \\
   \bigoplus_i\Ext^{p+q}(\gr_i, \gr_{i-p} ) & \text{otherwise,}
    \end{cases}
\end{equation*}
converging to $\Ext^{p+q}_-(F_\phi,F_\phi)$. By Lemma~\ref{lem:gr-ext} (b), $E_1^{p,q}$ vanishes when $p+q=2$, hence $\Ext_-^2(F_\phi,F_\phi)=0$, so $h_+$ is surjective.

To calculate $\Ext_+^1(F_\phi,F_\phi)$, we use the spectral sequence
\begin{equation*}
    E_1^{p,q}=\begin{cases}
    \bigoplus_i\Ext^{p+q}(\gr_i, \gr_{i-p} ) & \text{if $p<0$,} \\
   0 & \text{otherwise,}
    \end{cases}
\end{equation*}
converging to $\Ext^{p+q}_+(F_\phi,F_\phi)$.
By Lemma~\ref{lem:gr-ext}, $\Ext_+^0(F_\phi,F_\phi)= \Ext_+^2(F_\phi,F_\phi)=0$ and the spectral sequence degenerates on the first page, yielding
\begin{equation*}
    \Ext^1_+(F_\phi,F_\phi)\cong \bigoplus_{i< j}\Ext^1(\gr_i,\gr_j).
\end{equation*}
Using (\ref{eq:euler-pair}), we thus calculate
\begin{equation}\label{eq:ext^1_+}
\ext_+^1(F_\phi,F_\phi)=\sum_{i<j}\ext^1(\gr_i,\gr_j)=-\sum_{i<j}\chi(\gr_i,\gr_j) = \sum_{i<j} r_i r_j(\Delta_i + \Delta_j - P(\nu_j - \nu_i)),
\end{equation}
where $P$ is the Hilbert polynomial of $\mathcal{O}X$, as in \S~\ref{ss:chern-characters}.

To finish the proof of the proposition, we will show that
$\ext^1_+(F_\phi,F_\phi) \geqslant 2$ given the conditions $(\nu_i-\nu_j) \cdot D \leqslant 2$ for all $D \in \mathcal{D} \cup \{C+A\}$ and $(\nu_i-\nu_j) \cdot H \geqslant 0$ for all $i<j$.
As $\chi(\gr_i,\gr_j) \leqslant 0$ for each $i<j$, we see that
\[
    -\sum_{i < j} \chi(\gr_i,\gr_j) \geqslant -\chi(\gr_1,\oplus_{1 < j} \gr_j),
\]
which allows us to reduce to the case where the Harder-Narasimhan filtration has length 2. To see that the conditions (a) and (b) still hold, note that $\nu(\oplus_{1 < j} \gr_j) = \sum_{1 < j} r_j \nu_j/\sum_{1 < j} r_j$
is a weighted average of the $\nu_j$, hence we have $(\nu_1 - \nu(\oplus_{1 < j} \gr_j)) \cdot H \geqslant 0$ and $(\nu_1 - \nu(\oplus_{1 < j} \gr_j)) \cdot D \leqslant 2$ for all $D \in \mathcal{D} \cup \{C+A\}$. As we know $\Delta_1 \geqslant 0$ but no such inequality is guaranteed for the discriminant of $\oplus_{1 < j} \gr_j$, we need $\sum_{1 < j} r_j \geqslant r_1$ to apply the lemma below; if this does not hold, then a similar setup using $\oplus_{i < \ell} \gr_i$ and $\gr_\ell$, which satisfies $\sum_{i < \ell} r_i > r_\ell$ and $\Delta_\ell \geqslant 0$, meets the conditions of the lemma.
Thus, it suffices to prove the following:

\begin{lem} Assume the condition (\ref{eq:conditions-on-H}) on $H$. Suppose the class $(r,\nu,\Delta)$ is the sum of classes $(r_1,\nu_1,\Delta_1)$ and $(r_2,\nu_2,\Delta_2)$ of positive rank with the property that $(\nu_1 - \nu_2) \cdot D \leqslant 2$ for $D = A,C,A+C$, that $(\nu_1 - \nu_2) \cdot H \geqslant 0$, that $\Delta_1 \geqslant 0$ if $r_2 \geqslant r_1$, and that $\Delta_2 \geqslant 0$ if $r_1 \geqslant r_2$.
Then the condition
\begin{equation}\label{eq:discriminant-bound}
\Delta \geqslant \frac{(\lambda+1)^2}{4\lambda} + \frac{t}{8} + \frac{1}{r}, \qquad \text{where $\lambda = \frac{u}{v} + \frac{e}{2}$},
\end{equation}
is sufficient to ensure that $r_1 r_2(\Delta_1 + \Delta_2 - P(\nu_2 - \nu_1)) \geqslant 2$.
\end{lem}

\begin{proof}
Note that $r_1 + r_2 = r$, $r_1 \nu_1 + r_2 \nu_2 = r \nu$, and $r_1 \Delta_1 + r_2 \Delta_2 = r \Delta + \frac{r_1 r_2}{2r} (\nu_2 - \nu_1)^2$.
Using this, in the case $r_2 \geqslant r_1$, we write
\begin{align*}
    r_1 r_2(\Delta_1 + \Delta_2 - P(\nu_2 - \nu_1))
    &= r_1(r_1 \Delta_1 + r_2 \Delta_2) + (r_2-r_1)r_1 \Delta_1 - r_1 r_2 P(\nu_2 - \nu_1) \\
    & = (r_2-r_1)r_1 \Delta_1 + r_1 r \Delta - r_1 r_2(\tfrac{r_2}{2r} (\nu_2 - \nu_1)^2 - \tfrac{1}{2}(\nu_2-\nu_1) \cdot K + 1).
\end{align*}
As $\Delta_1 \geqslant 0$, it suffices to find an upper bound for $\frac{r_2}{2r} (\nu_2 - \nu_1)^2 - \frac{1}{2}(\nu_2 - \nu_1) \cdot K + 1$.
Setting $\xi = r/r_2$ and writing \[\nu_1 - \nu_2 = aA + bB - \sum_i e_i E_i,\] this equals
\[
    \xi^{-1} \left( (a-\xi - \frac{e}{2}b)(b-\xi) - \frac{1}{2} \sum_i e_i(e_i-\xi)  \right) + 1 - \xi.
\]
The lemma below shows that given (\ref{eq:conditions-on-H}), an upper bound is
\[
    \xi^{-1} \left( \frac{\xi^2(\lambda+1)^2}{4\lambda} + \frac{t \xi^2}{8} \right) + 1 - \xi,
\]
where $\lambda = \tfrac{u}{v} + \tfrac{e}{2} > 0$ since if $e=0$ then $H$ ample ensures that $u,v>0$.
This yields the bound
\[
    r_1 r_2(\Delta_1 + \Delta_2 - P(\nu_2 - \nu_1)) \geqslant r_1 r \Delta - r_1 r_2 \left( \frac{\xi (\lambda+1)^2}{4\lambda} + \frac{t \xi}{8}  + 1 - \xi \right),
\]
and, as $r_2 \xi = r$, we can guarantee that the right side is $\geqslant 2$ by assuming $\Delta$ satisfies (\ref{eq:discriminant-bound}).

The argument in the case $r_1 \geqslant r_2$ is similar.
\end{proof}

Before stating and proving the lemma below, we introduce some useful notation. Set $\lambda = \frac{u}{v} + \frac{e}{2}$ as above and consider the change of variables
\[
    J = a + \frac{u}{v}b.
\]
The conditions $(\nu_1 - \nu_2) \cdot D \leqslant 2$ can be written as $a \leqslant 2$, $b \leqslant 2$, $a+b \leqslant 2$, and the condition $(\nu_1 - \nu_2) \cdot H \geqslant 0$ is $J \geqslant \sum_i \tfrac{d_i}{v} e_i$. Moreover, thinking of $J$ as fixed,
\[
    (a-\xi-\frac{e}{2}b)(b-\xi) = -\lambda b^2 + (J+\xi(\lambda-1))b - \xi(J-\xi)
\]
is a quadratic function with maximum value
\begin{equation}\label{eq:parabola-max}
    \frac{(J-\xi(\lambda+1))^2}{4\lambda} \qquad \text{occurring at} \quad b=\frac{J+\xi(\lambda-1)}{2\lambda}.
\end{equation}
In particular, assuming $J \geqslant 0$ and the constraints $a \leqslant 2$, $b \leqslant 2$, and $a+b \leqslant 2$, an upper bound is obtained by taking $J=0$, as the constraints imply $J \leqslant 2 \max\{ u/v,1 \}$ and this upper bound for $J$ yields a smaller value since $\max\{u/v,1\} < \xi(\lambda+1)$.

\begin{lem}
Assume $H$ satisfies (\ref{eq:conditions-on-H}).
Given the constraints $J \geqslant \sum_i \tfrac{d_i}{v} e_i$, $a \leqslant 2$, $b \leqslant 2$, and $a+b \leqslant 2$, we have
\[
    (a-\xi - \frac{e}{2}b)(b-\xi) - \frac{1}{2} \sum_{i=1}^t e_i(e_i-\xi) \leqslant \frac{\xi^2(\lambda+1)^2}{4\lambda} + \frac{t \xi^2}{8}.
\]
\end{lem}

\begin{proof}
For $J \geqslant 0$, the result follows by bounding the first term on the left side by setting $J=0$ in (\ref{eq:parabola-max}), as well as the fact that $-\frac{1}{2}\sum_i e_i(e_i-\xi) \leqslant t\xi^2/8$
as the maximum value of $-e_i(e_i-\xi)$ is $\xi^2/4$.

Now suppose that $J < 0$. Then $\sum_i \tfrac{d_i}{v} e_i \leqslant J$ is also negative. By (\ref{eq:parabola-max}),
\[
    (a-\xi - \frac{e}{2}b)(b-\xi) \leqslant \frac{(J-\xi(\lambda+1))^2}{4\lambda} = \frac{\xi^2(\lambda+1)^2}{4\lambda} + \frac{\xi(\lambda+1) (-J)}{2\lambda} + \frac{J^2}{4\lambda}.
\]
Using the inequality $-J \leqslant \sum \tfrac{d_i}{v}(-e_i) \leqslant -\sum_{i \in I} \tfrac{d_i}{v} e_i$, where $I \subset \{1,\dots,t\}$ is the subset of indices for which $e_i$ is negative, as well as the Cauchy-Schwarz inequality to deduce $(\sum_{i \in I} \frac{d_i}{v} e_i)^2 \leqslant \#I \sum_{i \in I} (\frac{d_i}{v})^2 e_i^2$, the two terms involving $J$ are bounded above by
\begin{equation}\label{eq:terms-involving-J}
     -\frac{\xi(\lambda+1)}{2\lambda} \sum_{i \in I} \frac{d_i}{v} e_i + \frac{\#I}{4\lambda} \sum_{i \in I} \left(\frac{d_i}{v}\right)^2 e_i^2.
\end{equation}
Bounding each $-e_i(e_i-\xi)$ where $e_i \geqslant 0$ by its maximum value $\xi^2/4$, we see that $-\frac{1}{2} \sum_i e_i(e_i-\xi)$ is bounded above by
\begin{equation}\label{eq:exceptional-term}
    \frac{\xi}{2} \sum_{i \in I} e_i -\frac{1}{2} \sum_{i \in I} e_i^2 + \frac{(t-\#I) \xi^2}{8}.
\end{equation}
As $e_i < 0$ for $i \in I$, the conditions on $d_i/v$ in (\ref{eq:conditions-on-H}) imply that the sum of (\ref{eq:terms-involving-J}) and (\ref{eq:exceptional-term}) is bounded by $(t-\#I)\xi^2/8$, which completes the proof.
\end{proof}

By a similar argument, we can show the following statement whose proof will be sketched.
The statement is similar to \cite[Corollary 15.4.6]{LeP97}, but the proof is more involved since the Picard group is not as simple as $\mathbb{Z}$.
\begin{lem}\label{lem:strictly-ss}
Suppose $H$ is an
ample divisor satisfying (\ref{eq:H-pos-lin-combo}) and (\ref{eq:conditions-on-H}). %such that $H\cdot K<0$.
    Consider a complete family of $\{F_y\}_{y\in Y}$ of semistable sheaves of a fixed class on $X$, parametrized by a smooth algebraic variety $Y$. When the discriminant is large, say as in (\ref{eq:discriminant-bound}), the set of points $y\in Y$ such that $F_y$ is strictly semistable forms a closed subset of codimension $\geqslant 2$.
\end{lem}

\begin{proof}
    For $y\in Y$ such that $F_y$ is strictly semistable, consider one of its Jordan-H\"older filtrations and let $\gr_i$ be the corresponding sub-quotients and $\nu_i=\nu(\gr_i)$.
    The proof of Lemma~\ref{large-hn} also shows that the set
    \begin{equation*}
\left\{
y\in Y \middle\vert \begin{array}{l}
     F_y \mbox{ is strictly semistable and  }\\
    \mbox{$(\nu_i-\nu_j)\cdot D>2$ for some $i,j$ and some $D \in \mathcal{D} \cup \{C+A\}$}
\end{array}
\right\}\end{equation*}
has codimension $\geqslant 2$.

We next consider $y\in Y$ such that (a) $F_y$ is strictly semistable and (b) $(\nu_i-\nu_j)\cdot D\leqslant 2$, $\forall D\in \mathcal{D} \cup \{C+A\}$.
Note that $\Ext^2(\gr_i,\gr_j)\cong \Hom(\gr_j,\gr_i\otimes K_X)^\vee=0$, for all $i,j$.
Thus, $\Ext^2_-(F_y,F_y)=0$, which is calculated with respect to the fixed Jordan-H\"older filtration.
Consider the relative flag scheme $\operatorname{Flag}$ of filtrations of the same type as the Jordan-H\"older-filtration. We have
\[\operatorname{T}_p\operatorname{Flag}\to \operatorname{T}_yY\twoheadrightarrow{}\Ext_+^1(F_y,F_y).\]
Moreover, the vanishing of $\Ext^2(\gr_i,\gr_j)$ implies $\Ext^2_+(F_y,F_y)=0$.
The codimension of set of $y\in Y$ satifying conditions (a) and (b) is bounded below by
    $\ext_+^1(F_s,F_s)\geqslant \sum_{i<j}-\chi(\gr_i,\gr_j)\geqslant 2$.
\end{proof}
Using this lemma, we can immediately strengthen Proposition~\ref{prop:bad-locus-res} replacing ``semistable" by ``stable".

\subsection{Locus of semistable sheaves not admitting Gaeta resolutions}
This subsection is devoted to the proof of Theorem~\ref{thm:moduli-nongr-codim2}.

First, the subset $Z$ is indeed closed.
As in the construction of the moduli space using geometric invariant theory \cite{MumFogKir94}, let $\Quot(\oo_X(-m)^{\oplus N}, f)$ be the Quot scheme such that $M(f)$ is a good quotient of the semistable locus $\Quot^{\sstable}(\oo_X(-m)^{\oplus N}, f)$ with respect to the action by $\GL(N,\mathbb{C})$. According to Proposition~\ref{prop:special-GTR-criterion} and upper semicontinuity, the subset of quotient sheaves that do not admit Gaeta resolutions is closed and invariant under the action of $\GL(N,\mathbb{C})$. Under the good quotient map, the image of this subset is closed and is exactly $Z$.
\footnote{If the Jordan-H\"older grading of a semistable sheaf admits a Gaeta resolution, the sheaf does as well.}

Let
\[G=\prod_{i=1}^n\GL(a_i,\mathbb{C})
,\]
which acts on $\rcone$ (\cite[\S~4.3]{Ped21}).
Let $\bar{G}=G/\mathbb{C}^*(\id,\dots,\id)$. There is an induced action of $\bar{G}$ on $\rcone$.
The universal cokernel (\ref{eq:univ-res}) induces a map
\begin{align}\label{eq:don-f}
    \lambda_{\mathbb{F}}\colon \K(X)\to \Pic^{{G}}(\rcone), \quad w\mapsto \det\left(p_!\left([\mathbb{F}]\cdot q^*w\right)\right).
\end{align}
which will be shown to be an isomorphism.

Since $\rcone$ is an open subset in a vector space and its complement has codimension $\geqslant 2$, $\Pic^{G}(\rcone)$ is isomorphic to the character group $\Char(G)\cong \mathbb{Z}^n$ and $\Pic^{\Bar{G}}(\rcone)\cong \Char(\Bar{G})\cong \mathbb{Z}^{n-1}$.
The character groups can be explicitly described as follows:
\begin{align*}
    \mathbb{Z}^n&\xrightarrow{\cong} \Char(G),\\
    (x_1,\dots,x_n)&\mapsto [(M_1,\dots, M_n)\mapsto \prod _{i=1}^n\det(M_i)^{x_i}],
\end{align*}
and under this isomorphism,
\begin{align}\label{eq:char-quotient-gp}
    \Char(\Bar{G})=\left\{(x_1,\dots,x_n)\in \mathbb{Z}^n \,\middle\vert \,\sum_{i=1}^n a_i x_i=0\right\}.
\end{align}
\begin{prop}
    The map $\lambda_{\mathbb{F}}$ in (\ref{eq:don-f}) is an isomorphism and it induces an isomorphism $f^\perp \xrightarrow{\cong} \Pic^{\Bar{G}}(\rcone)$.
\end{prop}
This is similar to \cite[Lemma 18.5.1]{LeP97} and \cite[Proposition 4.2]{Ped21}.
\begin{proof}
    We calculate $\lambda_{\mathbb{F}}$ using the isomorphism $\Pic^{G}(\rcone)\cong \Char(G)\cong \mathbb{Z}^n$.
    In $\K(X)$, for $j=1,\dots,n$, let $\mathbf{e}_j$ be the class $[\E_j^\vee]$ of the dual bundle.
    Then the $\mathbf{e}_j$ form a basis of $\K(X)$.
    Let $V_j$ be a complex vector space of dimension $a_j$ for $j=1,\dots,n$.
    We can calculate $\lambda_{\mathbb{F}}(\mathbf{e}_j)$ using the $G$-equivariant short exact sequence (\ref{eq:univ-res}):
    \begin{align*}
        \lambda_{\mathbb{F}}(\mathbf{e}_j)
        =\det\left(-\sum_{j\leqslant i\leqslant d}[\oo_\rcone \otimes V_i\otimes \Hom(\E_j,\E_i)]+\sum_{\max\{j,d+1\}\leqslant i\leqslant n}[\oo_\rcone\otimes V_i\otimes \Hom(\E_j,\E_i)]\right).
        \end{align*}
    The map $\lambda_{\mathbb{F}}$ takes the following matrix form:
    \[
    \left[\operatorname{sgn}\left(i-d-\frac{1}{2}\right)\chi(\E_j,\E_i)\right]_{1\leqslant i,j\leqslant n}.
    \]
    Since the matrix is lower triangular with $\pm 1$ on the diagonal, $\lambda_{\mathbb{F}}$ is an isomorphism.
    For $w=\sum_i w_i\mathbf{e}_i\in \K(X)$, $w\in f^\perp$ if and only if $\sum_i w_i\chi (\E_i,f)=0$, if and only if $\lambda_\mathbb{F}(w)\in \Char(\Bar{G})$.
\end{proof}

Let $\rcone^{\sstable}\subset \rcone$ denote the subset of cokernels which are semistable.
According to Proposition~\ref{prop:bad-locus-res}, $\rcone \setminus \rcone^{\sstable}$ has codimension $\geqslant 2$ in $\rcone$.
Let $\rcone^{\stable}\subset \rcone^{\sstable}$ denote the subset of cokernels which are stable.
The coarse moduli property provides a map \[\pi\colon \rcone^{\stable} \to M(f),\]
which factors through $U=M(f)\setminus Z$.
According to Lemma~\ref{lem:strictly-ss}, the restriction map induces isomorphisms $\Pic^{{G}}(\rcone)\cong \Pic^{{G}}(\rcone^{\stable})$ and $\Pic^{\Bar{G}}(\rcone)\cong \Pic^{\Bar{G}}(\rcone^{\stable})$.
The Donaldson map $\lambda_M\colon f^\perp \to \Pic (M(f))$ is an isomorphism according to Theorem~\ref{thm:yoshioka}, which will be proved at the end of the subsection. By the functoriality of the determinant line bundle construction, we have the following commutative diagram:
\begin{equation*}
    \begin{tikzcd}
        f^\perp \ar[rr,"\lambda_{\mathbb{F}}","\cong"below] \ar[d,"\lambda_M"left,"\cong"right] & & \Pic^{\Bar{G}}(\rcone)\ar[d,"\cong"]\\
        \Pic(M(f)) \ar[r,two heads] & \Pic(U) \ar[r] & \Pic^{\Bar{G}}(\rcone^{\stable}).
    \end{tikzcd}
\end{equation*}
It is clear from the diagram that the restriction map $\Pic(M(f))\to \Pic(U)$ is also injective. Since the polarization is general, $M(f)$ is locally factorial (\cite[Corollary 3.4]{Yos96}). Therefore, $Z$ has codimension $\geqslant 2$ in $M(f)$. We have proven Theorem~\ref{thm:moduli-nongr-codim2}.

\begin{proof}[Proof of Corollary \ref{cor:gen-stable-prop}]
The result follows directly from Propositions \ref{prop:unirat-large-L},  \ref{prop:special-GTR-criterion}(b),  \ref{prop:cokernel-properties}(b), and Theorem~\ref{thm:moduli-nongr-codim2}.
\end{proof}

We are left to provide
\begin{proof}[Proof of Theorem~\ref{thm:yoshioka}]Let $f=(r,c_1,c_2)\in \K(S)$ be the corresponding class.
According to \cite[Corollary 3.4]{Yos96}, there is a surjective map $\mathbb{Z}\oplus \Pic S\cong f^\perp \twoheadrightarrow \Pic M(f)$. On the other hand, under our assumption on $c_2$, $\Pic M(f)$ contains a subgroup $\mathbb{Z}\oplus \Pic S$, as shown in \cite[Example 8.1.6]{HuyLeh10}.
\end{proof}

\section{Strange duality}\label{sect:sd}

In this section, we review Le Potier's strange duality conjecture for rational surfaces over $\mathbb{C}$ and use our study of Gaeta resolutions to prove Theorem \ref{thm:sd-injective}, which states that the strange morphism is injective in various cases on $\mathbb{P}^2$ and on $X$ an admissible two-step blowup of $\mathbb{F}_e$. The argument is similar to what was shown on $\bp^2$ in \cite{BerGolJoh16}.
We assume Theorem~\ref{thm:finite-quot-scheme}, which will be proved in \S~\ref{sect:finite-quot}.

\subsection{Strange morphism}\label{subsect:strange-mor}
Let $S$ be a smooth projective rational surface over $\mathbb{C}$ and $H$ be an ample divisor.
Let $\sigma$ and $\rho$ denote two classes in the Grothendieck group $\K(S)$.
On $\K(S)$, there is a pairing given by $\chi(\sigma\cdot\rho)$.

Let $M(\sigma)$ and $M(\rho)$ be the moduli spaces of $H$-semistable sheaves of class $\sigma$ and $\rho$ respectively.
For the moment, suppose there are no strictly semistable sheaves, namely $M(\sigma)=M^{\operatorname{s}}(\sigma)$, and there is a universal family $\mathcal{W}$ over $M(\sigma)\times S$. Let $F$ be a sheaf of class $\rho$, then we have a determinant line bundle on $M(\sigma)$,
\begin{equation*}
    \Theta_\rho:=\det\left({p}_{!}\left(\mathcal{W}\stackrel{L}{\otimes} {q}^*F\right)\right)^{*}.
\end{equation*}
Here, ${p}$ and ${q}$ are projections from $M(\sigma)\times S$ to the first and second factors respectively.
The isomorphism class of $\Theta_\rho$ does not depend on $F$ but only on its class in $\K(S)$, so we are justified in using the subscript $\rho$.
Two universal families may differ by a line bundle pulled back from $M(\sigma)$, but if we assume $\chi(\sigma\cdot \rho)=0$, then they will provide isomorphic determinant line bundles on the moduli space.
Thus, from now on, we assume $\chi(\sigma \cdot \rho)=0$, so that $\Theta_\rho$ is independent of the choice of $\mathcal{W}$.

Even if there does not exist a universal family, we can still define $\Theta_\rho$ by carrying out the construction on the Quot scheme coming from the GIT construction where there is a universal family, and then showing that it descends to $M(\sigma)$.
If there are strictly semistable sheaves of class $\sigma$, we need to further require $\rho\in \{1,h,h^2\}^{\perp \perp}$ for $h=[\oo_H]\in \K(S)$; see \cite[(2.9)]{LP92}.
These conditions will be satisfied in our setting.
Similarly, we can construct a determinant line bundle
\[\Theta_\sigma\to M(\rho).\]

Orthogonal classes $\sigma$ and $\rho$ are {\it candidates for strange duality} if  the moduli spaces $M(\sigma)$ and $M(\rho)$ are non-empty, and if the following conditions on pairs $(W,F) \in M(\sigma) \times M(\rho)$ are satisfied:
\begin{enumerate}[(a)]\label{cond:sd-candidates}
\item $h^2(W \otimes F) = 0$ and $\TOR^1(W,F)=\TOR^2(W,F)=0$ for all $(W,F)$ away from a codimension $\geqslant 2$ subset in $M(\sigma)\times M(\rho)$, and
\item $h^0(W \otimes F)=0$ for some $(W,F)$.
\end{enumerate}

Under these conditions, there is a line bundle
\[\Theta_{\sigma,\rho}\to M(\sigma)\times M(\rho)\]
with a canonical section whose zero locus is given by
\[\{\, (W,F) \mid h^0(W \otimes F) > 0 \,\} \subset M(\sigma) \times M(\rho), \]
(see \cite[Proposition 9]{LP05}).
The see-saw theorem implies that
\[\Theta_{\sigma,\rho}\cong \Theta_\rho\boxtimes \Theta_\sigma,\]
see \cite[Lemme 8]{LP05}. Then, using the K\"{u}nneth formula, the canonical section of $\Theta_{\sigma,\rho}$ induces a linear map
\[\operatorname{SD}_{\sigma,\rho}: H^0(M(\rho),\Theta_{\sigma}))^* \rightarrow H^0\big(M(\sigma),\Theta_\rho)\]
that is well-defined up to a non-zero scalar.
Following Le Potier, we call this the {\em strange morphism}.
\begin{conj}[Le Potier]
If $\mathrm{SD}_{\sigma,\rho}$ is nonzero, then it is an isomorphism.
\end{conj}

We focus on the case when $S$ is $\bp^2$ or $X$, an admissible blowup of $\mathbb{F}_e$. For the former, we take $H$ to be the hyperplane class, while for the latter, we assume $H$ is general and satisfies (\ref{eq:H-pos-lin-combo}, \ref{eq:conditions-on-H}).
Let
\[
    \rho = (1,0,1-\ell)
\]
be the numerical class of an ideal sheaf of $\ell \geqslant 1$ points, so that $M(\rho) = S^{[\ell]}$ is the Hilbert scheme of points on $S$.
Clearly $\rho\in \{1,h,h^2\}^{\perp \perp}$.
Every numerical class orthogonal to $\rho$ has the form
\[
    \sigma = (r,L,\chi=r\ell).
\]
We assume $r$ and $\ell$ are fixed, that $r \geqslant 2$, and that
\[
    \text{$L$ is sufficiently positive}.
\]
The condition on $r$ ensures that general sheaves in $M(\sigma)$ are locally free, and it is no restriction in the study of strange duality as strange duality is known over all surfaces in the case when $\sigma$ and $\rho$ both have rank one. Assumptions about the positivity of $L$ are required to apply many results in this paper, and we assume here that $L$ is sufficiently positive such that the positivity assumptions of every result we need are met. A precise statement of the positivity conditions we impose on $L$ can be found in the appendix.

In particular, since $L$ is sufficiently positive, $\sigma$ admits Gaeta resolutions by Proposition~\ref{prop:chern-has-gr}(a), the discriminant condition (\ref{eq:discriminant-bound}) holds, and $M(\sigma)$ is non-empty by Proposition \ref{prop:bad-locus-res},
so general sheaves in $M(\sigma)$ admit Gaeta resolutions by Proposition~\ref{prop:unirat-large-L} or Theorem \ref{thm:moduli-nongr-codim2}.
In this situation, conditions (a) and (b) on p.~\pageref{cond:sd-candidates} are satisfied. Let $W\in M(\sigma)$ and $I_Z\in M(\rho)$ be such that $\operatorname{sing} W\cap Z=\emptyset$, which holds away from codimension $r+1$. Under this condition $\TOR_i(W,I_Z)=0$ for $i=1,2$.
Furthermore,
\[h^2(W\otimes I_Z)=h^2(W)=\hom(W,K_S)\]
vanishes by semistability. Condition (b) also holds by Proposition~\ref{orth-ker-quot}. Thus, $\sigma$ and $\rho$ are candidates for strange duality, and we will study
the strange morphism
\[\SD_{\sigma,\rho}\colon H^0(S^{[\ell]},\Theta_\sigma)^*\to H^0(M(\sigma),\Theta_\rho).\]
We begin with the determinant line bundle $\Theta_\sigma$ on the Hilbert scheme of points.

\subsection{Determinant line bundles on the Hilbert scheme of points}\label{ss:det-line-bundles}
We first review some general results about line bundles on the Hilbert scheme of points on surfaces.
The Hilbert scheme of points, $S^{[\ell]}$, is a resolution of singularities of the symmetric product $S^{(\ell)}$. The resolution, which we denote by $\pi\colon S^{[\ell]}\to S\syml$, is called the {\em  Hilbert-Chow morphism}.
Given a line bundle $M$ on $S$, let $M^{\boxtimes \ell}$ be $\otimes_{i=1}^\ell \operatorname{pr}_i^*M$ on the $\ell$-fold product $S^\ell=S\times \cdots \times S$. There is an $\frak{S}_\ell$-action on $S^\ell$ such that $M^{\boxtimes \ell}$ is $\frak{S}_\ell$-equivariant. The line bundle $M^{\boxtimes \ell}$ descends onto $S^{(\ell)}$, giving a line bundle $M\syml$.  We denote its pullback to $S\hilbl$ as
\[M_\ell=\pi^* M\syml.\]
Via this construction, we can view $\Pic S$ as a subgroup of $\Pic S\hilbl$, by sending $M$ to $M_\ell$. Fogarty~\cite{Fog73} showed that under this inclusion
\[\Pic S\hilbl\cong \Pic S\oplus \mathbb{Z}\frac{E}{2}.\]
Here, $E$ is the exceptional divisor of the Hilbert-Chow morphism, which parametrizes non-reduced subschemes.
Furthermore, if $M$ is ample, then $M_\ell$ is nef, and the canonical divisor on $S^{[\ell]}$ is $K_{S^{[\ell]}} = (K_S)_{\ell}$.

The determinant line bundle on $S^{[\ell]}$ induced by $\sigma$ is
\[
    \Theta_{\sigma} = L_{\ell} - \frac{r}{2} E.
\]
By the Kodaira vanishing theorem, if $\Theta_\sigma-K_{S^{[\ell]}}=L_{\ell}-(K_{S})_{\ell}-\frac{r}{2}E$
is ample on $S^{[\ell]}$, then $\Theta_\sigma$ has no higher cohomology.
Results of Beltrametti, Sommese,
Catanese, and G\"{o}ttsche \cite{BelSom88,CatGot90} show that if $M$ is a line bundle on $S$, then $M_{\ell} - \frac{1}{2} E$ is nef if $M$ is ($\ell-1$)-very ample and very ample if $M$ is $\ell$-very ample. Thus, we deduce the following:

\begin{lem} \label{lem:high-coh-theta-0} Suppose $L$ is sufficiently positive, for instance $L -K_S$ is the tensor product of an $\ell$-very ample line bundle and $r-1$ ($\ell$-1)-very ample line bundles. Then $\Theta_{\sigma}$ has vanishing higher cohomology.
\end{lem}

Thus, the vanishing of the higher cohomology of $\Theta_\sigma$ follows from the assumption that $L$ is sufficiently positive. A precise statement of a sufficient condition on $L$ is (\ref{eq:last-condition}) in the appendix.

\subsection{Injectivity of the strange morphism}

To prove Theorem \ref{thm:sd-injective}, we will make use of a Quot scheme argument. As above, let $\sigma = (r,L,r\ell)$ and $\rho = (1,0,1-\ell)$ be orthogonal classes, where $\ell \geqslant 1$ and $r \geqslant 2$ are fixed, so sheaves of class $\sigma$ are expected to be locally free, and we assume $L$ is sufficiently positive in the sense explained in the appendix.
Consider the class \[v = \sigma + \rho.\] The first part of the argument is to show that if $V$ is a vector bundle of class $v$ that admits a general Gaeta resolution, then $\Quot(V^*,\rho)$ is finite and reduced, which will be proved in \S~\ref{sect:finite-quot}. The strategy is to show that the relative Quot scheme over the space of Gaeta resolutions has relative dimension 0. This is known for $\bp^2$ by \cite{BerGolJoh16}, so we write the argument for $X$, an admissible blowup of $\mathbb{F}_e$, though it works for $\bp^2$ as well.
Since $\sigma$ admits Gaeta resolutions, an easy calculation using Proposition \ref{prop:chern-has-gr}(a) shows that $v$ also admits Gaeta resolutions, with the same exponents $\gamma_i$ and $\gamma_j$ and with $\alpha_1,\alpha_2,\alpha_3$ each larger by $\ell$ and $\alpha_4$ smaller by $\ell-1$.

The starting point is the following result.

\begin{lem}\label{lem:isolated-quotient}
There is a vector bundle $V$ of class $v$ such that $V$ has a Gaeta resolution and $\Quot(V^*,\rho)$ contains an isolated point.
\end{lem}

\begin{proof} Choose a vector bundle $W$ of class $\sigma$ that admits a Gaeta resolution and a general ideal sheaf $I_Z$ of class $\rho$ such that $H^0(W \otimes I_Z) = 0$, which is possible by Proposition \ref{orth-ker-quot}. Since $h^0(W) = r\ell \geqslant \ell$ and $Z$ is general, we can choose a quotient $W \to \oo_Z$ such that $H^0(W) \to H^0(\oo_Z)$ is surjective. Let $J$ denote the kernel of $W \to \oo_Z$. Since $L$ is sufficiently positive, $\Ext^1(J,\oo)$ is large (for instance, (\ref{eq:f-admits-gr}) suffices). Then a general extension $V$ of $J$ by $\oo$ is locally free (see the proof of \cite[Lemma 5.9]{BerGolJoh16}) and there are short exact sequences
\begin{align*}
    &0 \to J \to W \to \oo_Z \to 0, \\
     &0\to \oo \to V \to J \to 0,
\end{align*}
which together give a long exact sequence $0 \to \oo \to V \to W \to \oo_Z \to 0$ whose dual
\[
    0 \to W^* \to V^* \to I_Z \to 0
\]
is a point of $\Quot(V^*,\rho)$. As the tangent space at this point is $\Hom(W^*,I_Z) \cong H^0(W \otimes I_Z)=0$, this is an isolated point of $\Quot(V^*,\rho)$. Thus, all that remains is to show that $V$ admits a Gaeta resolution.

We do this in two steps using Proposition~\ref{prop:special-GTR-criterion}, by first showing that $J$ admits a Gaeta resolution. First, we see that $H^p(J)=0$ for $p=1,2$ using the corresponding vanishings for $W$ and the fact that $H^0(W) \to H^0(\oo_Z)$ is surjective. Similarly, $H^1(J(E_i))=0$ for all $i\in S_0 \cup S_1$, as the induced map $H^0(W(E_i)) \to H^0(\oo_Z(E_i))$ is still surjective. Finally, the vanishing of $H^p(J(D))$ for $p \ne 1$ and $D$ one of the divisors appearing in Proposition~\ref{prop:special-GTR-criterion} (b.ii) follows from the same vanishing for $W$ and the vanishing of the higher cohomology of $\oo_Z$.

Second, since $H^p(\oo)=0$ for $p=1,2$, $H^p(\oo_{E_i})=0$.
Using the vanishings for $J$, we obtain that $H^p(V)=0$ for $p=1,2$ and that $H^1(V(E_i))=0$ for all $i \in S_0 \cup S_1$.
For each divisor $D$ appearing in Proposition~\ref{prop:special-GTR-criterion} (b.ii), $H^p(\oo(D))=0$ for all $p$, so $H^p(V(D))\cong H^p(J(D))$ and this vanishes for $p \ne 1$.
Therefore, $V$ also admits a Gaeta resolution.
\end{proof}

The lemma establishes that the relative Quot scheme contains a point with vanishing relative tangent space, which can be deformed to give an open set with this property. The main technical task is to prove that the relative Quot scheme cannot have any other components of the same dimension, which is carried out in the next section for $S = X$, and which was done for $\bp^2$ in \cite{BerGolJoh16}. The proof requires strong positivity assumptions on $L$.

The second part of the argument is to count the points of $\Quot(V^*,\rho)$ and compare the result to $h^0(S^{[\ell]},\Theta_\sigma)$. In previous work \cite{GolLin22}, we showed that if the Quot scheme is finite and reduced, then the number of points is
\[
    \# \Quot(V^*,\rho) = \int_{S^{[\ell]}} c_{2\ell}({V}^{[\ell]}),
\]
where ${V}^{[\ell]}$ is the tautological vector bundle defined as in (\ref{eq:taut-bundle}), whose rank is $(r+1)\ell$ and whose fiber at a point $[I_Z] \in S^{[\ell]}$ is $H^0(V \otimes \oo_Z)$. By the following theorem, this top Chern class calculates the Euler characteristic of the determinant line bundle:

\begin{thm}\label{numbers-match} Let $S$ be a smooth projective surface with $\chi(\oo_S)=1$. Then
\[
    \int_{S^{[\ell]}} c_{2\ell}({V}^{[\ell]}) = \chi(S^{[\ell]},\Theta_{\sigma}).
\]
\end{thm}
Over Enriques surfaces, this statement was proved by Marian-Oprea-Pandharipande \cite[Proposition 3.2]{MarOprPan22} . In general, it was conjectured by Johnson and shown to be equivalent to another conjecture (\cite[Conjecture 1.3, Theorem 4.1]{Joh18}). The second conjecture was proven by the works of Marian-Opera-Pandharipande \cite{MarOprPan19,MarOprPan22} and G\"ottsche-Mellit \cite{GotMel22} combined. See for example \cite[Corollary 1.2]{GotMel22}.

Now, using the fact that $L$ sufficiently positive ensures that $\Theta_\sigma$ has vanishing higher cohomology, we can deduce the injectivity of the strange morphism.

\begin{proof}[Proof of Theorem \ref{thm:sd-injective}] Above, we checked the conditions for the strange morphism $\mathrm{SD}_{\sigma,\rho}$ to be well-defined. Let $V$ be a vector bundle of class $\sigma + \rho$ that admits a general Gaeta resolution. By Theorem \ref{thm:finite-quot-scheme}, $\Quot(V^*,\rho)$ is finite and reduced and the points of $\Quot(V^*,\rho)$ are short exact sequences
\[
    0 \to W_i^* \to V^* \to I_{Z_i} \to 0
\]
for sheaves $W_i$ of class $\sigma$ that are locally free and semistable.
Since the Quot scheme is finite and reduced, the tangent space $\Hom(W_i^*,I_{Z_i}) \cong H^0(W_i \otimes I_{Z_i})$ is 0, while if $i \ne j$ then $\Hom(W_i^*,I_{Z_j}) \ne 0$ follows from semistability, as otherwise the induced map $W_i^* \to V^* \to I_{Z_j}$ would be zero, hence $W_i^* \to V^*$ would factor through $W_j^* \to V^*$, yielding an equality $W_i^*=W_j^*$ as subsheaves of $V^*$, identifying the two points of the Quot scheme.
It follows that the hyperplanes in $H^0(M(\rho),\Theta_{\sigma})$ determined by the points $I_{Z_i}$ map under $\mathrm{SD}_{\sigma,\rho}$ to linearly independent lines $\Theta_{Z_i} \in \bp H^0(M(\sigma),\Theta_{\rho})$.
Thus, the rank of $\mathrm{SD}_{\sigma,\rho}$ is at least
\[
    \# \Quot(V^*,\rho) = \int_{S^{[\ell]}} c_{2\ell}({V}^{[\ell]}) = \chi(S^{[\ell]},\Theta_{\sigma}) = h^0(S^{[\ell]},\Theta_{\sigma}),
\]
namely $\mathrm{SD}_{\sigma,\rho}$ is injective.
\end{proof}

\section{Finite Quot schemes}\label{sect:finite-quot}

This section is devoted to the proof of Theorem~\ref{thm:finite-quot-scheme}.
We will show that for $X$ an admissible blowup of $\mathbb{F}_e$ over $k$, under appropriate conditions, Quot schemes are indeed finite and reduced. This is an extension of the corresponding result in \cite{BerGolJoh16} for $\bp^2$. We then apply our previous work \cite{GolLin22} to enumerate the finite Quot scheme.

We first sketch the ideas. For a vector bundle $V$ admitting a Gaeta resolution, instead of directly studying ideal sheaf quotients $V^* \twoheadrightarrow I_Z$, we replace $V^*$ by the dual of the Gaeta resolution of $V$ and $I_Z$ by the canonical map $\oo\twoheadrightarrow \oo_Z$. Namely, we consider commutative diagrams of the form (\ref{quot-comm-diag}), which we can view as a family over an open subset $R$ of the resolution space $R_v$ (\ref{eq:dual-res-sp}). We prove in \S~\ref{subsect:dim-count} that the main component of the family has the same dimension as $R$, hence we deduce that an open subscheme of the main component is isomorphic to a relative Quot scheme over an open subset of $R$.
The fibers of the relative Quot scheme have dimension 0, and generically the relative Zariski tangent space has dimension $0$. We conclude that when $V$ is general with appropriate numerical constraints, the Quot scheme is finite and reduced.

We prove Theorem~\ref{thm:finite-quot-scheme} in \S~\ref{subsect:proof-finite-quot}, except for leaving the technical dimension count argument for \S~\ref{subsect:dim-count}.

\subsection{Finite Quot schemes}\label{subsect:proof-finite-quot}

Over $X$ an admissible blowup of $\mathbb{F}_e$, as in the previous section, let $\rho=(1,0,1-\ell)$ be the class of an ideal sheaf of $\ell$ points for some fixed $\ell \geqslant 1$ and $\sigma = (r,L,\chi=r\ell)$ be an orthogonal class for some fixed $r \geqslant 2$ and $L$ satisfying the positivity conditions in the appendix. In particular, the need for the positivity condition (\ref{eq:first-three-conditions}) will be seen in the proof of Proposition \ref{dim-comm-diag}. Then $v = \sigma + \rho$ admits Gaeta resolutions, with exponents denoted $\alpha_1,\{ \gamma_j \}_{j \in S_1}, \alpha_2,\alpha_3,\{ \gamma_i \}_{i \in S_0}, \alpha_4$ as in (\ref{gtr-blowup}). Since $L$ is sufficiently positive, all these exponents are large except for $\alpha_4 = (r-1)\ell + 1$.

\subsubsection{Morphisms vs. chain maps}
Consider a vector bundle $V$ of class $v$ that admits a general Gaeta resolution. Then $V$ is locally free, so its dual has a resolution
\begin{equation}\label{eq:dual-gr}
    0 \to V^* \to \Lambda \xrightarrow{\varphi} \Omega \to 0,
\end{equation}
where
\begin{align*}\Lambda &= \oo^{\alpha_4} \oplus \bigoplus_{i \in S_0} \oo(E_i)^{\gamma_i} \oplus \oo(A)^{\alpha_3} \oplus \oo(C)^{\alpha_2}, \mbox{ and }\\
\Omega &= \bigoplus_{j \in S_1} \oo(C+A-E_j)^{\gamma_j} \oplus \oo(C+A)^{\alpha_1}.
\end{align*}
We wish to study quotients of $V^*$ of class $\rho$, which are expected to be ideal sheaves of points. Instead of considering maps $V^* \twoheadrightarrow I_Z$, we replace them by chain maps of complexes using (\ref{eq:dual-gr}) and the short exact sequence $0 \to I_Z \to \oo \xrightarrow{1_Z} \oo_Z \to 0$. Namely, we consider commutative diagrams
\begin{equation}\label{quot-comm-diag}
\xymatrix{
    \Lambda \ar[r]^\varphi \ar[d]_\pi & \Omega \ar[d]^\psi \\
    \oo \ar[r]^{1_Z} & \oo_Z
}\end{equation}
where $\varphi$ is surjective and $Z$ has length $\ell$.
Letting $R \subset R_v$ denote the open subset of the space of Gaeta resolutions for which the cokernel $V$ is locally free, dualization yields an inclusion
\begin{equation}\label{eq:dual-res-sp}
    R \subset \mathbb{P}\Hom(\Lambda,\Omega)=:\mathbb{P}
\end{equation}
and a universal sequence
\begin{equation}\label{univ-dual-gtr}
    0\to \mathcal{V}^*\to \pi_X^*\Lambda\to \pi_R^*\mathcal{O}_{R}(1) \otimes \pi_X^*\Omega \to 0,
\end{equation}
where $\pi_R$ and $\pi_X$ denote the projections from $R\times X$ to the factors.

Let $\Omega^{*[\ell]}$ be the Fourier-Mukai transform of $\Omega^*$ over $X^{[\ell]}$ and
\[
    \Xi = \{ (\varphi,\psi,\pi) \mid \psi \circ \varphi = 1_Z \circ \pi \mbox{ up to a non-zero scalar}\} \subset \bp \times \mathbb{P}(\Omega^{*[\ell]}) \times \bp^{\alpha_4-1}
\]
be the locus of diagrams (\ref{quot-comm-diag}), which are {\em commutative}, with the reduced induced scheme structure.
We will show in Proposition~\ref{dim-comm-diag} that \[\dim \Xi=\dim \mathbb{P}.\]

\subsubsection{Relative Quot schemes}
On the other hand, we consider the relative Quot scheme $\Quot_{\pi_R}:=\Quot_{\pi_R}(\mathcal{V}^*, \rho)$. Letting
\[
    U_Q\subset \Quot_{\pi_R}
\]
denote the subset where the relative Zariski tangent space has dimension $0$, $U_Q$ is open (by upper semicontinuity), non-empty (by Lemma \ref{lem:isolated-quotient}), and smooth.

We claim that quotients in $U_Q$ can only be ideal sheaves of points. A sheaf $F$ of class $\rho$ could take the following forms: $F$ is isomorphic to an ideal sheaf $I_Z$, $F$ contains dimension $0$ torsion, or $F$ contains dimension $1$ torsion. The second case cannot occur as it violates our assumption on the relative tangent space.
In the third case, let $T$ denote the torsion subsheaf of $F$, so that $F/T$ is torsion free. Then the quotient $F\twoheadrightarrow F/T$ provides a nonzero morphism $V^*\to \oo(-D)$ where $D$ is a non-trivial effective curve given by $c_1(T)$.
But this cannot happen under the conditions in Proposition~\ref{prop:no-sec-van-curves}, which are guaranteed when $L$ is sufficiently positive, see the appendix.

\subsubsection{Finite Quot schemes}
We next define a regular map \[\iota\colon U_Q \to \Xi\] over $\mathbb{P}$ by associating to each quotient $V^* \twoheadrightarrow I_Z$ a diagram of the form (\ref{quot-comm-diag}).
It is enough to define it on the level of functor of points. Furthermore, it is enough to consider morphisms from an affine scheme. Let $S$ be the spectrum of some $k$-algebra.
A morphism $s\colon S\to \Quot_{\pi_R}$ is equivalent to a family of quotients $s_X^*\mathcal{V}^*\to \mathcal{I}_\mathscr{Z}$ over $S\times X$. Here, $s_X=s\times \id_X$ and $\mathcal{I}_\mathscr{Z}$ is the ideal sheaf of a subscheme $\mathscr{Z}\subset S\times X$. We denote the projection maps on $S \times X$ by $\pi_S$ and $\pi_X$.
Applying the functor $\lHom_{\pi_S}(-, \oo_{S\times X})$ to the pull-back of the sequence (\ref{univ-dual-gtr}) via $s_X$, we have the following exact sequence
\begin{align*}
     0&\to \lHom_{\pi_S}(\pi_X^*\Omega, \oo_{S\times X})\to \lHom_{\pi_S}(\pi_X^*\Lambda, \oo_{S\times X})\to\lHom_{\pi_S}(s_X^*\mathcal{V}^*, \oo_{S\times X}) \\
     &\to \lExt^1_{\pi_S}(\pi_X^*\Omega, \oo_{S\times X}).
\end{align*}
According to our choice of $\Omega$, the first term and the last term are zero. Therefore, the family $s_X^*\mathcal{V}\to \mathcal{I}_\mathscr{Z}$ can be completed to a commutative diagram
\begin{equation*}
    \begin{tikzcd}
    0 \arrow{r} & s_X^*\mathcal{V}^* \arrow{r}\arrow{d} & \pi_X^*\Lambda \arrow{r}\arrow{d} & \pi_X^*\Omega\arrow{r}\arrow{d} &0\\
    0 \arrow{r} & \mathcal{I}_\mathscr{Z} \arrow{r} & \oo_{S\times X} \arrow{r} & \oo_{\mathscr{Z}} \arrow{r} & 0
    \end{tikzcd}
\end{equation*}
Then the square on the right provides a morphism $S\to \Xi$.
Clearly, the square uniquely determines the left-most vertical morphism.
We have obtained an injective morphism $\iota\colon U_Q\to \Xi$ whose image is contained in the unique component with the maximal dimension.

The complement $\Xi\setminus \im \iota$ of the image has dimension $<\dim \mathbb{P}$. Then $U\subset R$, the complement of the image of $\Xi\setminus \im \iota\to R$, is non-empty and open in $R$. %We can shrink $U$, if necessary, such that
For each sheaf $V$ parametrized by $U$, $V^*$ has an isolated quotient.
By the definition of $U$, each fiber of $\Xi$ over $[V] \in U$, which parametrizes all non-zero maps of the form $V^* \to I_Z$, is entirely contained in the image of $\iota$. In particular, the maps have to be surjective. Therefore, $\iota$ induces an isomorphism $U_Q|_U\to \Xi|_U$.

We have thus proved Theorem \ref{thm:finite-quot-scheme}(a) for every sheaf $V$ parametrized by $U$.

\subsubsection{Genericity of kernels and cokernels}

We have seen that Quot schemes for $V$ in the nonempty open set $U \subset R_v$ are finite and reduced and that the quotient sheaves are all ideal sheaves $I_Z$. Moreover, in the proof of Lemma \ref{lem:isolated-quotient}, the isolated point $[V^* \twoheadrightarrow I_Z]$ of $\Quot(V^*,\rho)$ has the property that $Z$ is general, hence we can shrink $U$ further if necessary to ensure all the ideal sheaves arising as quotients are general.
Similarly, by Proposition~\ref{prop:unirat-large-L},
general sheaves in $M(\sigma)$ admit Gaeta resolutions, and since $M(\sigma)$ is nonempty, we may choose $W$ in the proof of Lemma \ref{lem:isolated-quotient} to be semistable, hence the isolated quotient $[V \twoheadrightarrow I_Z]$ has kernel $W^*$, which is also semistable. As semistability is an open condition in families, and as the relative Quot scheme can be viewed as a family of sheaves with invariants $\sigma^*$, there is a non-empty open set in the relative Quot scheme of quotients for which the kernel is semistable. Shrinking $U$ if necessary, we may assume that all kernels of the finite Quot schemes are semistable. This proves Theorem \ref{thm:finite-quot-scheme}(b).

\subsubsection{Length of finite Quot schemes}
According to \cite[Proposition 1.2]{GolLin22}, when the Quot scheme $\Quot(V^*,\rho)$ is finite and reduced, it is isomorphic to the moduli space $S(V^*,\rho)$ of limit stable pairs, which consist of torsion free sheaves $F$ of class $\rho$ together with nonzero morphisms $V^* \to F$. Then the virtual fundamental class of $S(V^*,\rho)$ agrees with the fundamental class, and its degree is as stated in Theorem \ref{thm:finite-quot-scheme}(c), by \cite[Theorem 1.1]{GolLin22}.
We refer the reader to \cite{Lin18} for foundational material on the moduli space of limit stable pairs.

\subsection{Dimension count}\label{subsect:dim-count}
We prove that $\dim \Xi = \dim \bp$ by counting diagrams of the form (\ref{quot-comm-diag}) for fixed $\pi$ (nonzero) and $\psi$ (surjective). For fixed $\pi$ and $\psi$, the locus of $\varphi$ is $(\psi_*)^{-1}(1_Z \circ \pi)$ in the exact sequence
\[
    0 \to \Hom(\Lambda,\ker \psi) \to \Hom(\Lambda,\Omega) \xrightarrow{\psi_*} \Hom(\Lambda,\oo_Z) \to \Ext^1(\Lambda,\ker \psi) \to 0.
\]
This locus, if it is non-empty, is an affine space in $\mathbb{P} = \mathbb{P}\Hom(\Lambda,\Omega)$ that is isomorphic to $\Hom(\Lambda,\ker \psi)$.
Since $\pi$ is fixed, we may assume it is nonzero only on the first $\oo$-summand of $\Lambda$. For the locus to be non-empty, we need $1_Z$ to be in the image of the middle map of
\[
    0 \to \Hom(\oo,\ker \psi) \to \Hom(\oo,\Omega) \to \Hom(\oo,\oo_Z) \to \Ext^1(\oo,\ker \psi) \to 0.
\]
Thus, for the purpose of a dimension count we may assume that the image of $\Hom(\oo,\Omega) \to \Hom(\oo,\oo_Z)$ contains $1_Z$.

The idea is to control the dimension of $\Ext^1(\Lambda,\ker \psi)$. In particular, if $\Ext^1(\Lambda,\ker \psi) = 0$, which is true for general $\psi$ as we show below, then the dimension of $(\psi_*)^{-1}(1_Z \circ \pi)$ is $\hom(\Lambda,\Omega) - \hom(\Lambda,\oo_Z)$. Adding to this the dimension of the parameter spaces for $\pi$ and $\psi$, we get
\begin{equation}\label{eq:dim-comm-diag}
    \hom(\Lambda,\Omega) - \hom(\Lambda,\oo_Z) + (\alpha_4 - 1) + (2\ell + \hom(\Omega,\oo_Z) - 1).
\end{equation}
Using the fact that $\hom(\Lambda,\oo_Z) = \ell \, \mathrm{rk}(\Lambda)$, $\hom(\Omega,\oo_Z) = \ell \, \mathrm{rk}(\Omega)$, $\rk(\Lambda) - \rk(\Omega) = r+1$, and $\alpha_4 = (r-1)\ell + 1$, we see that (\ref{eq:dim-comm-diag}) equals $\hom(\Lambda,\Omega) - 1$, which equals $\dim \bp$. Thus, the main technical difficulty is to control the dimension of $\Ext^1(\Lambda,\ker \psi)$ in the case when $\psi$ is not general.

\begin{prop}\label{dim-comm-diag} Let
\[
    M = \max\{\,m(\ell + r + 1 - m) \mid 1 \leqslant m \leqslant \ell \,\}
\]
and assume
\[
    \alpha_2,\alpha_3 \geqslant M, \quad \gamma_j \geqslant M+\ell \mbox{ for }j \in S_1, \quad \text{and} \quad \gamma_i \geqslant M + \ell+ \sum_{j \colon p_j \succ p_i} \gamma_j\mbox{ for }i \in S_0.
\]
Then the space $\Xi$ of commutative diagrams has the expected dimension, namely $\dim \Xi=\dim \mathbb{P}$. Furthermore, there is only one component with this dimension, while other components have strictly smaller dimension.
% \todo{Y: Should it be "dimensions"? \\ T: I think it is okay as is.}
\end{prop}
\begin{proof}
Choose an open set $U \subset X$ containing $Z$ and trivializations on $U$ of $\oo(C+A)$ and $\oo(C+A-E_j)$ for $j \in S_1$. The components of the map $\psi \colon \Omega \to \oo_Z$ generate subspaces of $\Hom(\oo(C+A),\oo_Z)$ and $\Hom(\oo(C+A-E_j),\oo_Z)$, which we identify with subspaces $T, T_{s+1},\dots, T_{t}%\{ T_j \}_{j \in S_1}
$
of $H^0(\oo_Z)$ using the trivializations on $U$.
For $D = 0,E_1,\dots,E_s,A,$ or $C$,
we calculate bounds on $\ext^1(\oo(D),\ker \psi)$ by considering the dimension of the image of
\[
    \Hom(\oo(D),\Omega) \to \Hom(\oo(D),\oo_Z).
\]
Picking a trivialization of $\oo(D)$ on $U$ (shrinking $U$ if necessary), this image is the image of the multiplication map
\[
    m_{D} \colon \left( H^0(\oo(C+A-D)) \otimes T \right) \oplus \bigoplus_{j \in S_1} \left(H^0(\oo(C+A-E_j-D)) \otimes T_j \right) \to H^0(\oo_Z)
\]
defined by $(\phi \otimes x,\sum \phi_j \otimes x_j) \mapsto \phi|_Z \cdot x + \sum \phi_j|_Z \cdot x_j$.
As $|C+A-D|$ is basepoint-free, there is a section $\phi_D \in H^0(\oo(C+A-D))$ that is nowhere zero on the support of $Z$, hence multiplying by $\phi_D$ is injective, which implies that the image of $H^0(\oo(C+A-D)) \otimes T$ has dimension $\geqslant \dim T$. Similarly, as $|C+A-E_i-E_j|$ is basepoint-free when $p_j \succ p_i$ (Example~\ref{ex:LS_ample}),
we can choose $\phi_j \in H^0(\oo(C+A-E_i-E_j))$ that is nowhere zero on the support of $Z$, so the image of $H^0(\oo(C+A-E_i-E_j)) \otimes T_j$ has dimension $\geqslant \dim T_j$.

Define the stratum
\[
    W_{\lambda,\{ \lambda_j \}_{j \in S_1}, \{ \lambda_i \}_{i \in S_0}} \subset \Hom(\Omega,\oo_Z)
\]
by the conditions that $\dim T = \ell-\lambda$, $\dim T_j = \ell - \lambda_j$ for $j \in S_1$, and the sum $\phi_{E_i}|_Z \cdot T + \sum_{j \colon p_j \succ p_i} \phi_j|_Z \cdot T_j$ has dimension $\ell-\lambda_i$ for  $i \in S_0$. By the above, the stratum is empty unless $0 \leqslant \lambda_i \leqslant \lambda,\lambda_j \leqslant \ell$ when $p_j \succ p_i$. For $\psi$ in this stratum, the image of $m_D$ has dimension $\geqslant \ell-\lambda$ and the image of $m_{E_i}$ has dimension $\geqslant \ell-\lambda_i$.
Thus, $\ext^1(\oo(D),\ker \psi) \leqslant \lambda$ and $\ext^1(\oo(E_i),\ker \psi) \leqslant \lambda_i$, which yields the estimate
\[
    \ext^1(\Lambda,\ker \psi) \leqslant \lambda(\alpha_2 + \alpha_3 + \alpha_4) + \sum_{i \in S_0} \lambda_i \gamma_i = \lambda\Big(\alpha_1 + \sum_{j \in S_1} \gamma_j + r + 1\Big) - \sum_{i \in S_0}(\lambda-\lambda_i)\gamma_i.
\]
For strata with $\lambda=0$, then also $\lambda_i=0$, so $\ext^1(\Lambda,\ker \psi)=0$ and the dimension of the space of commutative diagrams for $\psi$ in the union of such strata is equal to the expected dimension. Thus, it suffices to show that for strata with $\lambda > 0$, the codimension of the stratum in $\Hom(\Omega,\oo_Z)$ is at least as large as $\ext^1(\Lambda,\ker \psi)$.

When $\lambda > 0$, we will compare this estimate for $\ext^1(\Lambda,\ker \psi)$ to the codimension of the stratum in $\Hom(\Omega,\oo_Z)$. A map $\psi$ in this stratum can be obtained by first choosing a subspace $T_i \subset H^0(\oo_Z)$ of dimension $\ell-\lambda_i$, then choosing subspaces $T \subset (\phi_{E_i}|_Z)^{-1} \cdot T_i$ and $T_j \subset (\phi_j|_Z)^{-1} \cdot T_i$ of dimension $\ell-\lambda$ and $\ell-\lambda_j$, and finally using $T,T_j$ to define the map $\psi$. Comparing a dimension count based on this description to $\hom(\Omega,\oo_Z)$, we see that the stratum has codimension
\[
    \lambda(\alpha_1+\lambda-\ell) + \sum_{j \in S_1} \lambda_j (\gamma_j + \lambda_j - \ell) + \sum_{i \in S_0} \lambda_i \left( \lambda_i - \lambda + \sum_{j \colon p_j \succ p_i} (\ell-\lambda_j) \right).
\]

The difference between the upper bound for $\ext^1(\Lambda,\ker \psi)$ and this codimension is

\[
     \Delta = \lambda(\ell + r+1-\lambda) + \sum_{i \in S_0} \sum_{j \colon p_j \succ p_i} (\lambda_j-\lambda_i)(\ell-\lambda_j) - \sum_{j \in S_1} (\lambda_j-\lambda) \gamma_j - \sum_{i \in S_0}(\lambda - \lambda_i)(\gamma_i-\lambda_i).
\]
Since $\gamma_i \geqslant M+\ell + \sum_{j \colon p_j \succ p_i} \gamma_j$, we obtain
\[
    \Delta \leqslant \lambda(\ell + r+1-\lambda) - \sum_{i \in S_0} \sum_{j \colon p_j \succ p_i} (\lambda_j-\lambda_i)(\gamma_j+\lambda_j-\ell) - \sum_{i \in S_0} (\lambda-\lambda_i)(M+\ell-\lambda_i).
\]
As $\alpha_2,\alpha_3,\gamma_i \geqslant M \geqslant \lambda(\ell+r+1-\lambda)$, we may assume that the image of $m_D$ for $D=A,C$ is exactly $\phi_D|_Z \cdot T$ and that the image of $m_{E_i}$ has dimension exactly $\ell-\lambda_i$, as otherwise we can improve the upper bound for $\ext^1(\Lambda,\ker\psi)$ by subtracting $M$, leading to a new $\Delta$ that is non-positive. Similarly, as $M+\ell-\lambda_i$ and $\gamma_j+\lambda_j-\ell$ are each $\geqslant M$, $\Delta$ is non-positive unless $\lambda_i=\lambda_j=\lambda$. Thus, all that remains is the case when $\phi_{E_i}|_Z \cdot T = \phi_j|_Z \cdot T_j$ for all $i,j$ such that $p_j \succ p_i$ and the image of each $m_D$ has dimension $\ell-\lambda$. In other words, up to multiplication by units, $T$ is stable under multiplication by rational functions $\zeta$ in $H^0(\oo(C+A-D))$ for $D = A,C,E_i$ and $H^0(\oo(C+A-E_i-E_j)$ for $p_j \succ p_i$.

This final case cannot happen. Since $1_Z$ is in the image of $\Hom(\oo,\Omega) \to \Hom(\oo,\oo_Z)$, $T$ must contain an element $\alpha$ that restricts to a unit at each point of the support of $Z$. But the linear systems $|A|$, $|C|$, $|C+A-E_i|$ for $i \in S_0$, and $|C+A-E_i-E_j|$ for $p_j \succ p_i$ collectively separate points and tangents on $X$ (compare with the proof of Proposition \ref{prop:very-ample}), hence multiplying $\alpha$ by the functions $\zeta$ generates all of $H^0(\oo_Z)$, so the only way $T$ can be stable under multiplication by all $\zeta$ is if $T = H^0(\oo_Z)$. This cannot happen as $\lambda > 0$.

We finish the proof by arguing that there is a unique component with the maximal dimension. As shown above, over a general point $(\psi,\pi)\in \mathbb{P}\left(\Omega^{*[\ell]} \right) \times \mathbb{P}^{\alpha_4-1}$, the fiber of $\Xi$ is an affine space in $\bp$ isomorphic to $\Hom(\Lambda,\ker \psi)$. The fiber of any component with maximal dimension must contain a non-empty open set in this affine space for dimension reasons, but since any two non-empty open sets in an affine space must intersect, there can only be one such component.
\end{proof}

\section*{Appendix: Positivity Conditions}

% format the equation environment
\renewcommand{\theequation}{A.\arabic{equation}}

% reset the counter
\setcounter{equation}{0}

% fix for hyperlink issues after resetting counter
\renewcommand\theHequation{A.\arabic{equation}}

We rephrase the conditions on the exponents of Gaeta resolutions in Proposition \ref{prop:chern-has-gr}(a) to explain the connection with the positivity of the first Chern class. We then summarize the inequalities required for the proof of Theorems \ref{thm:sd-injective} and \ref{thm:finite-quot-scheme}.

Given the numerical class
\[
    f=\Big(r,\alpha A + \delta C - \sum_{i \in S_0} \gamma_i E_i - \sum_{j \in S_1} \gamma_j E_j,\chi \Big)
\]
on $X$ an admissible blowup of $\mathbb{F}_e$, the Euler characteristics in Proposition \ref{prop:chern-has-gr}(a) can be written more explicitly as
\begin{align*}
    \alpha_1 &= \delta + \alpha + r - \chi - \sum_{i \in S_0} \gamma_i - \sum_{j \in S_1} \gamma_j; \\
    \alpha_2 &= \alpha + r - \chi - \sum_{i \in S_0} \gamma_i; \\
    \alpha_3 &= \delta + r - \chi - \sum_{i \in S_0} \gamma_i.
\end{align*}
Thus, rephrasing Proposition \ref{prop:chern-has-gr}(a) in the case when $\chi \geqslant r$ and $\gamma_i \geqslant \sum_{j \colon p_j \succ p_i} \gamma_j$ for all $i$, we can say that $f$ admits Gaeta resolutions if and only if $\gamma_i,\gamma_j,\chi$ are all $\geqslant 0$ and
\begin{equation}\label{eq:f-admits-gr}
    \alpha, \delta \geqslant \sum_{i \in S_0} \gamma_i + \chi - r.
\end{equation}

Now, as in \S~\ref{sect:sd}, consider the orthogonal classes $\rho = (1,0,1-\ell)$ and
\[
    \sigma = \Big(r, L = \alpha A + \delta C - \sum_{i \in S_0} \gamma_i E_i - \sum_{j \in S_1} \gamma_j E_j, \chi = r\ell \Big),
\]
 where $\ell \geqslant 1$ and $r \geqslant 2$ are fixed, and set $v = \sigma + \rho$. The assumption that $L$ is sufficiently positive in \S~\ref{sect:sd} includes three conditions. Let
\[
 M = \max\{\,m(\ell + r + 1 - m) \mid 1 \leqslant m \leqslant \ell \,\}
\]
The first set of conditions used in the proof of Theorem \ref{thm:sd-injective} is:
\begin{align}\label{eq:first-three-conditions}
    \begin{split}
    \gamma_j &\geqslant M \qquad \qquad \quad \quad \: \: \text{for all $j \in S_1$}; \\
    \gamma_i &\geqslant M + \sum_{j \colon p_j \succ p_i} \gamma_j \qquad \: \text{for all $i \in S_0$}; \\
    \alpha,\delta &\geqslant M + \sum_{i \in S_0} \gamma_i + r(\ell - 1).
    \end{split}
\end{align}
These conditions ensure that $\sigma$ and $v$ admit Gaeta resolutions (Proposition \ref{prop:chern-has-gr}(a)), general Gaeta resolutions for $\sigma$ and $v$ are locally free (Proposition \ref{prop:cokernel-properties}), general Gaeta resolutions for $v$ have no sections vanishing on curves (Proposition \ref{prop:no-sec-van-curves}(b)), Weak Brill-Noether holds (Proposition \ref{orth-ker-quot}), and $\dim \Xi = \dim \bp$ (Proposition \ref{dim-comm-diag}).

The inequalities (\ref{eq:first-three-conditions}) imply that $L$ is $M$-very ample, as $L$ can be expressed as a tensor product of $M$ very ample line bundles by Proposition \ref{prop:very-ample} and $m$-very ampleness is additive under tensor products (\cite{HTT05}).

The second condition is the discriminant condition (\ref{eq:discriminant-bound}), which in this context is
\begin{equation}\label{eq:sd-discriminant-bound}
    P\left(\frac{1}{r}L\right) \geqslant \frac{(\lambda+1)^2}{4\lambda} + \frac{t}{8} + \frac{1}{r} + \ell, \qquad \text{where $\lambda = \frac{u}{v} + \frac{e}{2}$}.
\end{equation}
The last condition is that
\begin{equation}\label{eq:fifth-condition}
    \text{$\Theta_\sigma = L_\ell - \frac{r}{2}E$ has vanishing higher cohomology},
\end{equation}
which by Proposition \ref{prop:very-ample} and Lemma \ref{lem:high-coh-theta-0} can be ensured by the following inequalities:
\begin{align}\label{eq:last-condition}
    \begin{split}
    \gamma_j &\geqslant r(\ell-1) \hspace{4cm} \text{for all $j \in S_1$}; \\
    \gamma_i &\geqslant r(\ell-1) + \sum_{j \colon p_j \succ p_i} (\gamma_j+1) \hspace{1.32cm} \text{for all $i \in S_0$}; \\
    \delta &\geqslant r(\ell-1)-1 + \sum_{i \in S_0} (\gamma_i+1); \\
    \alpha &\geqslant r(\ell-1)-1 +e + \sum_{i \in S_0} (\gamma_i+1).
    \end{split}
\end{align}
The sufficiency of these conditions follows from computing $L-K_X$ and observing that it can be decomposed as the tensor product of an $\ell$-very ample line bundle and $r-1$ $(\ell-1)$-very ample line bundles. For $e,s,t$ not too large relative to $\ell$ and $r$, (\ref{eq:last-condition}) is implied by (\ref{eq:first-three-conditions}).

%%%%%%%%%%%%%%%%
% Bibliography %
%%%%%%%%%%%%%%%%

\bibliography{gaetabib}{}
\bibliographystyle{alpha}
\end{document}